\documentclass[11pt,a4paper]{article}
\usepackage{amssymb}
\usepackage{amsmath}
\usepackage{amsthm}
\usepackage{verbatim}
\usepackage{pgf}
\usepackage{listings}
\usepackage{graphicx}
\usepackage{tabularx}
\usepackage{psfrag}
\usepackage{epstopdf}
\usepackage{multirow}
\usepackage{ifthen}
\usepackage{booktabs}
\usepackage{mathtools}
\usepackage{algorithm}
\usepackage{algorithmicx}
\usepackage{algpseudocode}
\usepackage{authblk}
\usepackage[justification=raggedleft,singlelinecheck=false]{caption}
\usepackage{geometry}
\numberwithin{equation}{section}

\newtheorem{thm}{Theorem}[section]  
\newtheorem{lem}[thm]{Lemma}
\newtheorem{cor}[thm]{Corollary}

\theoremstyle{definition}

\newtheorem{asn}[thm]{Assumption}
\newtheorem{prp}[thm]{Proposition}

\newtheorem{rem}[thm]{Remark}

\theoremstyle{remark}
\newtheorem*{prf}{Proof}

\begin{document}

\title{Reconstruction of Refractive Indices from Spectral Measurements of Monodisperse Aerosols}

\author[a]{Tobias Kyrion}

\affil[a]{RWTH Aachen University, Aachen, Germany}

\renewcommand\Authands{ and }

\maketitle

\begin{abstract}
For the investigation of two-component aerosols one needs to know the refractive indices of the two aerosol components. One problem is that they depend on temparature and pressure, so one needs for their determination a robust measurement instrument such as the FASP device, which can cope with rigid environmental conditions. In this article we show that the FASP device is capable of measuring the needed refractive indices, if monodisperse aerosols of the pure components are provided. We determine the particle radii of the monodisperse aerosols needed for this task and investigate how accurate the measurements have to be in order to retrieve refractive indices in a sufficient quality, such that they are suitable for investigations of two-component aerosols.  
\end{abstract}

\noindent \textbf{key words: } nonlinear inverse problem; refractive index reconstruction; Mie theory; nonlinear Tikhonov regularization; nonlinear regression

\vspace{0.2cm}


\section{Setup of the Experiment}
\label{experiment}

This paper provides an algorithm for the reconstruction of refractive indices from spectral measurements of monodisperse aerosols. The experiments we are conducting are similar to the experiments presented in \cite{REDR07} with the difference that we are using air as surrounding medium for the aerosol particles and that temperature and pressure may approach $200^\circ\text{C}$ and $8$ bar respectively. For these rigid conditions reliable databases for refractive indices do not exist up to now. These refractive index databases are needed for the measurement of particle size distributions of polydisperse aerosols using the FASP.

As outlined in \cite{AK16} the FASP measures light intensities $I_{long}(l)$ and $I_{short}(l)$ having passed a long and a short measurement path length $L_{long}$ and $ L_{short}$ respectively. The evaluations of the FASP measurements are based on the relation
\begin{equation}
\int_{0}^{\infty} k(r,l) n(r) dr = e(l) \quad \text{with} \quad e(l) = - \frac{\log (I_{long}(l)) - \log(I_{short}(l))}{L_{long} - L_{short}},
\label{operator_equation}
\end{equation}
where $k(r, l) := \pi r^{2}Q_{ext}(m_{med}(l), m_{part}(l), r, l)$ is the so-called kernel function, $l$ is the wavelength of the incident light, $r$ is the radius of the spherical scattering particle and $m_{med}(l)$ and $m_{part}(l)$ are the refractive indices of the surrounding medium and the particle material depending on the wavelength $l$. The function $Q_{ext}(m_{med}(l), m_{part}(l), r, l)$ is the \textit{Mie extinction efficiency} from \cite{FS01}. The function $n(r)$ is the size distribution of the scattering particels. The right-hand side $e(l)$ in \eqref{operator_equation} is denoted as the \textit{spectral extinction}. 

Now if $n(r)$ is the size distribution of a \textit{monodisperse aerosol}, where all particles possess the same radius $r_{m}$, it is given by $n(r) = n\delta(r - r_{m})$, where $n$ is the total number of particles and $\delta(r - r_{m})$ is a Dirac delta distribution truncated on the positive half-axis. Inserting this into \eqref{operator_equation} gives
\begin{equation}
n \pi r_{m}^{2}Q_{ext}(m_{med}(l), m_{part}(l), r_{m}, l) = e(l), 
\end{equation} 
hence the Mie extinction efficiency is measured directly at the radius $r_{m}$.

The Mie extinction efficiency is given as an infinite series, i.e.
\begin{equation*}
Q_{ext}(m_{med}(l), m_{part}(l), r, l) = \sum_{n=1}^{\infty}q_{n}(m_{med}(l), m_{part}(l), r, l).
\end{equation*}
The computation of the coefficient functions $q_{n}(m_{med}(l), m_{part}(l), r, l)$ will be discussed in Section \ref{Mie_theory}. 

It is clear that in practical computations $Q_{ext}(m_{med}(l), m_{part}(l), r, l)$ can only be approximated by a truncated series, because only the computation of a finite number of the $q_{n}(m_{med}(l), m_{part}(l), r, l)$'s is practically feasible.

We now fix a wavelength $l$. The complex refractive index $m_{part}(l)$ for the wavelength $l$ is reconstructed from FASP measurements of several monodisperse aerosols with particle radii $r_{1}$, ..., $r_{N}$. Let $q(r_{1},l)$, ..., $q(r_{N},l)$ denote the measured spectral extinctions $e(l)$ corresponding to the particle radii $r_{1}$, ..., $r_{N}$. We assume that they are contaminated by additive Gaussian noise, i.e. $q(r_{i},l) = q_{true}(r_{i},l) + \delta_{i}$ with $\delta_{i} \sim \mathcal{N}(0,s_{i}^{2})$ for $i = 1, ..., N$. Furthermore we assume that the standard deviations $s_{i}$ can be estimated from measurements sufficiently accurately, such that we can regard them as known. We have
\begin{equation*}
q_{true}(r_{i},l) = n_{i} \pi r_{i}^{2}\sum_{n=1}^{\infty}q_{n}(m_{med}(l), m_{part}(l), r_{i}, l),
\end{equation*}
where $n_{i}$ is the number of particles having the same radius $r_{i}$. Then a reconstrution of $m_{part}(l)$ is obtained from the set of solutions $M(l)$ of the nonlinear regression problem
\begin{equation}
M(l) := \underset{m \in \mathbb{C}}{\mathrm{argmin}}\; \sum_{i = 1}^{N}\frac{1}{2\left(\frac{s_{i}}{n_{i}}\right)^{2}}\left(\pi r_{i}^{2}\sum_{n = 1}^{N_{tr}}q_{n}(m_{med}(l), m, r_{i}, l) - \frac{q(r_{i},l)}{n_{i}}\right)^{2}.
\label{nonlinregprob}
\end{equation}
Note that $M(l)$ contains in general more than one solution, especially when $q(r_{i},l)$ is perturbed by measurement noise. We discuss nonlinear regression problems with truncated series expansions such as \eqref{nonlinregprob} in Section \ref{regression_with_truncated_series}.

For solving \eqref{nonlinregprob} we use a global optimization strategy presented in Section \ref{recon_alg} to generate reasonable candidates for start values for a local solver for a regularized version of \eqref{nonlinregprob}. Section \ref{smoothest_reconstructions} provides a selection method to find a unique start value out of the candidates. In order to apply a gradient-based local solver we must know the derivatives of the Mie extinction efficiency series, which are discussed in Section \ref{derivatives}.

\section{Mie Theory}
\label{Mie_theory}

We recapitulate Mie theory in an absorbing medium as presented in \cite{FS01}. Our first step is to introduce the complex-valued \textit{Riccati-Bessel-functions} $\xi_n: \mathbb{C} \rightarrow \mathbb{C}$ and $\psi_n: \mathbb{C} \rightarrow \mathbb{C}$ given by
\begin{equation}
\label{funcXi}
\xi_n(z) = \sqrt{\textstyle{\frac{\pi}{2}}} \sqrt{z}J_{n + \frac{1}{2}}(z) 
\end{equation}
and
\begin{equation}
\label{funcPsi}
\psi_n(z) = \sqrt{\textstyle{\frac{\pi}{2}}} \sqrt{z}J_{n + \frac{1}{2}}(z) + i\sqrt{\textstyle{\frac{\pi}{2}}} \sqrt{z}Y_{n + \frac{1}{2}}(z),  
\end{equation}
with the Bessel functions $J_{n + \frac{1}{2}}: \mathbb{C} \rightarrow \mathbb{C}$ and $Y_{n + \frac{1}{2}}: \mathbb{C} \rightarrow \mathbb{C}$ of order $n + \frac{1}{2}$ of first and second kind. We define the \textit{size parameter} $\rho = 2\pi\frac{r}{l}$.
Then we set $z_{med} := \rho \cdot m_{med}$ and $z_{part} := \rho \cdot m_{part}$. Here and in the following we omit the wavelength dependence of $m_{med}$ and $m_{part}$ for better readability. We introduce the notation $m_{med} = n_{med} + ik_{med}$ and $m_{part} = n_{part} + ik_{part}$.

We introduce the so-called \textit{Mie coefficients}:
\begin{equation}
\label{MieCoeff}
\begin{split}
a_n & := \frac{m_{part}\dot{\xi_n}(z_{med})\xi_n(z_{part}) - m_{med}\xi_n(z_{med})\dot{\xi_n}(z_{part})}{m_{part}\dot{\psi_n}(z_{med})\xi_n(z_{part}) - m_{med}\psi_n(z_{med})\dot{\xi_n}(z_{part})} \\
& \\
b_n & := \frac{m_{part}\xi_n(z_{med})\dot{\xi_n}(z_{part}) - m_{med}\dot{\xi_n}(z_{med})\xi_n(z_{part})}{m_{part}\psi_n(z_{med})\dot{\xi_n}(z_{part}) - m_{med}\dot{\psi_n}(z_{med})\xi_n(z_{part})} \\
& \\
c_n & := \frac{m_{part}\psi_n(z_{med})\dot{\xi_n}(z_{med}) - m_{part}\dot{\psi_n}(z_{med})\xi_n(z_{med})}{m_{part}\psi_n(z_{med})\dot{\xi_n}(z_{part}) - m_{med}\dot{\psi_n}(z_{med})\xi_n(z_{part})} \\
& \\
d_n & := \frac{m_{part}\dot{\psi_n}(z_{med})\xi_n(z_{med}) - m_{part}\psi_n(z_{med})\dot{\xi_n}(z_{med})}{m_{part}\dot{\psi_n}(z_{med})\xi_n(z_{part}) - m_{med}\psi_n(z_{med})\dot{\xi_n}(z_{part})} 
\end{split}
\end{equation}

With the Mie coefficients we can express the coefficient functions 
\begin{align}
A_n(\rho, m_{med}, m_{part}) & := \frac{l}{2 \pi m_{part}}\left(\left|c_n\right|^2\xi_n(z_{part})\overline{\dot{\xi_n}(z_{part})} - \left|d_n\right|^2\dot{\xi_n}(z_{part})\overline{\xi_n(z_{part})}\right) \label{CoeffA} \\
& \nonumber \\
& \text{and} \nonumber \\
& \nonumber \\
B_n(\rho, m_{med}, m_{part}) & := \frac{l}{2 \pi m_{med}}\left(\left|a_n\right|^2\dot{\psi_n}(z_{med})\overline{\psi_n(z_{med})} - \left|b_n\right|^2\psi_n(z_{med})\overline{\dot{\psi_n}(z_{med})}\right), \label{CoeffB}
\end{align}
which finally occur in the series expansion
\begin{equation}
\label{ExSeries}
Q_{ext}(r, l, m_{med}, m_{part}) = \frac{l}{2c I(r, l)}\sum_{n = 1}^{\infty}(2n + 1)\mathrm{Im} \big(A_n(\rho, m_{med}, m_{part}) + B_n(\rho, m_{med}, m_{part})\big).
\end{equation}

Here the quantity $I(r, l)$ is the \textit{average incident intensity} of light with wavelength $l$ for a spherical particle with radius $r$ and $c$ denotes the speed of light in vacuum. The function $I(r, l)$ is given by
\begin{equation}
\label{Intensity}
\begin{matrix}
I(r, l) = \displaystyle{\frac{l^2}{8\pi (k_{med})^2} \frac{n_{med}}{2c}}\bigg(1 + \left(4\pi k_{med} \frac{r}{l} - 1\right)e^{\displaystyle{4\pi k_{med} \frac{r}{l}}}\bigg), & \; \mathrm{for} \;\; k_{med} \neq 0 \\
& \\
I(r, l) = \pi r^{2} \displaystyle{\frac{n_{med}}{2c}}, & \; \; \mathrm{for} \;\; k_{med} = 0.
\end{matrix}
\end{equation}

Obviously we cannot evaluate \eqref{ExSeries} exactly, because we cannot compute an infinite sum due to limited processing resources. Therefore we have to truncate this series expansion. In \cite{Wi80} a commonly used truncation index $N_{trunc}$ is presented, which is given by
\begin{equation}
\begin{split}
N_{trunc} & = \big \lceil \left| M \right. + 4.05 \cdot M^{\frac{1}{3}} + \left. 2\right| \big \rceil, \\
\mathrm{with} \;\; M & = \max \lceil \left| \rho \right|, | \rho \cdot m_{med} |, | \rho \cdot m_{part} | \rceil.
\end{split}
\end{equation} 

\section{Derivatives of the Truncated Mie Efficiency Series}
\label{derivatives}

The Bessel functions $J_{\alpha}(z)$ and $Y_{\alpha}(z)$ for an arbitrary weight $\alpha$ fulfill the recurrence relations
\begin{equation}
\frac{d}{dz}\big(z^{\alpha}J_{\alpha}(z)\big) = z^{\alpha}J_{\alpha - 1}(z) \quad \text{and} \quad \frac{d}{dz}\big(z^{\alpha}Y_{\alpha}(z)\big) = z^{\alpha}Y_{\alpha - 1}(z),
\label{diff_recurrence}
\end{equation}
cf. \cite{AS72}. For the Bessel functions occurring in the Riccati-Bessel-functions $\xi_n(z)$ and $\psi_n(z)$ follows from this with the weight $\alpha = n + \frac{1}{2}$ that
\begin{align}
\dot{\xi_n}(z) & = \sqrt{\textstyle{\frac{\pi}{2}}}\sqrt{z}\left(J_{n - \frac{1}{2}}(z) - \frac{n}{z}J_{n + \frac{1}{2}}(z)\right) \label{xi_first_diff}\\
\text{and} \quad \dot{\psi_n}(z) & = \sqrt{\textstyle{\frac{\pi}{2}}}\sqrt{z}\left(J_{n - \frac{1}{2}}(z) - \frac{n}{z}J_{n + \frac{1}{2}}(z)\right) + \sqrt{\textstyle{\frac{\pi}{2}}}\sqrt{z}\left(Y_{n - \frac{1}{2}}(z) - \frac{n}{z}Y_{n + \frac{1}{2}}(z)\right)i \label{first_diff_psi}
\end{align}
for $z \neq 0$. 

We apply \eqref{diff_recurrence} a second time to get $\ddot{\xi_n}(z)$, which yields 
\begin{equation*}
\ddot{\xi_n}(z) =  \sqrt{\textstyle{\frac{\pi}{2}}}\frac{\sqrt{z}}{z^{2}}\left(n(n + 1)J_{n + \frac{1}{2}}(z) + (1 - 2n)J_{n - \frac{1}{2}}(z) + z^{2}J_{n - \frac{3}{2}}(z)\right).
\end{equation*}
For an arbitrary weight $\alpha$ we have the recurrence relation
\begin{equation*}
J_{\alpha - 1}(z) = \frac{2\alpha}{z}J_{\alpha}(z) - J_{\alpha + 1}(z),
\end{equation*}
see \cite{AS72}, and we use it to eliminate the term $J_{n - \frac{3}{2}}(z)$ in the expression for $\ddot{\xi_n}(z)$. Then also $J_{n - \frac{1}{2}}(z)$ cancels out, such that we obtain the representation 
\begin{equation}
\ddot{\xi_n}(z) = \sqrt{\textstyle{\frac{\pi}{2}}}\frac{\sqrt{z}}{z^2}\left(n(n + 1) - z^{2}\right)J_{n + \frac{1}{2}}(z) \label{xi_second_diff}
\end{equation}
only involving $J_{n + \frac{1}{2}}(z)$.

We recapitulate the Cauchy-Riemann equations in its complex form. For a holomorphic function $f: \mathbb{C} \rightarrow \mathbb{C}$ with $f(z) = f(x + iy) = u(x,y) + iv(x,y)$ holds 
\begin{equation}
\dot{f}(z) = \frac{d}{dz}f(z) = \frac{\partial}{\partial x}f(x + iy) = -i\frac{\partial}{\partial y}f(x + iy). \label{CauchyRiemann}
\end{equation}
From this follows
\begin{equation}
u_{x} = \mathrm{Re}\big(\dot{f}(z)\big), \quad u_{y} = -\mathrm{Im}\big(\dot{f}(z)\big), \quad v_{x} = \mathrm{Im}\big(\dot{f}(z)\big) \quad \text{and} \quad v_{y} = \mathrm{Re}\big(\dot{f}(z)\big) \label{partial_derivatives}.
\end{equation}

Using the latter relations we can deduce using $z_{part} = \rho\left(n_{part} + ik_{part}\right)$
\begin{align*}
\frac{\partial}{\partial n_{part}}\mathrm{Re}\big(\xi_{n}(z_{part})\big) & = \rho\mathrm{Re}\big(\dot{\xi_{n}}(z_{part})\big), &
\frac{\partial}{\partial k_{part}}\mathrm{Re}\big(\xi_{n}(z_{part})\big) & = -\rho\mathrm{Im}\big(\dot{\xi_{n}}(z_{part})\big), \\ 
\frac{\partial}{\partial n_{part}}\mathrm{Im}\big(\xi_{n}(z_{part})\big) & = \rho\mathrm{Im}\big(\dot{\xi_{n}}(z_{part})\big), &
\frac{\partial}{\partial k_{part}}\mathrm{Im}\big(\xi_{n}(z_{part})\big) & = \rho\mathrm{Re}\big(\dot{\xi_{n}}(z_{part})\big)
\end{align*}
and analogously
\begin{align*}
\frac{\partial}{\partial n_{part}}\mathrm{Re}\big(\dot{\xi_{n}}(z_{part})\big) & = \rho\mathrm{Re}\big(\ddot{\xi_{n}}(z_{part})\big), &
\frac{\partial}{\partial k_{part}}\mathrm{Re}\big(\dot{\xi_{n}}(z_{part})\big) & = -\rho\mathrm{Im}\big(\ddot{\xi_{n}}(z_{part})\big), \\ 
\frac{\partial}{\partial n_{part}}\mathrm{Im}\big(\dot{\xi}_{n}(z_{part})\big) & = \rho\mathrm{Im}\big(\ddot{\xi_{n}}(z_{part})\big), &
\frac{\partial}{\partial k_{part}}\mathrm{Im}\big(\dot{\xi}_{n}(z_{part})\big) & = \rho\mathrm{Re}\big(\ddot{\xi_{n}}(z_{part})\big).
\end{align*}

The Bessel function values 
\begin{align*}
& J_{0 + \frac{1}{2}}(z_{med}), ..., J_{N_{trunc} + \frac{1}{2}}(z_{med}), \quad Y_{0 + \frac{1}{2}}(z_{med}), ..., Y_{N_{trunc} + \frac{1}{2}}(z_{med}) \\
\text{and} \quad & J_{0 + \frac{1}{2}}(z_{part}), ..., J_{N_{trunc} + \frac{1}{2}}(z_{part}). 
\end{align*}
already computed for a function evaluation of the truncated Mie extinction efficiency can be reused for their derivatives.

Now everything is prepared to differentiate the squared magnitude of the Mie coefficient $a_{n}$ with respect to $n_{part}$ and $k_{part}$. It is sufficient only to discuss $\left|a_{n}\right|^{2}$ in the following, because the differentiation of $\left|b_{n}\right|^{2}$, $\left|c_{n}\right|^{2}$ and $\left|d_{n}\right|^{2}$ works analogously.

First we write the squared norm of the Mie coefficient $a_n$ as $\left|a_{n}\right|^{2} = a_n\overline{a_n}$, which gives
\begin{align*}
\frac{\partial}{\partial n_{part}}\left|a_{n}\right|^{2} & = \left(\frac{\partial}{\partial n_{part}}a_n\right)\overline{a_n} + a_n\left(\frac{\partial}{\partial n_{part}}\overline{a_n}\right) \\
\text{and} \quad \frac{\partial}{\partial k_{part}}\left|a_{n}\right|^{2} & = \left(\frac{\partial}{\partial k_{part}}a_n\right)\overline{a_n} + a_n\left(\frac{\partial}{\partial k_{part}}\overline{a_n}\right).
\end{align*}
We write 
\begin{align*}
a_n = \frac{E}{D} \quad \text{with} \quad E & := m_{part}\dot{\xi_n}(z_{med})\xi_n(z_{part}) - m_{med}\xi_n(z_{med})\dot{\xi_n}(z_{part}) \\
\text{and} \quad D & := m_{part}\dot{\psi_n}(z_{med})\xi_n(z_{part}) - m_{med}\psi_n(z_{med})\dot{\xi_n}(z_{part}),
\end{align*}
which yields
\begin{align*}
\frac{d}{d m_{part}}a_n & = \frac{1}{D^{2}}\left(\left(\frac{d}{d m_{part}}E\right)D - E\left(\frac{d}{d m_{part}}D\right)\right) \\
\text{with} \quad \frac{d}{d m_{part}}E & = \dot{\xi_n}(z_{med})\xi_n(z_{part}) + m_{part}\dot{\xi_n}(z_{med})\rho\dot{\xi_n}(z_{part}) - m_{med}\xi_n(z_{med})\rho\ddot{\xi_n}(z_{part}) \\
\text{and} \quad \frac{d}{d m_{part}}D & = \dot{\psi_n}(z_{med})\xi_n(z_{part}) + m_{part}\dot{\psi_n}(z_{med})\rho\dot{\xi_n}(z_{part}) - m_{med}\psi_n(z_{med})\rho\ddot{\xi_n}(z_{part}).
\end{align*}
Furthermore follow from \eqref{CauchyRiemann} the relations
\begin{align*}
\frac{\partial}{\partial n_{part}}a_n & = \frac{d}{d m_{part}}a_n \\
\text{and} \quad \frac{\partial}{\partial k_{part}}a_n & = \left(\frac{d}{d m_{part}}a_n\right)i.
\end{align*}

Although $\overline{a_n}$ is not holomorphic with respect to $m_{part}$, we can still compute the partial derivatives $\displaystyle{\frac{\partial}{\partial n_{part}}\overline{a_n}}$ and $\displaystyle{\frac{\partial}{\partial k_{part}}\overline{a_n}}$. We obtain using \eqref{partial_derivatives} the relations
\begin{align*}
\frac{\partial}{\partial n_{part}}\overline{a_n} & = \overline{\frac{\partial}{\partial n_{part}}a_n} \\
\text{and} \quad \frac{\partial}{\partial k_{part}}\overline{a_n} & = \overline{\frac{\partial}{\partial k_{part}}a_n}.
\end{align*}
This completes the computations of $\displaystyle{\frac{\partial}{\partial n_{part}}\left|a_{n}\right|^{2}}$ and $\displaystyle{\frac{\partial}{\partial k_{part}}\left|a_{n}\right|^{2}}$.

\vspace*{0.2cm}

For a holomorphic function $f(x + iy)$ we can easily deduce from \eqref{partial_derivatives} that
\begin{align*}
\frac{\partial}{\partial x} \mathrm{Im}\left(f(x + iy)\right) = \mathrm{Im}\left(\frac{\partial}{\partial x} f(x + iy)\right) \quad \text{and} \quad \frac{\partial}{\partial y} \mathrm{Im}\left(f(x + iy)\right) = \mathrm{Im}\left(\frac{\partial}{\partial y} f(x + iy)\right).
\end{align*}
This gives with respect to \eqref{CauchyRiemann}
\begin{align*}
\frac{\partial}{\partial n_{part}}\mathrm{Im}\left(A_{n}\right) & = \mathrm{Im}\left(\frac{\partial}{\partial n_{part}} A_{n}\right) \\
& \\
\text{with} \quad \frac{\partial}{\partial n_{part}} A_{n} & = \frac{l}{2\pi} \left(\frac{\partial}{\partial n_{part}}\left(\left|c_{n}\right|^{2}\right)U_{1} + \left|c_{n}\right|^{2}\frac{\partial}{\partial n_{part}}U_{1} \right. \\
& \hspace*{5cm} \left. - \frac{\partial}{\partial n_{part}}\left(\left|d_{n}\right|^{2}\right)U_{2} - \left|d_{n}\right|^{2}\frac{\partial}{\partial n_{part}}U_{2}\right), \\
& \\
\text{where} \quad U_{1} & = \frac{\xi_n(z_{part})\overline{\dot{\xi_n}(z_{part})}}{m_{part}}, \\
& \\
\frac{\partial}{\partial n_{part}}U_{1} & = \frac{1}{m_{part}^{2}}\left(\left(\rho\dot{\xi_n}(z_{part})\overline{\dot{\xi_n}(z_{part})}+ \xi_n(z_{part})\rho\overline{\ddot{\xi_n}(z_{part})}\right)m_{part} \right. \\
& \hspace*{8cm} \left. - \xi_n(z_{part})\overline{\dot{\xi_n}(z_{part})}\right), \\
& \\
\text{and} \quad U_{2} & = \frac{\dot{\xi}_n(z_{part})\overline{\xi_n(z_{part})}}{m_{part}}, \\
& \\
\frac{\partial}{\partial n_{part}}U_{2} & = \frac{1}{m_{part}^{2}}\left(\left(\rho\ddot{\xi_n}(z_{part})\overline{\xi_n(z_{part})}+ \dot{\xi}_n(z_{part})\rho\overline{\dot{\xi_n}(z_{part})}\right)m_{part} \right. \\
& \hspace*{8cm} \left. - \dot{\xi}_n(z_{part})\overline{\xi_n(z_{part})}\right).
\end{align*}
Analogously we obtain
\begin{align*}
\frac{\partial}{\partial k_{part}}\mathrm{Im}\left(A_{n}\right) & = \mathrm{Im}\left(\frac{\partial}{\partial k_{part}} A_{n}\right) \\
& \\
\text{with} \quad \frac{\partial}{\partial k_{part}} A_{n} & = \frac{l}{2\pi} \left(\frac{\partial}{\partial k_{part}}\left(\left|c_{n}\right|^{2}\right)U_{1} + \left|c_{n}\right|^{2}\frac{\partial}{\partial k_{part}}U_{1} \right. \\
& \hspace*{5cm} \left. - \frac{\partial}{\partial n_{part}}\left(\left|d_{n}\right|^{2}\right)U_{2} - \left|d_{n}\right|^{2}\frac{\partial}{\partial k_{part}}U_{2}\right), \\
& \\
\text{where} \quad \frac{\partial}{\partial k_{part}}U_{1} & = \frac{1}{m_{part}^{2}}\bigg(\left(\rho\left(\dot{\xi_n}(z_{part})\right)i\:\:\overline{\dot{\xi_n}(z_{part})}+ \xi_n(z_{part})\rho\overline{\left(\ddot{\xi_n}(z_{part})\right) i}\right)m_{part}\\
& \hspace*{7.5cm} - \left(\xi_n(z_{part})\overline{\dot{\xi_n}(z_{part})}\right)i\bigg), \\
& \\
\text{and} \quad \frac{\partial}{\partial k_{part}}U_{2} & = \frac{1}{m_{part}^{2}}\bigg(\left(\rho\left(\ddot{\xi_n}(z_{part})\right)i\:\:\overline{\xi_n(z_{part})}+ \dot{\xi}_n(z_{part})\rho\overline{\left(\dot{\xi_n}(z_{part})\right) i}\right)m_{part} \\
& \hspace*{7.5cm} - \left(\dot{\xi}_n(z_{part})\overline{\xi_n(z_{part})}\right)i\bigg).
\end{align*}

The derivatives of $\mathrm{Im}\left(B_{n}\right)$ are much easier to compute, since the dependence on $n_{part}$ and $k_{part}$ lies only in $\left|a_n\right|^{2}$ and $\left|b_n\right|^{2}$ here. So we can deduce
\begin{align*}
\frac{\partial}{\partial n_{part}}\mathrm{Im}\left(B_{n}\right) & = \mathrm{Im}\left(\frac{\partial}{\partial n_{part}} B_{n}\right) \\
& \\
& \\
& \\
& \\
& \\
\text{with} \quad \frac{\partial}{\partial n_{part}} B_{n} & = \frac{l}{2\pi}\left(\frac{\partial}{\partial n_{part}}\left(\left|a_{n}\right|^{2}\right)\frac{\dot{\psi_n}(z_{med})\overline{\psi_n(z_{med})}}{m_{med}} \right. \\
& \hspace*{5cm} \left. - \frac{\partial}{\partial n_{part}}\left(\left|b_{n}\right|^{2}\right)\frac{\psi_n(z_{med})\overline{\dot{\psi_n}(z_{med})}}{m_{med}}\right).
\end{align*}
Analogously we get
\begin{align*}
\frac{\partial}{\partial k_{part}}\mathrm{Im}\left(B_{n}\right) & = \mathrm{Im}\left(\frac{\partial}{\partial k_{part}} B_{n}\right) \\
& \\
\text{with} \quad \frac{\partial}{\partial k_{part}} B_{n} & = \frac{l}{2\pi}\left(\frac{\partial}{\partial k_{part}}\left(\left|a_{n}\right|^{2}\right)\frac{\dot{\psi_n}(z_{med})\overline{\psi_n(z_{med})}}{m_{med}} \right. \\
& \hspace*{5cm} \left. - \frac{\partial}{\partial k_{part}}\left(\left|b_{n}\right|^{2}\right)\frac{\psi_n(z_{med})\overline{\dot{\psi_n}(z_{med})}}{m_{med}}\right).
\end{align*}

\section{Nonlinear Regression using Truncated Series Expansions}
\label{regression_with_truncated_series}

We wish to reconstruct the refractive indices of a particle material from spectral measurements by solving a nonlinear regression problem of the form
\begin{equation}
\label{nonlinreg}
X_{t, \delta} := \underset{\boldsymbol{x} \in \mathbb{R}^{D}}{\mathrm{argmin}} \sum_{i=1}^{N}\frac{1}{2\sigma_{i}^{2}}\left(\sum_{n=1}^{t}a_{n}^{i}(\boldsymbol{x}) - \sum_{n=1}^{\infty}a_{n}^{i}(\boldsymbol{x}_{true}) - \delta_{i}\right)^{2}, 
\end{equation}
where $t \in \mathbb{N}$ is a finite truncation index and $\delta_{i} \sim \mathcal{N}(0, s_{i}^{2})$. Remember that $N$ represents the number of particle radii $r_{i}$ of the different monodisperse aerosols we are investigating. We still assume that for each radius $r_{i}$ the standard deviations $s_{i}$ are determined well enough from a set of experiments, such that they can be regarded as known. We have to confine ourselves to a finite truncation index $t$, because it is practically not feasible to compute all coefficient functions $a_{n}^{i}(\boldsymbol{x})$ for $i = 1, ..., N$. Throughout this paper we assume that the feasible set $\Omega$ is compact.

We define the functions $\boldsymbol{f}_{t}: \mathbb{R}^{D} \rightarrow \mathbb{R}^{N}$ and $\boldsymbol{f}: \mathbb{R}^{D} \rightarrow \mathbb{R}^{N}$ by
\begin{align*}
\boldsymbol{f}_{t}(\boldsymbol{x}) & := \left(\sum_{n=1}^{t}a^{1}_{n}(\boldsymbol{x}), \; ..., \; \sum_{n=1}^{t}a^{N}_{n}(\boldsymbol{x})\right)^{T} \\
\text{and} \quad \boldsymbol{f}(\boldsymbol{x}) & := \left(\sum_{n=1}^{\infty}a^{1}_{n}(\boldsymbol{x}), \; ..., \; \sum_{n=1}^{\infty}a^{N}_{n}(\boldsymbol{x})\right)^{T}. 
\end{align*}
We set $\boldsymbol{e} := \boldsymbol{f}(\boldsymbol{x}_{true}) + \boldsymbol{\delta}$ with $\boldsymbol{\delta} := \left(\delta_{1}, ..., \delta_{N}\right)^{T}$. Then the observed probability density is given by
\begin{equation*}
p_{observed}(\boldsymbol{e} | \boldsymbol{x}) :=  (2\pi)^{-\frac{N_{l}}{2}}\big|\mathrm{det}(\boldsymbol{\Sigma_{\sigma}})\big|^{-\frac{1}{2}}\exp(-\textstyle{\frac{1}{2}}\|\boldsymbol{\Sigma_{\sigma}}^{-\frac{1}{2}}(\boldsymbol{f}_{t}(\boldsymbol{x}) - \boldsymbol{e})\|_{2}^{2})
\end{equation*}
with the covariance matrix $\boldsymbol{\Sigma_{\sigma}} := \mathrm{diag}\left(\sigma_{1}^{2}, ..., \sigma_{N}^{2}\right)$. We know a priori that the vector $\boldsymbol{x}$ specifying our model $\boldsymbol{f}_{t}(\boldsymbol{x})$ lies within the set $\Omega$. This knowledge can be expressed with the prior probability density
\begin{equation*}
p_{prior}(\boldsymbol{x}) := \left(\mathrm{vol}(\Omega)\right)^{-1}I_{\Omega}(\boldsymbol{x}),
\end{equation*} 
where $I_{\Omega}$ is the indicator function of $\Omega$. Now $X_{t, \delta}$ is the set of MAP-estimators of the posterior probability density, i.e.
\begin{equation}
\label{MAP_estimator}
\begin{split}
X_{t, \delta} := \underset{\boldsymbol{x}}{\mathrm{argmax}} & \; p_{posterior}(\boldsymbol{x} | \boldsymbol{e}) \\
 \text{with} \quad p_{posterior}(\boldsymbol{x} | \boldsymbol{e}) \propto p_{observed}(\boldsymbol{e} | \boldsymbol{x})p_{prior}(\boldsymbol{x}) & \propto \exp(-\textstyle{\frac{1}{2}}\|\boldsymbol{\Sigma_{\sigma}}^{-\frac{1}{2}}(\boldsymbol{f}_{t}(\boldsymbol{x}) - \boldsymbol{e})\|_{2}^{2})I_{\Omega}(\boldsymbol{x}).
\end{split}
\end{equation}
We carry out all the following investigations under the next assumption on the covariance matrix:
\begin{asn}
\label{covariance}
\textit{The covariance matrix $\boldsymbol{\Sigma_{\sigma}}$ has the simple form} 
\begin{equation*}
\boldsymbol{\Sigma_{\sigma}} = \delta^{2}\cdot \mathrm{diag}(\sigma_{1}^{2}, ..., \sigma_{N}^{2}) =: \delta^{2} \cdot \boldsymbol{\Sigma},
\end{equation*}
\textit{where $\delta \geq 0$ is an arbitrary but fixed noise level and $\sigma_{1}$, ..., $\sigma_{N}$ are fixed.}
\end{asn}

To simplify notations we introduce the two functions $\boldsymbol{f}_{t}: \mathbb{R}^{D} \rightarrow \mathbb{R}^{N}$ and $\boldsymbol{g}_{t}: \mathbb{R}^{D} \rightarrow \mathbb{R}^{N}$ depending on the truncation index $t$ and defined by
\begin{align*}
\left(\boldsymbol{f}_{t}(\boldsymbol{x})\right)_{i} & := \sum_{n=1}^{\lfloor t \rfloor}a_{n}^{i}(\boldsymbol{x}) + \big(t - \lfloor t \rfloor \big)a_{\lfloor t \rfloor + 1}^{i}(\boldsymbol{x}) \\
\text{and} \quad \left(\boldsymbol{g}_{t}(\boldsymbol{x})\right)_{i} & := \left(\boldsymbol{f}(\boldsymbol{x})\right)_{i} - \left(\boldsymbol{f}_{t}(\boldsymbol{x})\right)_{i}, \quad \text{for} \quad i = 1, ..., N. 
\end{align*}

In the following we will investigate how an element $\boldsymbol{x}_{t, \delta}$ of the set $X_{t, \delta}$ depends on the truncation index $t$. We change to a continuous truncation index here, i.e. we change from now on from \eqref{nonlinreg} to the new regression problem  
\begin{equation}
\begin{split}
\label{nonlinreg_continuous}
X_{t, \delta} & := \underset{\boldsymbol{x} \in \mathbb{R}^{D}}{\mathrm{argmin}} \; F_{t, \delta}(\boldsymbol{x}) \quad \text{s.t.} \quad \boldsymbol{x} \in \Omega,\\
\text{with} \quad F_{t, \delta}(\boldsymbol{x}) & := \|\boldsymbol{\Sigma}^{-\frac{1}{2}}(\boldsymbol{f}_{t}(\boldsymbol{x}) - \boldsymbol{f}(\boldsymbol{x}_{true}) - \boldsymbol{\delta})\|_{2}^{2}
\end{split}
\end{equation}
where the truncation index $t \geq 0$ is allowed to be non-integer. 

As a preparation we prove the following technical lemma, which will form the basis of our continuity and convergence results.
\begin{lem} 
\label{MinimaCurve}
\textit{Let the twice continuously differentiable function $F: \mathbb{R}^{N} \rightarrow \mathbb{R}$ have a strict local minimum $\boldsymbol{x}_{0}$ inside a compact set $S \subset \mathbb{R}^{N}$. Let the function $h: \mathbb{R}^{N} \times \mathbb{R} \rightarrow \mathbb{R}$ have the property $\lim_{\varepsilon \to 0} h(\boldsymbol{x}, \varepsilon) = 0$ for all $\boldsymbol{x} \in S$ and let $ h(\boldsymbol{x}, \varepsilon)$ be twice continuously differentiable with respect to $\boldsymbol{x}$ and continuous in $\varepsilon$. Furthermore we assume that the local minima $\boldsymbol{x}_{\varepsilon}$ of $F_{\varepsilon}(\boldsymbol{x}) := F(\boldsymbol{x}) + h(\boldsymbol{x}, \varepsilon)$ are strict for any $\varepsilon > 0$. Then there exists a sequence of local minima $\boldsymbol{x}_{\varepsilon}$ of $F_{\varepsilon}(\boldsymbol{x})$ with $\lim_{\varepsilon \to 0}\boldsymbol{x}_{\varepsilon} = \boldsymbol{x}_{0}$.}
\end{lem}

\begin{prf}
The strategy of the proof is to construct for given $\varepsilon$ a neighborhood of $\boldsymbol{x}_{0}$ which must contain a local minimizer $\boldsymbol{x}_{\varepsilon}$ of the perturbed function $F_{\varepsilon}(\boldsymbol{x})$. By sending $\varepsilon$ to $0$, this neighborhood shrinks down to the local minimum $\boldsymbol{x}_{0}$ itself, thus yielding the convergence of $\boldsymbol{x}_{\varepsilon}$ to $\boldsymbol{x}_{0}$. To have this neighborhood shrink down to $\boldsymbol{x}_{0}$, it is crucially important that $\boldsymbol{x}_{0}$ must be a strict local minimum.

We define $d(\varepsilon) := \sup_{\boldsymbol{x} \in S}|h(\boldsymbol{x}, \varepsilon)|$. From $\lim_{\varepsilon \to 0} h(\boldsymbol{x}, \varepsilon) = 0$ for all $\boldsymbol{x} \in S$ follows $\lim_{\varepsilon \to 0}d(\varepsilon) = 0$. Let us now introduce the function $F^{-}(\boldsymbol{x}) := F(\boldsymbol{x}) - d(\varepsilon)$. Obviously $\boldsymbol{x}_{0}$ is also a local minimum of $F^{-}(\boldsymbol{x})$, so for $\varepsilon$ sufficiently small there exists a neighborhood $U_{2d(\varepsilon)}(\boldsymbol{x}_{0}) \subset S$ of $\boldsymbol{x}_{0}$ with
\begin{equation*}
F^{-}(\boldsymbol{x}) \geq F^{-}(\boldsymbol{x}_{0}) \quad \text{and} \quad F^{-}(\boldsymbol{x}) - F^{-}(\boldsymbol{x}_{0}) \leq 2d(\varepsilon) \quad \text{for all} \quad \boldsymbol{x} \in U_{2d(\varepsilon)}(\boldsymbol{x}_{0}).
\end{equation*}
In particular we have
\begin{equation*}
\forall \boldsymbol{x} \in \partial U_{2d(\varepsilon)}(\boldsymbol{x}_{0}): \; F^{-}(\boldsymbol{x}) = F^{-}(\boldsymbol{x}_{0}) + 2d(\varepsilon) = F(\boldsymbol{x}_{0}) + d(\varepsilon).
\end{equation*}

Let us assume that there exists an $\boldsymbol{x} \in \partial U_{2d(\varepsilon)}(\boldsymbol{x}_{0})$ with
\begin{equation*}
 F_{\varepsilon}(\boldsymbol{x}) < F(\boldsymbol{x}_{0}) + d(\varepsilon) =   F^{-}(\boldsymbol{x}) =  F(\boldsymbol{x}) - d(\varepsilon).
\end{equation*}
Then $F_{\varepsilon}(\boldsymbol{x}) = F(\boldsymbol{x}) + h(\boldsymbol{x}, \varepsilon)$ implies $- d(\varepsilon) > h(\boldsymbol{x}, \varepsilon)$, hence $- h(\boldsymbol{x}, \varepsilon) > d(\varepsilon) \geq - h(\boldsymbol{x}, \varepsilon)$ by definition of $d(\varepsilon)$, contradiction. Therefore we conclude
\begin{equation}  
\label{border_inequality}
\forall \boldsymbol{x} \in \partial U_{2d(\varepsilon)}(\boldsymbol{x}_{0}): \; F_{\varepsilon}(\boldsymbol{x}) \geq F(\boldsymbol{x}_{0}) + d(\varepsilon). 
\end{equation}

Since $F_{\varepsilon}(\boldsymbol{x})$ is continuous and $\overline{U}_{2d(\varepsilon)}(\boldsymbol{x}_{0})$ is compact for $\varepsilon$ small enough, there exists an $\boldsymbol{x}_{\varepsilon} \in \overline{U}_{2d(\varepsilon)}(\boldsymbol{x}_{0})$ with 
\begin{equation*}
F_{\varepsilon}(\boldsymbol{x}_{\varepsilon}) = \min_{\boldsymbol{x} \in \overline{U}_{2d(\varepsilon)}(\boldsymbol{x}_{0})} F_{\varepsilon}(\boldsymbol{x}).
\end{equation*}
Let us assume $F_{\varepsilon}(\boldsymbol{x}_{\varepsilon}) > F(\boldsymbol{x}_{0}) + d(\varepsilon)$. Then by definition of $\boldsymbol{x}_{\varepsilon}$ we get in particular
\begin{align*}
F(\boldsymbol{x}_{0}) + h(\boldsymbol{x}_{0}, \varepsilon) = F_{\varepsilon}(\boldsymbol{x}_{0}) \geq F_{\varepsilon}(\boldsymbol{x}_{\varepsilon}) >  F(\boldsymbol{x}_{0}) + d(\varepsilon),
\end{align*}
i.e. $h(\boldsymbol{x}_{0}, \varepsilon) > d(\varepsilon) \geq h(\boldsymbol{x}_{0}, \varepsilon)$, contradiction. It follows
\begin{equation}
\label{minimizer_upper_bound}
F_{\varepsilon}(\boldsymbol{x}_{\varepsilon}) \leq F(\boldsymbol{x}_{0}) + d(\varepsilon) \quad \text{and} \quad F_{\varepsilon}(\boldsymbol{x}_{0}) \leq F(\boldsymbol{x}_{0}) + d(\varepsilon),
\end{equation}
where the latter follows with a proof by contradiction as well.

If it happens to hold that $F_{\varepsilon}(\boldsymbol{x}_{\varepsilon}) = F(\boldsymbol{x}_{0}) + d(\varepsilon)$, then we also have $F_{\varepsilon}(\boldsymbol{x}_{0}) = F(\boldsymbol{x}_{0}) + d(\varepsilon)$. Otherwise we have $F_{\varepsilon}(\boldsymbol{x}_{\varepsilon}) < F(\boldsymbol{x}_{0}) + d(\varepsilon)$ and then \eqref{border_inequality} implies that $\boldsymbol{x}_{\varepsilon}$ cannot lie on $\partial U_{2d(\varepsilon)}(\boldsymbol{x}_{0})$, thus it must lie within the interior of $U_{2d(\varepsilon)}(\boldsymbol{x}_{0})$. So in any case \eqref{minimizer_upper_bound} gives that $U_{2d(\varepsilon)}(\boldsymbol{x}_{0})$ must contain a local minimizer $\boldsymbol{x}_{\varepsilon}$ of $F_{\varepsilon}(\boldsymbol{x})$.

Now $\lim_{\varepsilon \to 0} d(\varepsilon) = 0$ gives $\lim_{\varepsilon \to 0}\boldsymbol{x}_{\varepsilon} = \boldsymbol{x}_{0}$. The existence of the last limit is guaranteed by the fact that $\boldsymbol{x}_{0}$ is strict and the claim is proved.

\hfill $\square$
\end{prf}

\begin{prp}
\label{MinimaContinuity}
\textit{Let all coefficient functions $a_{n}^{i}(\boldsymbol{x})$ be twice continuously differentiable and bounded on $\Omega$. We assume that each local minimum $\boldsymbol{x}_{t, \delta}$ of the right hand side function $F_{t, \delta}(\boldsymbol{x})$ in \eqref{nonlinreg_continuous} is strict and lies in the interior of $\Omega$. Then each local minimum depends continuously on the truncation index $t$.}
\end{prp}

\begin{prf}
To prove the claim, one could be tempted to apply the implicit function theorem on the equation $\boldsymbol{d}(t, \boldsymbol{x}_{t,\delta}) = 0$ with $\boldsymbol{d}(s, \boldsymbol{x}) := \nabla \boldsymbol{f}_{s}(\boldsymbol{x})$. This would give that the local minima are parameterized by a function $\boldsymbol{m}(s)$ with the property $\boldsymbol{m}(t) = \boldsymbol{x}_{t,\delta}$, where $s$ is from an environment $U(t)$ of $t$. The problem with this approach is that it requires continuous differentiability of $\boldsymbol{d}(s, \boldsymbol{x})$ in the truncation parameter $s$. Thus the continuous truncation we are using would need more complicated methods such as spline interpolation of the partial sums, which would increase the overall computational effort.

Therefore we use in the following a more direct approach to prove the claim. Let $\varepsilon > 0$ be arbitrary. First we consider an integer truncation index $t \in \mathbb{N}$, i.e. we have $t = \lfloor t \rfloor$. Now for $\varepsilon$ small enough, we get $\lfloor t + \varepsilon \rfloor = t$ and $\lfloor t - \varepsilon \rfloor = t - 1$. This gives
\begin{align*}
\left(\boldsymbol{f}_{t + \varepsilon}(\boldsymbol{x})\right)_{i} & = \left(\boldsymbol{f}_{t}(\boldsymbol{x})\right)_{i} + \varepsilon a_{t + 1}^{i}(\boldsymbol{x}) \\
\text{and} \quad \left(\boldsymbol{f}_{t - \varepsilon}(\boldsymbol{x})\right)_{i} & = \left(\boldsymbol{f}_{t - 1}(\boldsymbol{x})\right)_{i} + (1 - \varepsilon)a_{t}^{i}(\boldsymbol{x})\\
& = \left(\boldsymbol{f}_{t}(\boldsymbol{x})\right)_{i} - \varepsilon a_{t}^{i}(\boldsymbol{x}).
\end{align*}
As next step we turn to an noninteger truncation index $t$. In this case, we can always select $\varepsilon$ small enough such that $\lfloor t + \varepsilon \rfloor = \lfloor t \rfloor$ and $\lfloor t - \varepsilon \rfloor = \lfloor t \rfloor$ respectively hold. This yields
\begin{align*}
\left(\boldsymbol{f}_{t + \varepsilon}(\boldsymbol{x})\right)_{i} & = \left(\boldsymbol{f}_{t}(\boldsymbol{x})\right)_{i} + \varepsilon a_{\lfloor t \rfloor + 1}^{i}(\boldsymbol{x}) \\
\text{and} \quad \left(\boldsymbol{f}_{t - \varepsilon}(\boldsymbol{x})\right)_{i} & = \left(\boldsymbol{f}_{t}(\boldsymbol{x})\right)_{i} - \varepsilon a_{\lfloor t \rfloor + 1}^{i}(\boldsymbol{x}).
\end{align*}
Now we introduce the function
\begin{align*}
\boldsymbol{a}(\boldsymbol{x})
:=
\begin{cases}
\big( a_{t}^{1}(\boldsymbol{x}), ...,  a_{t}^{N}(\boldsymbol{x})\big)^{T}, & \text{for} \quad t - \varepsilon, \; t \in \mathbb{N} \\
\big( a_{\lfloor t \rfloor + 1}^{1}(\boldsymbol{x}), ...,  a_{\lfloor t \rfloor + 1}^{N}(\boldsymbol{x})\big)^{T}, & \text{else}.
\end{cases}
\end{align*}
For $F_{t, \delta}(\boldsymbol{x}) = \|\boldsymbol{\Sigma}^{-\frac{1}{2}}(\boldsymbol{f}_{t}(\boldsymbol{x}) - \boldsymbol{f}(\boldsymbol{x}_{true}) - \boldsymbol{\delta})\|_{2}^{2}$ this yields 
\begin{align*}
F_{t + \varepsilon, \delta}(\boldsymbol{x}) & = F_{t, \delta}(\boldsymbol{x}) + 2\varepsilon \big\langle \boldsymbol{\Sigma}^{-\frac{1}{2}}\boldsymbol{a}(\boldsymbol{x}), \; \boldsymbol{\Sigma}^{-\frac{1}{2}}(\boldsymbol{f}_{t}(\boldsymbol{x}) - \boldsymbol{f}(\boldsymbol{x}_{true}) - \boldsymbol{\delta})\big\rangle + \varepsilon^{2} \|\boldsymbol{\Sigma}^{-\frac{1}{2}}\boldsymbol{a}(\boldsymbol{x})\|_{2}^{2} \\
\text{and} \quad F_{t - \varepsilon, \delta}(\boldsymbol{x}) & = F_{t, \delta}(\boldsymbol{x}) - 2\varepsilon \big\langle \boldsymbol{\Sigma}^{-\frac{1}{2}}\boldsymbol{a}(\boldsymbol{x}), \; \boldsymbol{\Sigma}^{-\frac{1}{2}}(\boldsymbol{f}_{t}(\boldsymbol{x}) - \boldsymbol{f}(\boldsymbol{x}_{true}) - \boldsymbol{\delta})\big\rangle + \varepsilon^{2} \|\boldsymbol{\Sigma}^{-\frac{1}{2}}\boldsymbol{a}(\boldsymbol{x})\|_{2}^{2}.
\end{align*}
Therefore we obtain both for $F_{t + \varepsilon, \delta}(\boldsymbol{x})$ and $F_{t - \varepsilon, \delta}(\boldsymbol{x})$ a decomposition of the form $F_{t + \varepsilon, \delta}(\boldsymbol{x}) = F_{t, \delta}(\boldsymbol{x}) + h^{\varepsilon}_{t, \delta}(\boldsymbol{x})$ and $F_{t - \varepsilon, \delta}(\boldsymbol{x}) = F_{t, \delta}(\boldsymbol{x}) + h^{\varepsilon}_{t, \delta}(\boldsymbol{x})$ respectively, where the function $h^{\varepsilon}_{t, \delta}(\boldsymbol{x})$ is appropriately selected according to above findings. We can readily check $\lim_{\varepsilon \to 0} |h^{\varepsilon}_{t, \delta}(\boldsymbol{x})| = 0$ for all $\boldsymbol{x} \in \Omega$ from the boundedness of the $a_{n}^{i}(\boldsymbol{x})$'s. Then the result follows from Lemma \ref{MinimaCurve}.

\hfill $\square$
\end{prf}

\begin{cor}
\label{continuation}
\textit{Let $t_{1}$ and $t_{2}$ be truncation indices with $t_{1} < t_{2}$. Let $\boldsymbol{x}_{t_{1}, \delta}$ be a local minimizer of \eqref{nonlinreg}. Let $\gamma \in [0, 1]$ and define $t_{\gamma} := t_{1} + \gamma \left(t_{2} - t_{1}\right)$. Then beginning at $\gamma = 0$ one can successively find local minimizers $\boldsymbol{x}_{t_{\gamma}, \delta}$ for the truncation index $t_{\gamma}$ using numerical continuation, see \cite{DH02}. Here for $\gamma_{1} < \gamma_{2}$ the minimizer $\boldsymbol{x}_{t_{\gamma_{1}}, \delta}$ is used as a start vector to compute the next minimizer $\boldsymbol{x}_{t_{\gamma_{2}}, \delta}$. The next parameter $\gamma_{2}$ has to be sufficiently close to $\gamma_{1}$, such that the start vector $\boldsymbol{x}_{t_{\gamma_{1}}, \delta}$ still lies within the domain of convergence for Newton's method.} 
\hfill $\square$
\end{cor}

We use Corollary \ref{continuation} to compute $\boldsymbol{x}_{t, \delta}$ for increasing truncation index $t$ in a stable way. If we would keep $t$ as integer and increase it in integer steps, we might leave the domain of convergence in the continuation method. Therefore we increase them using a smaller step width.

In the following we investigate how well the minimizers $\boldsymbol{x}_{t, \delta}$ of the noise-contaminated regression problem \eqref{nonlinreg_continuous} with truncated series expansions approximate the minimizers $\boldsymbol{x}_{\infty, 0}$ of the noise-free and untruncated problem
\begin{equation}
\label{nonlinreg_untruncated}
\boldsymbol{x}_{\infty, 0} := \underset{\boldsymbol{x} \in \mathbb{R}^{D}}{\mathrm{argmin}} \sum_{i=1}^{N}\frac{1}{2\sigma_{i}^{2}}\left(\sum_{n=1}^{\infty}a_{n}^{i}(\boldsymbol{x}) - \sum_{n=1}^{\infty}a_{n}^{i}(\boldsymbol{x}_{true})\right)^{2} \quad \text{s.t.} \quad \boldsymbol{x} \in \Omega.
\end{equation}

\begin{prp}
\label{minimizer_convergence}
\textit{Let the noise vector $\boldsymbol{\delta}$ fulfill $\lim_{\delta \to 0}\|\boldsymbol{\Sigma}^{-\frac{1}{2}}\boldsymbol{\delta}\|_{2}^{2} = 0$ and let the functions $\boldsymbol{f}(\boldsymbol{x})$ and $\boldsymbol{f}_{t_{\delta}}(\boldsymbol{x})$ be bounded on $\Omega$. Assume $\lim_{\delta \to 0}\|\boldsymbol{\Sigma}^{-\frac{1}{2}}\boldsymbol{g}_{t_{\delta}}(\boldsymbol{x})\|_{2}^{2} = 0$ for all $\boldsymbol{x} \in \Omega$. Then for any strict minimizer $\boldsymbol{x}_{\infty, 0}$ of the right hand side function of \eqref{nonlinreg_untruncated} in the interior of $\Omega$ exist minimizers $\boldsymbol{x}_{t_{\delta}, \delta}$ of \eqref{nonlinreg_continuous} with $\lim_{\delta \to 0} \boldsymbol{x}_{t_{\delta}, \delta} = \boldsymbol{x}_{\infty, 0}$. Here we also assume the $\boldsymbol{x}_{t_{\delta}, \delta}$'s to be strict for all $\delta > 0$.}
\end{prp}

\begin{prf}
With the notation introduced before we can write
\begin{align*}
\boldsymbol{x}_{\infty, 0} \in X_{\infty, 0} & := \underset{\boldsymbol{x} \in \mathbb{R}^{D}}{\mathrm{argmin}} \; F_{\infty, 0}(\boldsymbol{x}) \quad \text{s.t.} \quad \boldsymbol{x} \in \Omega \\
\text{with} \quad F_{\infty, 0}(\boldsymbol{x}) & := \|\boldsymbol{\Sigma}^{-\frac{1}{2}}(\boldsymbol{f}(\boldsymbol{x}) - \boldsymbol{f}(\boldsymbol{x}_{true}))\|_{2}^{2}.
\end{align*}
From the decomposition $\boldsymbol{f}_{t_{\delta}}(\boldsymbol{x}) = \boldsymbol{f}(\boldsymbol{x}) \: - \: \boldsymbol{g}_{t_{\delta}}(\boldsymbol{x})$ we obtain
\begin{multline*}
F_{t_{\delta}, \delta}(\boldsymbol{x}) = \; \|\boldsymbol{\Sigma}^{-\frac{1}{2}}\left(\boldsymbol{f}(\boldsymbol{x}) - \boldsymbol{f}(\boldsymbol{x}_{true}) - \boldsymbol{\delta}\right)\|_{2}^{2} - 2\big\langle \boldsymbol{\Sigma}^{-\frac{1}{2}} \boldsymbol{g}_{t_{\delta}}(\boldsymbol{x}), \; \boldsymbol{\Sigma}^{-\frac{1}{2}}\left(\boldsymbol{f}(\boldsymbol{x}) - \boldsymbol{f}(\boldsymbol{x}_{true}) - \boldsymbol{\delta}\right) \big\rangle \\ 
+ \|\boldsymbol{\Sigma}^{-\frac{1}{2}}\boldsymbol{g}_{t_{\delta}}(\boldsymbol{x})\|_{2}^{2}.
\end{multline*}
Then a further decomposition of the first term on the right hand side yields
\begin{align*}
F_{t_{\delta}, \delta}(\boldsymbol{x}) = \; & F_{\infty, 0}(\boldsymbol{x}) - 2\big\langle \boldsymbol{\Sigma}^{-\frac{1}{2}} \boldsymbol{\delta}, \; \boldsymbol{\Sigma}^{-\frac{1}{2}}\left(\boldsymbol{f}(\boldsymbol{x}) - \boldsymbol{f}(\boldsymbol{x}_{true})\right) \big\rangle + \|\boldsymbol{\Sigma}^{-\frac{1}{2}}\boldsymbol{\delta}\|_{2}^{2} \\
- & 2\big\langle \boldsymbol{\Sigma}^{-\frac{1}{2}} \boldsymbol{g}_{t_{\delta}}(\boldsymbol{x}), \; \boldsymbol{\Sigma}^{-\frac{1}{2}}\left(\boldsymbol{f}(\boldsymbol{x}) - \boldsymbol{f}(\boldsymbol{x}_{true}) - \boldsymbol{\delta}\right) \big\rangle + \|\boldsymbol{\Sigma}^{-\frac{1}{2}}\boldsymbol{g}_{t_{\delta}}(\boldsymbol{x})\|_{2}^{2} \\
=: & \; F_{\infty, 0}(\boldsymbol{x}) + H_{t_{\delta}, \delta}(\boldsymbol{x}).
\end{align*}
From the limit $\lim_{\delta \to 0}\|\boldsymbol{\Sigma}^{-\frac{1}{2}}\boldsymbol{g}_{t_{\delta}}(\boldsymbol{x})\|_{2}^{2} = 0$, the limit $\lim_{\delta \to 0} \|\boldsymbol{\Sigma}^{-\frac{1}{2}}\boldsymbol{\delta}\|_{2}^{2} = 0$ and the boundedness of $\boldsymbol{f}(\boldsymbol{x})$ follows $\lim_{\delta \to 0} |H_{t_{\delta}, \delta}(\boldsymbol{x})| = 0$ for arbitrary but fixed $\boldsymbol{x} \in \Omega$.

Then the existence of the $\boldsymbol{x}_{t_{\delta}, \delta}$'s follows from Lemma \ref{MinimaCurve}.

\hfill $\square$
\end{prf}

At last we study how the minimizers $\boldsymbol{x}_{t_{\delta},\delta}$ of \eqref{nonlinreg_continuous} behave for $\delta \rightarrow 0$. We begin with a preparing corollary.

\begin{cor}
\label{minima_behavior}
\textit{Let the assumptions of Proposition \ref{minimizer_convergence} hold. Then we have for any local minimizer $\boldsymbol{x}_{t_{\delta}, \delta}$ of \eqref{nonlinreg_continuous} approximating a local minimizer $\boldsymbol{x}_{\infty, 0}$ of \eqref{nonlinreg_untruncated} for $\delta \rightarrow 0$ with $\|\boldsymbol{\Sigma}^{-\frac{1}{2}}(\boldsymbol{f}(\boldsymbol{x}_{\infty, 0}) - \boldsymbol{f}(\boldsymbol{x}_{true}))\|_{2} = 0$ that}
\begin{equation*}
\lim_{\delta \to 0}\|\boldsymbol{\Sigma}^{-\frac{1}{2}}(\boldsymbol{f}_{t_{\delta}}(\boldsymbol{x}_{t_{\delta}, \delta}) - \boldsymbol{f}(\boldsymbol{x}_{true}) - \boldsymbol{\delta})\|_{2} = 0.
\end{equation*}
\end{cor}

\begin{prf}
The assumptions of Proposition \ref{minimizer_convergence} give 
\begin{equation*}
\lim_{\delta \to 0}\|\boldsymbol{\Sigma}^{-\frac{1}{2}}\boldsymbol{\delta}\|_{2}^{2} = 0 \quad \text{and} \quad \lim_{\delta \to 0}\|\boldsymbol{\Sigma}^{-\frac{1}{2}}\boldsymbol{g}_{t_{\delta}}(\boldsymbol{x}_{t_{\delta}, \delta})\|_{2}^{2} = 0. 
\end{equation*}
We have by continuity of $\boldsymbol{f}(\boldsymbol{x})$ that $\lim_{\delta \to 0}\|\boldsymbol{\Sigma}^{-\frac{1}{2}}(\boldsymbol{f}(\boldsymbol{x}_{t_{\delta}, \delta}) - \boldsymbol{f}(\boldsymbol{x}_{true}))\|_{2} = 0$. Then
\begin{multline*}
\|\boldsymbol{\Sigma}^{-\frac{1}{2}}(\boldsymbol{f}(\boldsymbol{x}_{t_{\delta}, \delta}) - \boldsymbol{f}(\boldsymbol{x}_{true}))\|_{2} + \|\boldsymbol{\Sigma}^{-\frac{1}{2}}\boldsymbol{g}_{t_{\delta}}(\boldsymbol{x}_{t_{\delta}, \delta})\|_{2} + \|\boldsymbol{\Sigma}^{-\frac{1}{2}}\boldsymbol{\delta}\|_{2} \\
\geq \|\boldsymbol{\Sigma}^{-\frac{1}{2}}(\boldsymbol{f}_{t_{\delta}}(\boldsymbol{x}_{t_{\delta}, \delta}) - \boldsymbol{f}(\boldsymbol{x}_{true}) - \boldsymbol{\delta})\|_{2} 
\end{multline*}
gives the desired result.
\hfill $\square$
\end{prf}

\begin{prp}
\label{minima_isolated}
\textit{Let the assumptions of Proposition \ref{minimizer_convergence} hold. Assume that the local minimizers $\boldsymbol{x}_{\infty, 0}$ of \eqref{nonlinreg_untruncated} with $\|\boldsymbol{\Sigma}^{-\frac{1}{2}}(\boldsymbol{f}(\boldsymbol{x}_{\infty, 0}) - \boldsymbol{f}(\boldsymbol{x}_{true}))\|_{2} = 0$ form a discrete set $S_{\infty, 0}$. Then the set $L_{\infty, 0}$ consisting of the limits $\lim_{\delta \to 0}\boldsymbol{x}_{t_{\delta}, \delta}$ of local minimizers $\boldsymbol{x}_{t_{\delta}, \delta}$ of \eqref{nonlinreg_continuous} with $\lim_{\delta \to 0}\|\boldsymbol{\Sigma}^{-\frac{1}{2}}(\boldsymbol{f}_{t_{\delta}}(\boldsymbol{x}_{t_{\delta}, \delta}) - \boldsymbol{f}(\boldsymbol{x}_{true}) - \boldsymbol{\delta})\|_{2} = 0$ coincides with $S_{\infty, 0}$ and there exists a noise level $\delta_{max}$ such that all minimizers $\boldsymbol{x}_{t_{\delta}, \delta}$ approximating $S_{\infty, 0}$ are isolated for all $\delta \leq \delta_{max}$.}
\end{prp}

\begin{prf}
On the one hand from Proposition \ref{minimizer_convergence} we know that there exists a sequence $\boldsymbol{x}_{t_{\delta}, \delta}$ of minimizers of \eqref{nonlinreg_continuous} with $\lim_{\delta \to 0} \boldsymbol{x}_{t_{\delta}, \delta} = \boldsymbol{x}_{\infty, 0}$. Then $\|\boldsymbol{\Sigma}^{-\frac{1}{2}}(\boldsymbol{f}(\boldsymbol{x}_{\infty, 0}) - \boldsymbol{f}(\boldsymbol{x}_{true}))\|_{2} = 0$ and Corollary \ref{minima_behavior} give $\lim_{\delta \to 0}\|\boldsymbol{\Sigma}^{-\frac{1}{2}}(\boldsymbol{f}_{t_{\delta}}(\boldsymbol{x}_{t_{\delta}, \delta}) - \boldsymbol{f}(\boldsymbol{x}_{true}) - \boldsymbol{\delta})\|_{2} = 0$, which implies $S_{\infty, 0} \subseteq L_{\infty, 0}$. 

On the other hand holds for $\boldsymbol{x}_{t_{\delta}, \delta}$ with $\lim_{\delta \to 0}\|\boldsymbol{\Sigma}^{-\frac{1}{2}}(\boldsymbol{f}_{t_{\delta}}(\boldsymbol{x}_{t_{\delta}, \delta}) - \boldsymbol{f}(\boldsymbol{x}_{true}) - \boldsymbol{\delta})\|_{2} = 0$ that
\begin{multline*}
\|\boldsymbol{\Sigma}^{-\frac{1}{2}}(\boldsymbol{f}_{t_{\delta}}(\boldsymbol{x}_{t_{\delta}, \delta}) - \boldsymbol{f}(\boldsymbol{x}_{true}) - \boldsymbol{\delta})\|_{2} + \|\boldsymbol{\Sigma}^{-\frac{1}{2}}\boldsymbol{g}_{t_{\delta}}(\boldsymbol{x}_{t_{\delta}, \delta})\|_{2} + \|\boldsymbol{\Sigma}^{-\frac{1}{2}}\boldsymbol{\delta}\|_{2} \\ \geq \|\boldsymbol{\Sigma}^{-\frac{1}{2}}(\boldsymbol{f}(\boldsymbol{x}_{t_{\delta}, \delta}) - \boldsymbol{f}(\boldsymbol{x}_{true}))\|_{2}
\end{multline*}
which implies $\lim_{\delta \to 0}\|\boldsymbol{\Sigma}^{-\frac{1}{2}}(\boldsymbol{f}(\boldsymbol{x}_{t_{\delta}, \delta}) - \boldsymbol{f}(\boldsymbol{x}_{true}))\|_{2} = 0$. In particular this means by continuity of $\boldsymbol{f}(\boldsymbol{x})$ that the vector $\lim_{\delta \to 0}\boldsymbol{x}_{t_{\delta}, \delta}$ must be a local minimizer of $\|\boldsymbol{\Sigma}^{-\frac{1}{2}}(\boldsymbol{f}(\boldsymbol{x}) - \boldsymbol{f}(\boldsymbol{x}_{true}))\|_{2}$. Thus we have also shown $L_{\infty, 0} \subseteq S_{\infty, 0}$. 

In the following we number all elements of $S_{\infty, 0}$ with the index $k$, i.e. we write $\boldsymbol{x}_{\infty, 0}^{k}$ for $k = 1, ..., |S_{\infty, 0}|$. Similarly we number all minimizers $\boldsymbol{x}_{t_{\delta}, \delta}$ with $\lim_{\delta \to 0}\|\boldsymbol{\Sigma}^{-\frac{1}{2}}(\boldsymbol{f}_{t_{\delta}}(\boldsymbol{x}_{t_{\delta}, \delta}) - \boldsymbol{f}(\boldsymbol{x}_{true}) - \boldsymbol{\delta})\|_{2} = 0$ approximating the $\boldsymbol{x}_{\infty, 0}^{k}$'s with $\boldsymbol{x}_{t_{\delta}, \delta}^{k}$, i.e. $\lim_{\delta \to 0}\boldsymbol{x}_{t_{\delta}, \delta}^{k} = \boldsymbol{x}_{\infty, 0}^{k}$ for $k = 1, ..., |S_{\infty, 0}|$. Define 
\begin{equation*}
D_{min} := \min_{i \neq j} \|\boldsymbol{x}_{\infty, 0}^{i} - \boldsymbol{x}_{\infty, 0}^{j}\|_{2}.
\end{equation*}
Since $\lim_{\delta \to 0}\boldsymbol{x}_{t_{\delta}, \delta}^{k} = \boldsymbol{x}_{\infty, 0}^{k}$, we can find an error levels $\delta_{max}^{k}$ such that
\begin{equation*}
\|\boldsymbol{x}_{t_{\delta}, \delta}^{k} - \boldsymbol{x}_{\infty, 0}^{k}\|_{2} < \textstyle{\frac{1}{2}}D_{min} \quad \text{for} \quad k = 1, ..., |S_{\infty, 0}|,
\end{equation*}
which holds for all $0 \leq \delta \leq \delta_{max}^{k}$ for each $k$. Then for all $0 \leq \delta \leq \delta_{max} := \min_{k}\{\delta_{max}^{k}\}$ the $\boldsymbol{x}_{t_{\delta}, \delta}^{k}$'s must have pairwise mutual distances greater than zero.

\hfill $\square$
\end{prf}

Now Proposition \ref{minima_isolated} gives that the number of local minima $\boldsymbol{x}_{t_{\delta}, \delta}^{k}$ remains constant if the noise level $\delta$ is small enough. It also yields that these local minima then form a set of separated continuous curves parametrized in $\delta$.

At last we wish to have an estimate of the convergence of the local minima $\boldsymbol{x}_{t_{\delta}, \delta}^{k}$ of the truncated and noise contaminated problem to the local minima $\boldsymbol{x}_{\infty, 0}^{k}$ of the noise-free and untruncated problem, which is useful for practical computations.

\begin{prp}
\label{ApproxEstimate}
\textit{For the noise level $\delta$ small enough, we can bound for any local minimum $\boldsymbol{x}_{t_{\delta}, \delta}^{k}$ the approximation error $\|\boldsymbol{x}_{\infty, 0}^{k} - \boldsymbol{x}_{t_{\delta}, \delta}^{k}\|_{2}$ with a positively weighted linear combination of the residual $\|\boldsymbol{\Sigma}^{-\frac{1}{2}}(\boldsymbol{f}_{t_{\delta}}(\boldsymbol{x}_{t_{\delta}, \delta}^{k}) - \boldsymbol{f}(\boldsymbol{x}_{true}) - \boldsymbol{\delta})\|_{2}$, the truncation error $\|\boldsymbol{\Sigma}^{-\frac{1}{2}}\boldsymbol{g}_{t_{\delta}}(\boldsymbol{x}_{t_{\delta}, \delta}^{k})\|_{2}$ and the noise estimate $\|\boldsymbol{\Sigma}^{-\frac{1}{2}}\boldsymbol{\delta}\|_{2}$.}
\end{prp}

\begin{prf}
The first order necessary conditions for a local minimum of $F_{t_{\delta}, \delta}(\boldsymbol{x}) = F_{\infty, 0}(\boldsymbol{x}) + H_{t_{\delta}, \delta}(\boldsymbol{x})$ at $\boldsymbol{x}_{t_{\delta}, \delta}^{k}$ and a local minimum of $F_{\infty, 0}(\boldsymbol{x})$ at $\boldsymbol{x}_{\infty, 0}^{k}$ yield in particular
\begin{align*}
\big\langle \nabla F_{\infty, 0}(\boldsymbol{x}_{t_{\delta}, \delta}^{k}) + \nabla H_{t_{\delta}, \delta}(\boldsymbol{x}_{t_{\delta}, \delta}^{k}), \; \boldsymbol{x}_{\infty, 0}^{k} - \boldsymbol{x}_{t_{\delta}, \delta}^{k} \big\rangle & \geq 0 \\
\text{and} \quad \big\langle \nabla F_{\infty, 0}(\boldsymbol{x}_{\infty, 0}^{k}), \; \boldsymbol{x}_{t_{\delta}, \delta}^{k} - \boldsymbol{x}_{\infty, 0}^{k} \big\rangle & \geq 0,
\end{align*}
Adding the last two inequalities yields
\begin{equation}
\label{gradient_inequality}
\big\langle \nabla H_{t_{\delta}, \delta}(\boldsymbol{x}_{t_{\delta}, \delta}^{k}), \; \boldsymbol{x}_{\infty, 0}^{k} - \boldsymbol{x}_{t_{\delta}, \delta}^{k} \big\rangle \geq \big\langle \nabla F_{\infty, 0}(\boldsymbol{x}_{\infty, 0}^{k}) - \nabla F_{\infty, 0}(\boldsymbol{x}_{t_{\delta}, \delta}^{k}), \; \boldsymbol{x}_{\infty, 0}^{k} - \boldsymbol{x}_{t_{\delta}, \delta}^{k} \big\rangle.
\end{equation}
Since $\nabla F_{\infty, 0}(\boldsymbol{x})$ is totally differentiable at $\boldsymbol{x}_{\infty, 0}^{k}$ we obtain
\begin{equation*}
\nabla F_{\infty, 0}(\boldsymbol{x}_{t_{\delta}, \delta}^{k}) = \nabla F_{\infty, 0}(\boldsymbol{x}_{\infty, 0}^{k}) + \mathrm{Hess}_{F_{\infty, 0}}(\boldsymbol{x}_{\infty, 0}^{k})\left(\boldsymbol{x}_{t_{\delta}, \delta}^{k} - \boldsymbol{x}_{\infty, 0}^{k}\right) + \boldsymbol{w}_{\infty, 0}(\boldsymbol{x}_{t_{\delta}, \delta}^{k}, \boldsymbol{x}_{\infty, 0}^{k}),
\end{equation*}
where $\boldsymbol{w}_{\infty, 0}(\boldsymbol{x}, \boldsymbol{x}_{\infty, 0}^{k})$ fulfills 
\begin{align*}
\|\boldsymbol{w}_{\infty, 0}(\boldsymbol{x}, \boldsymbol{x}_{\infty, 0}^{k})\|_{2} \leq \| \boldsymbol{x} - \boldsymbol{x}_{\infty, 0}^{k}\|_{2}\epsilon_{\infty, 0}(\boldsymbol{x}, \boldsymbol{x}_{\infty, 0}^{k}) \quad \text{with} \quad \lim_{\boldsymbol{x} \to \boldsymbol{x}_{\infty, 0}^{k}} \epsilon_{\infty, 0}(\boldsymbol{x}, \boldsymbol{x}_{\infty, 0}^{k}) = 0.
\end{align*}

Since $\mathrm{Hess}_{F_{\infty, 0}}(\boldsymbol{x}_{\infty, 0}^{k})$ is positive definite, the expression $\left(\big\langle \boldsymbol{x}, \mathrm{Hess}_{F_{\infty, 0}}(\boldsymbol{x}_{\infty, 0}^{k})\:\boldsymbol{x} \big\rangle \right)^{\frac{1}{2}}$ gives a norm on $\mathbb{R}^{D}$. Because of the equivalence of all norms in $\mathbb{R}^{D}$, there exists a constant $C_{\infty, 0}^{k} > 0$ with 
\begin{align*}
\left(\big\langle \boldsymbol{x}, \mathrm{Hess}_{F_{\infty, 0}}(\boldsymbol{x}_{\infty, 0}^{k})\:\boldsymbol{x} \big\rangle\right)^{\frac{1}{2}} \geq C_{\infty, 0}^{k}\|\boldsymbol{x}\|_{2} \quad \text{for all} \quad \boldsymbol{x} \in \mathbb{R}^{D}.
\end{align*}

Since $\lim_{\delta \to 0} \boldsymbol{x}_{t_{\delta}, \delta}^{k} = \boldsymbol{x}_{\infty, 0}^{k}$ we can find a noise level $\rho_{max}^{k}$ such that $|\epsilon_{\infty, 0}(\boldsymbol{x}_{t_{\delta}, \delta}^{k}, \boldsymbol{x}_{\infty, 0}^{k})| \leq d_{\infty, 0}^{k}$ for all $\delta \leq \rho_{max}^{k}$, where $d_{\infty, 0}^{k}$ is a constant with $0 \leq d_{\infty, 0}^{k} < \left(C_{\infty, 0}^{k}\right)^{2}$. Then using 
\begin{equation*}
\big\langle \boldsymbol{w}_{\infty, 0}(\boldsymbol{x}_{t_{\delta}, \delta}^{k}, \boldsymbol{x}_{\infty, 0}^{k}), \boldsymbol{x}_{\infty, 0}^{k} - \boldsymbol{x}_{t_{\delta}, \delta}^{k} \big\rangle \leq \epsilon_{\infty, 0}(\boldsymbol{x}_{t_{\delta}, \delta}^{k}, \boldsymbol{x}_{\infty, 0}^{k})\|\boldsymbol{x}_{\infty, 0}^{k} - \boldsymbol{x}_{t_{\delta}, \delta}^{k}\|_{2}^{2} 
\end{equation*}
and \eqref{gradient_inequality} we can estimate 
\begin{align*}
& \|\nabla H_{t_{\delta}, \delta}(\boldsymbol{x}_{t_{\delta}, \delta}^{k})\|_{2}\|\boldsymbol{x}_{\infty, 0}^{k} - \boldsymbol{x}_{t_{\delta}, \delta}^{k}\|_{2} \\
\geq \; & \big\langle \nabla F_{\infty, 0}(\boldsymbol{x}_{\infty, 0}^{k}) - \nabla F_{\infty, 0}(\boldsymbol{x}_{t_{\delta}, \delta}^{k}), \boldsymbol{x}_{\infty, 0}^{k} - \boldsymbol{x}_{t_{\delta}, \delta}^{k} \big\rangle \\
= \; & \big\langle \mathrm{Hess}_{F_{\infty, 0}}(\boldsymbol{x}_{\infty, 0}^{k})\left(\boldsymbol{x}_{\infty, 0}^{k} - \boldsymbol{x}_{t_{\delta}, \delta}^{k}\right) - \boldsymbol{w}_{\infty, 0}(\boldsymbol{x}_{t_{\delta}, \delta}^{k}, \boldsymbol{x}_{\infty, 0}^{k}), \; \boldsymbol{x}_{\infty, 0}^{k} - \boldsymbol{x}_{t_{\delta}, \delta}^{k} \big\rangle \\
\geq \; & \left(\big(C_{\infty, 0}^{k}\big)^{2} - \epsilon_{\infty, 0}(\boldsymbol{x}_{t_{\delta}, \delta}^{k}, \boldsymbol{x}_{\infty, 0}^{k})\right)\|\boldsymbol{x}_{\infty, 0}^{k} - \boldsymbol{x}_{t_{\delta}, \delta}^{k}\|_{2}^{2},
\end{align*}
i.e. this gives
\begin{equation}
\label{distance_estimate} 
\|\boldsymbol{x}_{\infty, 0}^{k} - \boldsymbol{x}_{t_{\delta}, \delta}^{k}\|_{2} \leq \left(\big(C_{\infty, 0}^{k}\big)^{2} - d_{\infty, 0}^{k}\right)^{-1}\|\nabla H_{t_{\delta}, \delta}(\boldsymbol{x}_{t_{\delta}, \delta}^{k})\|_{2}
\end{equation}
for all $\delta \leq \rho_{max}^{k}$.

We have
\begin{align*}
\nabla H_{t_{\delta}, \delta}(\boldsymbol{x}) = \; 2\Big(&\mathrm{Jac}_{\boldsymbol{g}_{t_{\delta}}}^{T}(\boldsymbol{x})\boldsymbol{\Sigma}^{-1}\boldsymbol{g}_{t_{\delta}}(\boldsymbol{x}) - \mathrm{Jac}_{\boldsymbol{g}_{t_{\delta}}}^{T}(\boldsymbol{x})\boldsymbol{\Sigma}^{-1}(\boldsymbol{f}(\boldsymbol{x}) - \boldsymbol{f}(\boldsymbol{x}_{true}) - \boldsymbol{\delta})\\
- &\mathrm{Jac}_{\boldsymbol{f}}^{T}(\boldsymbol{x})\boldsymbol{\Sigma}^{-1}\boldsymbol{g}_{t_{\delta}}(\boldsymbol{x}) - \mathrm{Jac}_{\boldsymbol{f}}^{T}(\boldsymbol{x})\boldsymbol{\Sigma}^{-1}\boldsymbol{\delta}\Big),
\end{align*}
i.e. we can find constants $K_{1}^{k}$ and $K_{2}^{k}$ such that
\begin{align*}
\|\nabla H_{t_{\delta}, \delta}(\boldsymbol{x}_{t_{\delta}, \delta}^{k})\|_{2} \leq \; & 2\Big(K_{1}^{k}\|\mathrm{Jac}_{\boldsymbol{g}_{t_{\delta}}}^{T}(\boldsymbol{x}_{t_{\delta}, \delta}^{k})\boldsymbol{\Sigma}^{-\frac{1}{2}}\|_{\infty}\|\boldsymbol{\Sigma}^{-\frac{1}{2}}\boldsymbol{g}_{t_{\delta}}(\boldsymbol{x}_{t_{\delta}, \delta}^{k})\|_{2} \\
& + K_{1}^{k}\|\mathrm{Jac}_{\boldsymbol{g}_{t_{\delta}}}^{T}(\boldsymbol{x}_{t_{\delta}, \delta}^{k})\boldsymbol{\Sigma}^{-\frac{1}{2}}\|_{\infty}\big(\|\boldsymbol{\Sigma}^{-\frac{1}{2}}(\boldsymbol{f}_{t_{\delta}}(\boldsymbol{x}_{t_{\delta}, \delta}^{k}) - \boldsymbol{f}(\boldsymbol{x}_{true}) - \boldsymbol{\delta})\|_{2} \\
& + \|\boldsymbol{\Sigma}^{-\frac{1}{2}}\boldsymbol{g}_{t_{\delta}}(\boldsymbol{x}_{t_{\delta}, \delta}^{k})\|_{2}\big) \\
& + K_{2}^{k}\|\mathrm{Jac}_{\boldsymbol{f}}^{T}(\boldsymbol{x}_{t_{\delta}, \delta}^{k})\boldsymbol{\Sigma}^{-\frac{1}{2}}\|_{\infty}\|\boldsymbol{\Sigma}^{-\frac{1}{2}}\boldsymbol{g}_{t_{\delta}}(\boldsymbol{x}_{t_{\delta}, \delta}^{k})\|_{2} \\
& + K_{2}^{k}\|\mathrm{Jac}_{\boldsymbol{f}}^{T}(\boldsymbol{x}_{t_{\delta}, \delta}^{k})\boldsymbol{\Sigma}^{-\frac{1}{2}}\|_{\infty}\|\boldsymbol{\Sigma}^{-\frac{1}{2}}\boldsymbol{\delta}\|_{2}\Big),
\end{align*}
which gives the result.
\hfill $\square$
\end{prf}

\begin{cor}
\label{BoundOnExpectedValue}
\textit{Let the truncation indices $t_{\delta}$ depend on the vector of independent Gaussian random variables $\boldsymbol{\delta}$ with $\lim_{\delta \to 0} \mathbb{E}\big(\|\boldsymbol{\Sigma}^{-\frac{1}{2}}\boldsymbol{\delta}\|_{2}^{2}\big) = 0$ such that $\lim_{\delta \to 0} \mathbb{E}\big(\|\boldsymbol{\Sigma}^{-\frac{1}{2}}\boldsymbol{g}_{t_{\delta}}(\boldsymbol{x})\|_{2}^{2}\big) = 0$ holds for all arbitrary but fixed $\boldsymbol{x} \in \Omega$. Then we have for all minimizers $\boldsymbol{x}_{\infty, 0}^{k}$ of \eqref{nonlinreg_untruncated} with $\|\boldsymbol{\Sigma}^{-\frac{1}{2}}(\boldsymbol{f}(\boldsymbol{x}_{\infty, 0}^{k}) - \boldsymbol{f}(\boldsymbol{x}_{true}))\|_{2} = 0$ that for $\delta$ sufficiently small there exist minimizers $\boldsymbol{x}_{t_{\delta}, \delta}^{k}$ of \eqref{nonlinreg_continuous} with} 
\begin{equation*}
\lim_{\delta \to 0}\mathbb{E}\big(\|\boldsymbol{x}_{\infty, 0}^{k} - \boldsymbol{x}_{t_{\delta}, \delta}^{k}\|_{2}\big) = 0.
\end{equation*}
\end{cor}

\begin{prf}
Proposition \ref{minimizer_convergence} establishes the existence of the $\boldsymbol{x}_{t_{\delta}, \delta}^{k}$'s. Corollary \ref{minima_behavior} gives
\begin{align*}
\lim_{\delta \to 0}\mathbb{E}\big(\|\boldsymbol{\Sigma}^{-\frac{1}{2}}(\boldsymbol{f}_{t_{\delta}}(\boldsymbol{x}_{t_{\delta}, \delta}^{k}) - \boldsymbol{f}(\boldsymbol{x}_{true}) - \boldsymbol{\delta})\|_{2}\big) = 0. 
\end{align*}
Proposition \ref{minima_isolated} gives that for $\delta$ sufficiently small there exists a constant $K_{\infty, 0}^{k}$ with
\begin{equation*}
\mathbb{E}\big(\|\boldsymbol{x}_{\infty, 0}^{k} - \boldsymbol{x}_{t_{\delta}, \delta}^{k}\|_{2}\big) \leq K_{\infty, 0}^{k}\mathbb{E}\big(\|\nabla H_{t_{\delta}, \delta}(\boldsymbol{x}_{t_{\delta}, \delta}^{k})\|_{2}\big).
\end{equation*}
Set $S_{1}^{k} := \sup_{\boldsymbol{x} \in \Omega}\|\mathrm{Jac}_{\boldsymbol{g}_{t_{\delta}}}^{T}(\boldsymbol{x})\boldsymbol{\Sigma}^{-\frac{1}{2}}\|_{\infty} < \infty$ and $S_{2}^{k} := \sup_{\boldsymbol{x} \in \Omega}\|\mathrm{Jac}_{\boldsymbol{f}}^{T}(\boldsymbol{x})\boldsymbol{\Sigma}^{-\frac{1}{2}}\|_{\infty} < \infty$. Then the estimate for $\|\nabla H_{t_{\delta}, \delta}(\boldsymbol{x}_{t_{\delta}, \delta}^{k})\|_{2}$ in the proof of Proposition \ref{minima_isolated} gives
\begin{align*}
\mathbb{E}\big(\|\nabla H_{t_{\delta}, \delta}(\boldsymbol{x}_{t_{\delta}, \delta}^{k})\|_{2}\big) \leq \; & 2\Big(K_{1}^{k}S_{1}^{k}\mathbb{E}\big(\|\boldsymbol{\Sigma}^{-\frac{1}{2}}\boldsymbol{g}_{t_{\delta}}(\boldsymbol{x}_{t_{\delta}, \delta}^{k})\|_{2}\big) \\
& + K_{1}^{k}S_{1}^{k}\Big(\mathbb{E}\big(\|\boldsymbol{\Sigma}^{-\frac{1}{2}}(\boldsymbol{f}_{t_{\delta}}(\boldsymbol{x}_{t_{\delta}, \delta}^{k}) - \boldsymbol{f}(\boldsymbol{x}_{true}) - \boldsymbol{\delta})\|_{2}\big) \\
& + \mathbb{E}\big(\|\boldsymbol{\Sigma}^{-\frac{1}{2}}\boldsymbol{g}_{t_{\delta}}(\boldsymbol{x}_{t_{\delta}, \delta}^{k})\|_{2}\big)\Big) \\
& + K_{2}^{k}S_{2}^{k}\mathbb{E}\big(\|\boldsymbol{\Sigma}^{-\frac{1}{2}}\boldsymbol{g}_{t_{\delta}}(\boldsymbol{x}_{t_{\delta}, \delta}^{k})\|_{2}\big) \\
 & + K_{2}^{k}S_{2}^{k}\mathbb{E}\big(\|\boldsymbol{\Sigma}^{-\frac{1}{2}}\boldsymbol{\delta}\|_{2}\big)\Big),
\end{align*}
which proves the claim since Assumption \ref{covariance} gives $\mathbb{E}\big(\|\boldsymbol{\Sigma}^{-\frac{1}{2}}\boldsymbol{\delta}\|_{2}\big) \leq \sqrt{N}\delta$.

\hfill $\square$
\end{prf}

The strategy for our retrieval algorithm is to start with an initial guess for the truncation index $t_{start}$ and try to find all local minima $\boldsymbol{x}_{t_{start}, \delta}^{k}$. Then the truncation index is gradually increased and starting from $\boldsymbol{x}_{t_{start}, \delta}^{k}$ the continuation method is applied to find finally the local minima $\boldsymbol{x}_{t_{\delta}, \delta}^{k}$. Motivated by Propositions \ref{ApproxEstimate} and \ref{BoundOnExpectedValue} only those local minima are considered to be possible approximations to our sought-after refractive index, where the residual $\|\boldsymbol{\Sigma}^{-\frac{1}{2}}(\boldsymbol{f}_{t_{\delta}}(\boldsymbol{x}_{t_{\delta}, \delta}^{k}) - \boldsymbol{f}(\boldsymbol{x}_{true}) - \boldsymbol{\delta})\|_{2}$ and an estimate of the truncation error $\|\boldsymbol{\Sigma}^{-\frac{1}{2}}\boldsymbol{g}_{t_{\delta}}(\boldsymbol{x}_{t_{\delta}, \delta}^{k})\|_{2}$ are both reasonably small. The latter serves also as a stopping criterion for the continuation method 

The initial guess $t_{start}$ has to be selected with care. On the one hand if it is to small, the model is to inaccurate and the retrieval of the sought-after local minima can not be guaranteed. On the other hand if it is too big, computational effort is wasted, since too many Mie coefficient functions with almost vanishing magnitudes and thus essentially not changing the local minima are computed.

\section{The Reconstruction Algorithm}
\label{recon_alg}

We now return to our regression problem \eqref{nonlinregprob}. For $i = 1, ..., N$ we see that the measured extinctions normalized by the number of particles $n_{i}$ with radius $r_{i}$, i.e. the quantity $\frac{e_{i}}{n_{i}}$, is Gaussian-distributed with mean $\frac{1}{n_{i}}q_{true}(r_{i},l)$ and standard deviation $\sigma_{i} := \frac{s_{i}}{n_{i}}$. In the following we fix a wavelength $l$, i.e. we reconstruct the sought-after particle refractive index $m_{part}(l)$ wavelength by wavelength. In the following the unit both for particle radii and light wavelengths is $\mu$m.

We make use of the function $\boldsymbol{q}_{N_{tr}}: \mathbb{R}^{2} \rightarrow \mathbb{R}^{N}$ defined by
\begin{equation}
\boldsymbol{q}_{N_{tr}}(\boldsymbol{x}) := \left(\pi r_{1}^{2}\sum_{n=1}^{N_{tr}}q_{n}(\boldsymbol{x}, r_{1}, l), \; ..., \; \pi r_{N}^{2}\sum_{n=1}^{N_{tr}}q_{n}(\boldsymbol{x}, r_{N}, l)\right)^{T},
\end{equation}
where $R = \{r_{1}, ..., r_{N}\}$ is the particle radius grid and $N_{tr}$ the truncation index to be used. We allow non-integer truncation indices $N_{tr}$ as well, where the non-integer truncation is done like in Proposition \ref{MinimaContinuity}. Here the expression $q_{n}(\boldsymbol{x}, r_{k}, l)$ is a short notation of $q_{n}(m_{med}(l), (\boldsymbol{x})_{1} + (\boldsymbol{x})_{2}i, r_{k}, l)$ from Section \ref{experiment}, where the sought-after refrative index $m_{part}(l)$ is identified with the vector $\boldsymbol{x}$ here, i.e. $m_{part}(l) =  (\boldsymbol{x})_{1} + (\boldsymbol{x})_{2}i$. So its computation follows Section \ref{Mie_theory}.

In the following the refractive index search area is given by the rectangle $\Omega := [0, 20]\times[0,40]$, which means that we only consider refractive indices of particle materials whose real parts lie in the interval $[0,20]$ and its imaginary parts in the interval $[0,40]$. This rather large search area makes the algorithm suitable for a wide range of aerosol materials.

\begin{algorithm}[H]
\caption{Reconstruction of Refractive Indices}
\label{reconstruction_algorithm}
\begin{algorithmic}[1]
\State $b_{real} = 20$
\State $b_{imag} = 40$
\State $N_{real} = 81$
\State $N_{imag} = 161$
\State $c_{i} = (i - 1)\frac{b_{real}}{N_{real} - 1}$ for $i = 1, ..., N_{real}$ \label{SearchGridReal}
\State $d_{i} = (i - 1)\frac{b_{imag}}{N_{imag} - 1}$ for $i = 1, ..., N_{imag}$ \label{SearchGridImag}
\State $R = \{0.1, 0.2, 0.3\}$
\State $N = 3$
\State $N_{tr} = 3$
\State $S_{start} = \{\}$
\State estimate $\sigma_{1}^{2}$, ..., $\sigma_{N}^{2}$ from the sample means approximating the standard deviations of $\frac{e_{1}}{n_{1}}$, ..., $\frac{e_{N}}{n_{N}}$. \vspace{0.05cm}  
\State $\delta^{2} := \mathrm{max}\big\{\sigma_{1}^{2}, ..., \sigma_{N}^{2}\big\}$\vspace{0.05cm}
\State $\boldsymbol{\Sigma} := \delta^{-2} \cdot \mathrm{diag}\big(\sigma_{1}^{2}, ..., \sigma_{N}^{2}\big)$\vspace{0.05cm}
     \For{$i = 1 \; \textbf{to} \; N_{real}$} \label{FirstLoopStart}	
	\For{$j = 1 \; \textbf{to} \; N_{imag}$}	
	\State compute the Hessian $\boldsymbol{H}(c_{i}, d_{j})$ of $F(\boldsymbol{x}) := \frac{1}{2}\|\boldsymbol{\Sigma}^{-\frac{1}{2}}\left(\boldsymbol{q}_{N_{tr}}(\boldsymbol{x}) - \boldsymbol{e}_{real}\right)\|_{2}^{2}$ \newline \hspace*{1.2cm} at $\boldsymbol{x} = (c_{i}, d_{j})^{T}$ \vspace{0.05cm} \label{FitFunction}
	      \If{$\boldsymbol{H}(c_{i}, d_{j})$ is positive definite} \vspace{0.05cm}
	      \State use $ (c_{i}, d_{j})^{T}$ as start vector to compute \vspace{0.05cm}
	      \State $\boldsymbol{x}_{new} = \underset{\boldsymbol{x} \in \mathbb{R}^{2}}{\mathrm{argmin}}\frac{1}{2}\|\boldsymbol{\Sigma}^{-\frac{1}{2}}\left(\boldsymbol{q}_{N_{tr}}(\boldsymbol{x}) - \boldsymbol{e}\right)\|_{2}^{2} \quad \text{s.t.} \quad \boldsymbol{x} \in [0, b_{real}] \times [0, b_{imag}]$
		\If{$S_{start}$ is empty $\lor$ $\frac{\|\boldsymbol{x} - \boldsymbol{x}_{new}\|_{2}}{\|\boldsymbol{x}\|_{2}} \geq 10^{-2} \quad \forall \boldsymbol{x} \in S_{start}$}\vspace{0.05cm} \label{CheckMinimum}
		\State $S_{start} = S_{start} \cup \{\boldsymbol{x}_{new}\}$
		\EndIf
                \EndIf
	\EndFor
     \EndFor \label{FirstLoopEnd}
\State $N_{start} = \left|S_{start}\right|$
\State $S_{out} = \{\}$
      \algstore{myalg_a}
  \end{algorithmic}
\end{algorithm}  


\begin{algorithm}
  \begin{algorithmic}
      \algrestore{myalg_a}
\State $\tau = 3$
     \For{$i = 1 \; \textbf{to} \; N_{start}$} \label{SecondLoopStart}
     \State $c = N_{tr}$
     \State $Tol_{rel} = 10^{-3}$
     \State $D_{rel} = \infty$
	\While{$D_{rel} > Tol_{rel}$}
	     \For{$p = 1 \; \textbf{to} \; 10$} 
	     \State use the vector $S_{start}(i)$ as start vector to compute
	     \State $\boldsymbol{x}_{new} = \underset{\boldsymbol{x} \in \mathbb{R}^{2}}{\mathrm{argmin}}\frac{1}{2}\|\boldsymbol{\Sigma}^{-\frac{1}{2}}\left(\boldsymbol{q}_{c + \frac{p}{10}}(\boldsymbol{x}) - \boldsymbol{e}\right)\|_{2}^{2} \quad \text{s.t.} \quad \boldsymbol{x} \in [0, b_{real}] \times [0, b_{imag}]$ \label{LocalSolver}
	     \State $Res_{cur} = \|\boldsymbol{\Sigma}^{-\frac{1}{2}}\left(\boldsymbol{q}_{c + \frac{p-1}{10}}(S_{start}(i)) - \boldsymbol{e}\right)\|_{2}^{2}$
	     \State $Res_{new} = \|\boldsymbol{\Sigma}^{-\frac{1}{2}}\left(\boldsymbol{q}_{c + \frac{p}{10}}(\boldsymbol{x}_{new}) - \boldsymbol{e}\right)\|_{2}^{2}$\vspace{0.05cm}
	     \State $D_{rel} = \frac{\left|Res_{cur} - Res_{new}\right|}{Res_{cur}}$\vspace{0.05cm}
	     \State $S_{start}(i) = \boldsymbol{x}_{new}$
	     \EndFor
           \State $c = c + 1$
	\EndWhile
	\If{$Res_{new} < \tau N \delta^{2}$} \label{CheckNewMinimumStart}
	     \If{$S_{out}$ is empty $\lor$ $\frac{\|\boldsymbol{x} - \boldsymbol{x}_{new}\|_{2}}{\|\boldsymbol{x}\|_{2}} \geq 10^{-2} \quad \forall \boldsymbol{x} \in S_{out}$}\vspace{0.05cm}
	     \State $S_{out} = S_{out} \cup \{\boldsymbol{x}_{new}\}$
	     \EndIf
	 \EndIf \label{CheckNewMinimumEnd}
     \EndFor \label{SecondLoopEnd}
\end{algorithmic}
\end{algorithm}

In the first loop from lines \ref{FirstLoopStart} - \ref{FirstLoopEnd} a search for local minima of the fit function $F(\boldsymbol{x})$ defined in line \ref{FitFunction} for the truncation index $N_{tr} = 3$ is performed. The loop runs through all grid points $ (c_{i}, d_{j})^{T}$ of the search grid defined in lines \ref{SearchGridReal} - \ref{SearchGridImag}. If the Hessian of $F(\boldsymbol{x})$ at some grid point $ (c_{i}, d_{j})^{T}$ is positive definite, this point might lie in the vicinity of a local minimum. The Hessian is computed exactly, where the second partial derivatives of the Mie extinction efficiency with respect to the real and imaginary part of the scattering material needed here are computed using the product rule approach from Section \ref{derivatives}. So we use $ (c_{i}, d_{j})^{T}$ as start point for a local solver in this case. In line \ref{CheckMinimum} we only accept a new local minimum if it is sufficiently different from the local minima already found. Then it is stored in the container $S_{start}$. This simple global search strategy can find all local minima if the search grid is fine enough. 

\begin{rem}
Of course other well established global search heuristics can be applied here as well. In test runs we compared genetic algorithms with our sequential search strategy on the two-dimensional refractive index search area, but their computational effort and reliability remained the same or were even worse. In \cite{EHR11} the technique of simulated annealing was used to retrieve aerosol refractive indices, which could be a promising alternative here. In our study the measurement noise was so high, that a unique global minimum of our fit function could not be determined. Instead our focus lied on effectively finding all local minima with small values of the fit function and we regarded them all as possible approximations to a thought-after refractive index. 
\end{rem}

The second loop from lines \ref{SecondLoopStart} - \ref{SecondLoopEnd} uses the local minima found in the first loop as start points for the continuation method following Proposition \ref{MinimaContinuity} and Corollary \ref{continuation}. We found that a step width of $0.1$ is for our problem a well-balanced choice between too big step widths rendering the continuation method unstable and too small step widths making it computationally inefiicient. With the stopping criterion $D_{rel} \leq Tol_{rel}$ of the while-loop it is approximately checked if the magnitude of the remainder term is small enough. Finally in line \ref{CheckNewMinimumStart} it is checked if the residual is small enough. In our implementation we did another run of lines \ref{CheckNewMinimumStart} - \ref{CheckNewMinimumEnd} with $\tau = 5$ and $\tau = 7$ respectively, if none of the reconstructions had a squared residual smaller than $\tau N_{r} \delta^{2}$ for the previous $\tau$. This had to be done, because the parameter $\tau$ has to be selected carefully in order to estimate the bound on $\mathbb{E}\big(\|\boldsymbol{x}_{\infty, 0}^{k} - \boldsymbol{x}_{t_{\delta}, \delta}^{k}\|_{2}\big)$ derived in the proof of Corollary \ref{BoundOnExpectedValue} correctly. 

\section{Comparison of the Numerical Continuation Approach with Established Truncation Index Heuristics}
\label{comparison_with_trunc_ind}

As solution of the forward problem we generated for a discrete set of wavelengths $l_{1}$, ..., $l_{N_{l}}$ unperturbed spectral extinctions normalized with the number of particles of the monodisperse aerosol by computing 
\begin{equation*}
\left(\boldsymbol{e}_{true}\right)_{i,j} := \pi r_{i}^{2}\sum_{n = 1}^{N_{tr}}q_{n}(m_{med}(l_{j}), m_{part}(l_{j}), r_{i}, l_{j}), \quad \text{for} \quad i = 1, ..., N, \quad j = 1, ..., N_{l} 
\end{equation*}
with $m_{part}(l_{i})$ taken as the refractive indices of Ag, $\mathrm{H}_{2}$O and CsI. Here we used the truncation index
\begin{equation}
\begin{split}
& \rho = 2\pi \frac{r}{l}, \\
& M = \max\{|\rho|, |\rho \cdot m_{part}(l)|, |\rho \cdot m_{med}(l)|\}, \\
& N_{tr} := \lceil|M + 4.05 \cdot M^{\frac{1}{3}} + 2|\rceil
\end{split}
\label{truncation_index}
\end{equation}
introduced in \cite{Wi80}.

For particle size distribution reconstructions as outlined in \cite{AK16} we need particle refractive indices for five optical windows, see \cite{GY89}, so the wavelength grid of interest consists of five ranges. These ranges are given by $8$ linearly spaced wavelengths from $0.6 - 0.8 \; \mu$m, $8$ from $1.1 - 1.3 \; \mu$m, $8$ from $1.6 - 1.8 \; \mu$m, $16$ from $2.1 - 2.5 \; \mu$m and $8$ from $3.1 - 3.3 \; \mu$m, so we have in total $N_{l} = 48$ wavelengths.

For each of the $48$ wavelengths we generated noisy spectral extinctions $\boldsymbol{e}$ by adding zero-mean Gaussian noise to $\boldsymbol{e}_{true}$, i.e.
\begin{equation*}
\left(\boldsymbol{e}\right)_{i,j} = \left(\boldsymbol{e}_{true}\right)_{i,j} + \delta_{i,j} \quad \text{ with } \; \delta_{i,j} \sim \mathcal{N}(0,(0.05 \cdot \left(\boldsymbol{e}_{true}\right)_{i,j})^{2}), \quad i = 1, ..., N \quad j = 1, ..., N_{l}.
\end{equation*}
Here the standard deviations were taken to be $5\%$ of the original extinction values. We computed each mean $\left(\boldsymbol{e}_{real}\right)_{i,j}$ of the noisy spectral extinctions with a sample size of $N_{s} = 300$ .

In the following Algortihm \ref{reconstruction_algorithm} is referred to as method 1. On the same simulated spectral extinctions we let Algorithm \ref{reconstruction_algorithm} run up to line \ref{FirstLoopEnd}, but with the difference that at each evaluation of $\boldsymbol{q}_{N_{tr}}(\boldsymbol{x})$ we directly took the trunction index from \eqref{truncation_index}. We denote this approach with method 2. We now display the average run times of method 1 and method 2 for $10$ sweeps through all $48$ wavelenghs. 

\subsubsection{Results for Ag}

\begin{figure}[h!]
\begingroup
\sbox0{\includegraphics[width =1.0\textwidth]{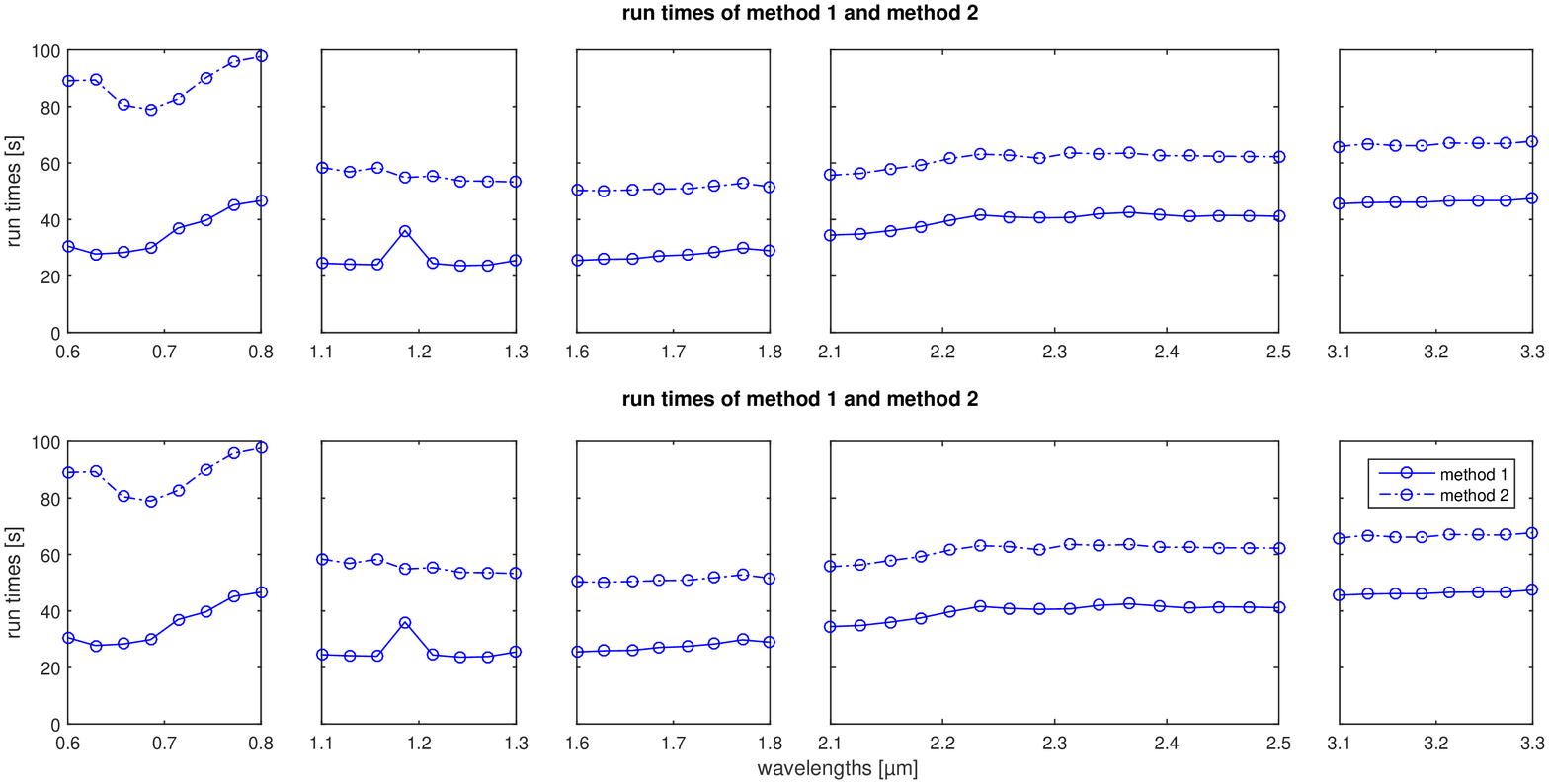}}
\includegraphics[clip,trim=0 0 0 {.5\wd0},width=1.0\textwidth]{RunTimesMethods1And2Silver.eps}
\endgroup
\label{RunTimesMethods1And2Silver}
\end{figure}

\newpage

\subsubsection{results for CsI}

\begin{figure}[h!]
\begingroup
\sbox0{\includegraphics[width =1.0\textwidth]{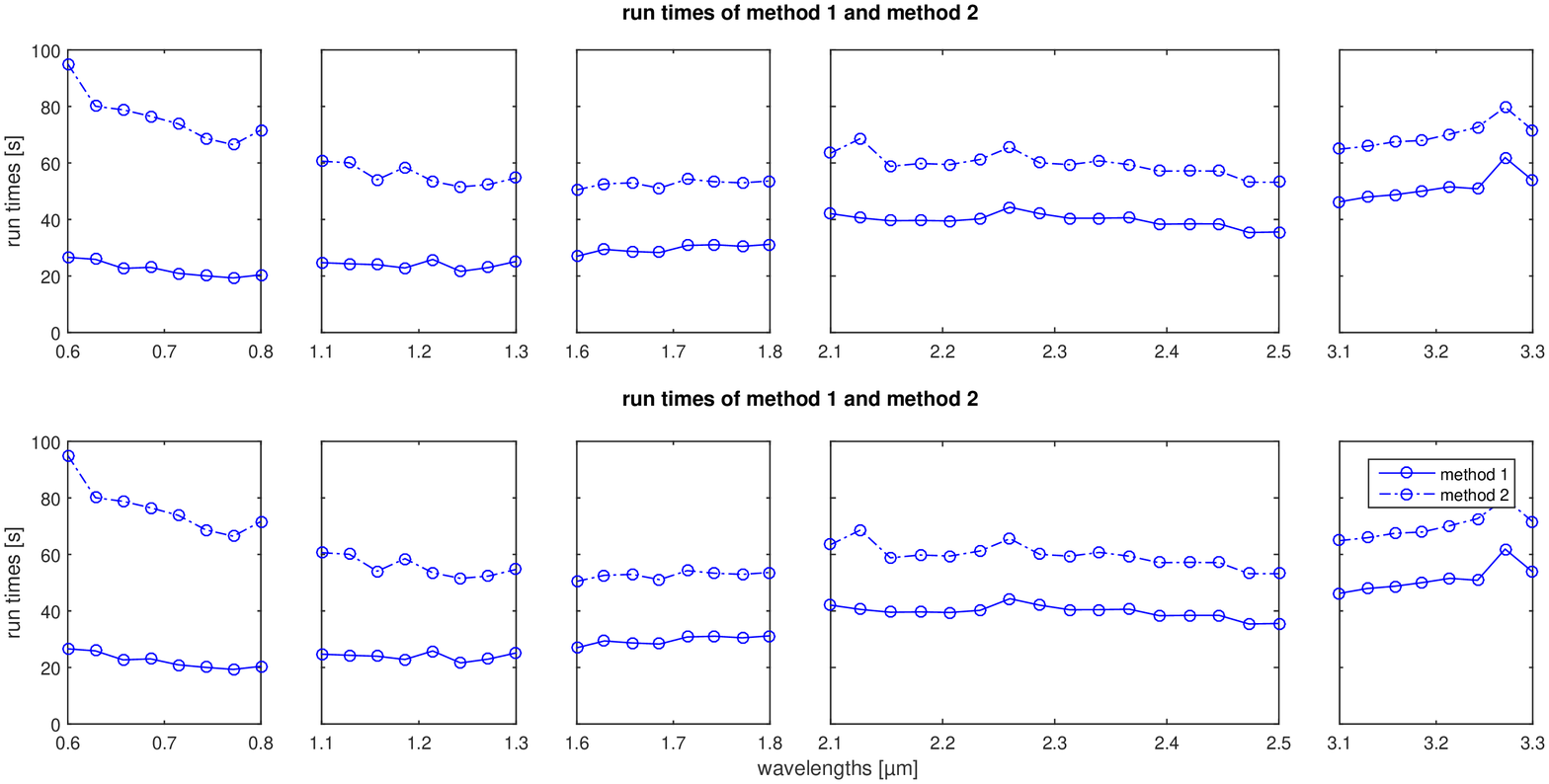}}
\includegraphics[clip,trim=0 0 0 {.5\wd0},width=1.0\textwidth]{RunTimesMethods1And2CsI.eps}
\endgroup
\label{RunTimesMethods1And2CsI}
\end{figure}

\subsubsection{Results for $\mathrm{H}_{2}$O}

\begin{figure}[h!]
\begingroup
\sbox0{\includegraphics[width =1.0\textwidth]{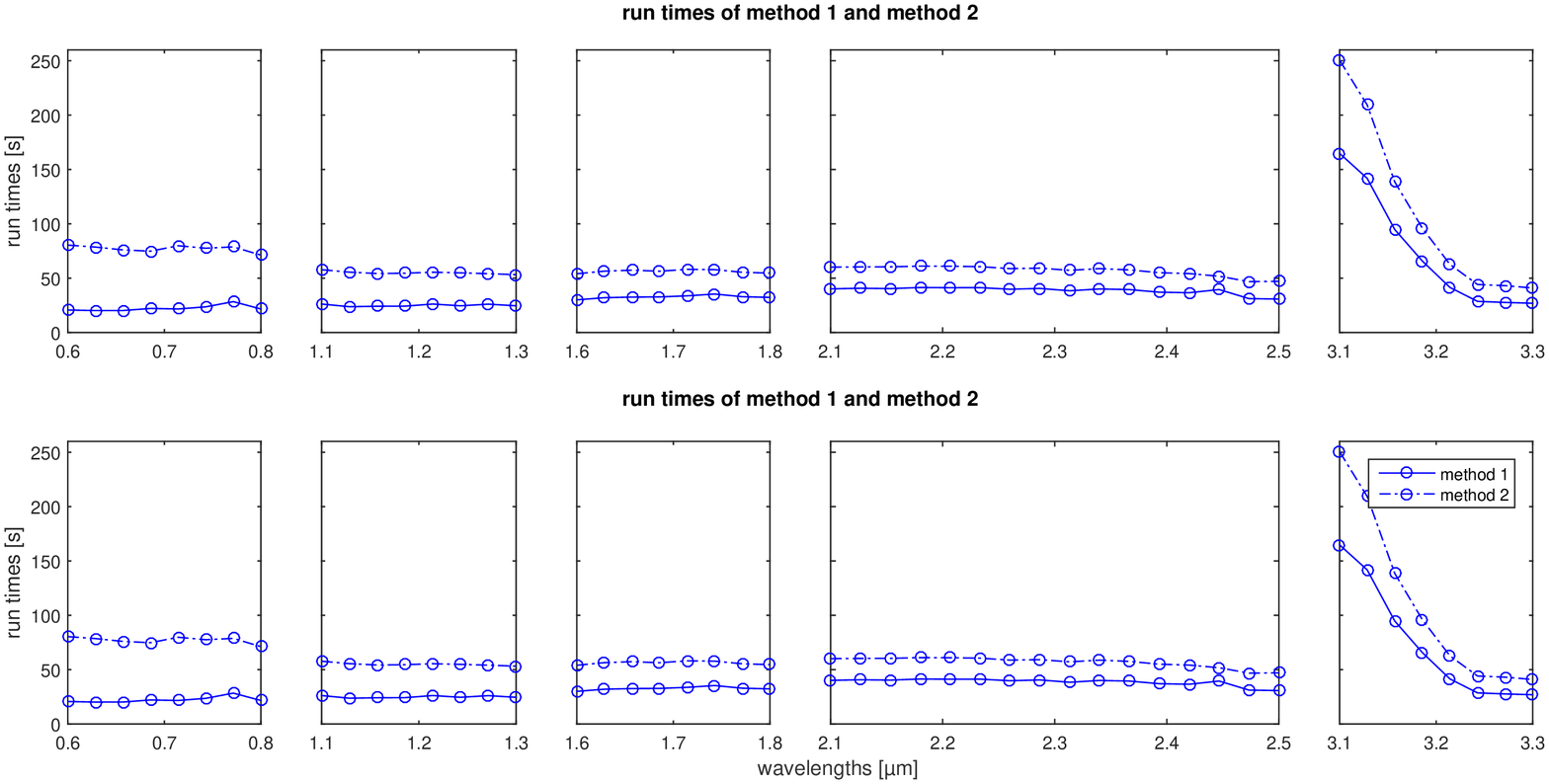}}
\includegraphics[clip,trim=0 0 0 {.5\wd0},width=1.0\textwidth]{RunTimesMethods1And2Water.eps}
\endgroup
\label{RunTimesMethods1And2Water}
\end{figure}

\subsection{Maximal Relative Deviations}

For the $10$ simulation runs we list the maximal relative deviations 
\begin{equation*}
100 \cdot \frac{\big\|\left(n_{part}^{1}(l_{j}), k_{part}^{1}(l_{j})\right)^{T} - \left(n_{part}^{2}(l_{j}), k_{part}^{2}(l_{j})\right)^{T}\big\|_{2}}{\big\|\left(n_{part}^{1}(l_{i}), k_{part}^{1}(l_{j})\right)^{T}\big\|_{2}}
\end{equation*}
of the refractive index reconstructions $\left(n_{part}^{2}(l_{j}), k_{part}^{2}(l_{j})\right)^{T}$ from method 2 from $\left(n_{part}^{1}(l_{j}), k_{part}^{1}(l_{j})\right)^{T}$  of method 1 for $j = 1, ..., 48$. At each wavelength, multiple local minima can be detected by both methods. For the relative deviations we always selected the local minima forming the smoothest reconstructions on each optical window in the sense of Section \ref{smoothest_reconstructions}.

\subsubsection{Results for Ag}

\begin{figure}[h!]
\begingroup
\sbox0{\includegraphics[width =1.0\textwidth]{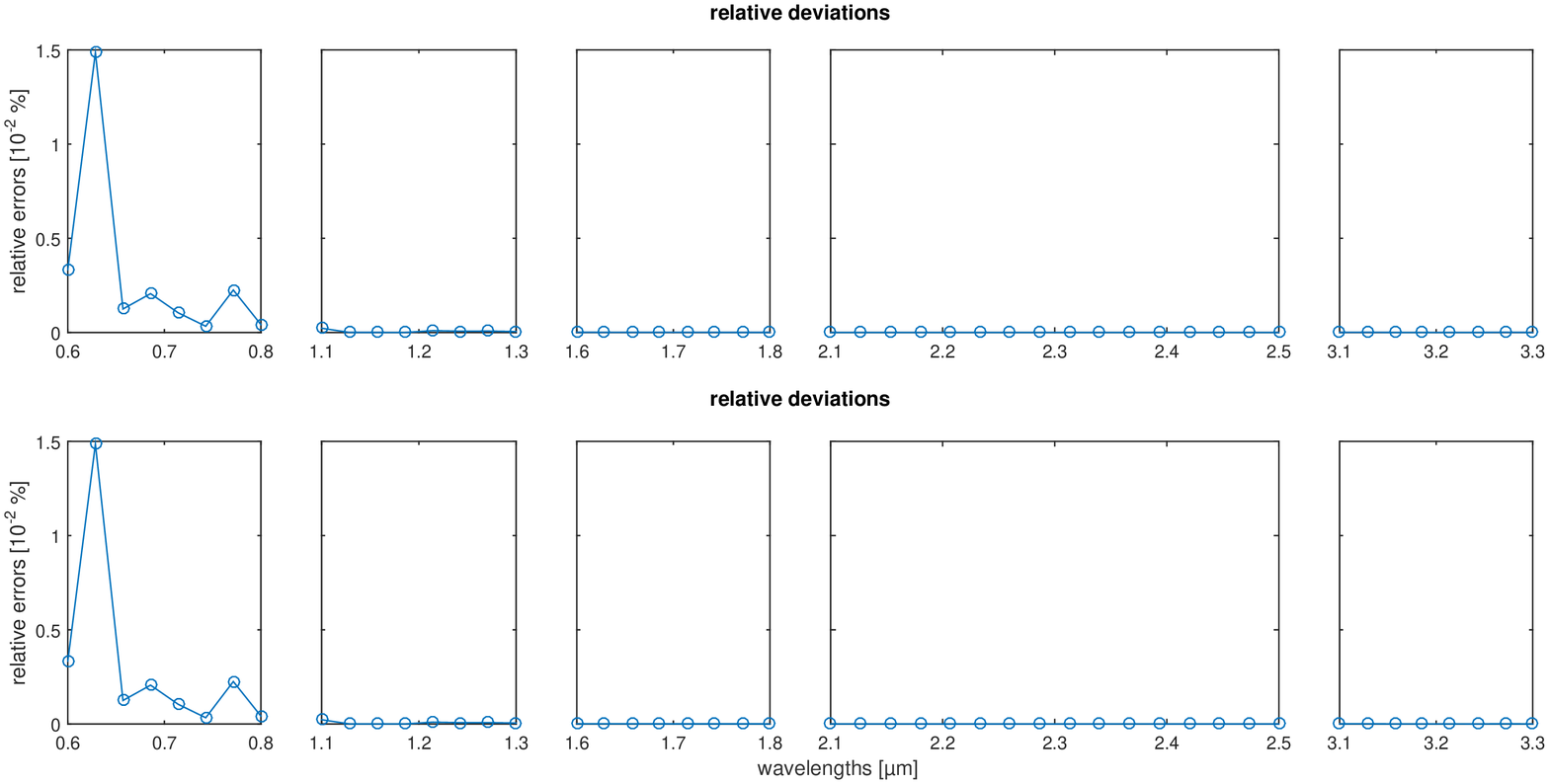}}
\includegraphics[clip,trim=0 0 0 {.5\wd0},width=1.0\textwidth]{RelDevSilver.eps}
\endgroup
\label{RelDevSilver}
\end{figure}

\newpage

\subsubsection{results for CsI}

\begin{figure}[h!]
\begingroup
\sbox0{\includegraphics[width =1.0\textwidth]{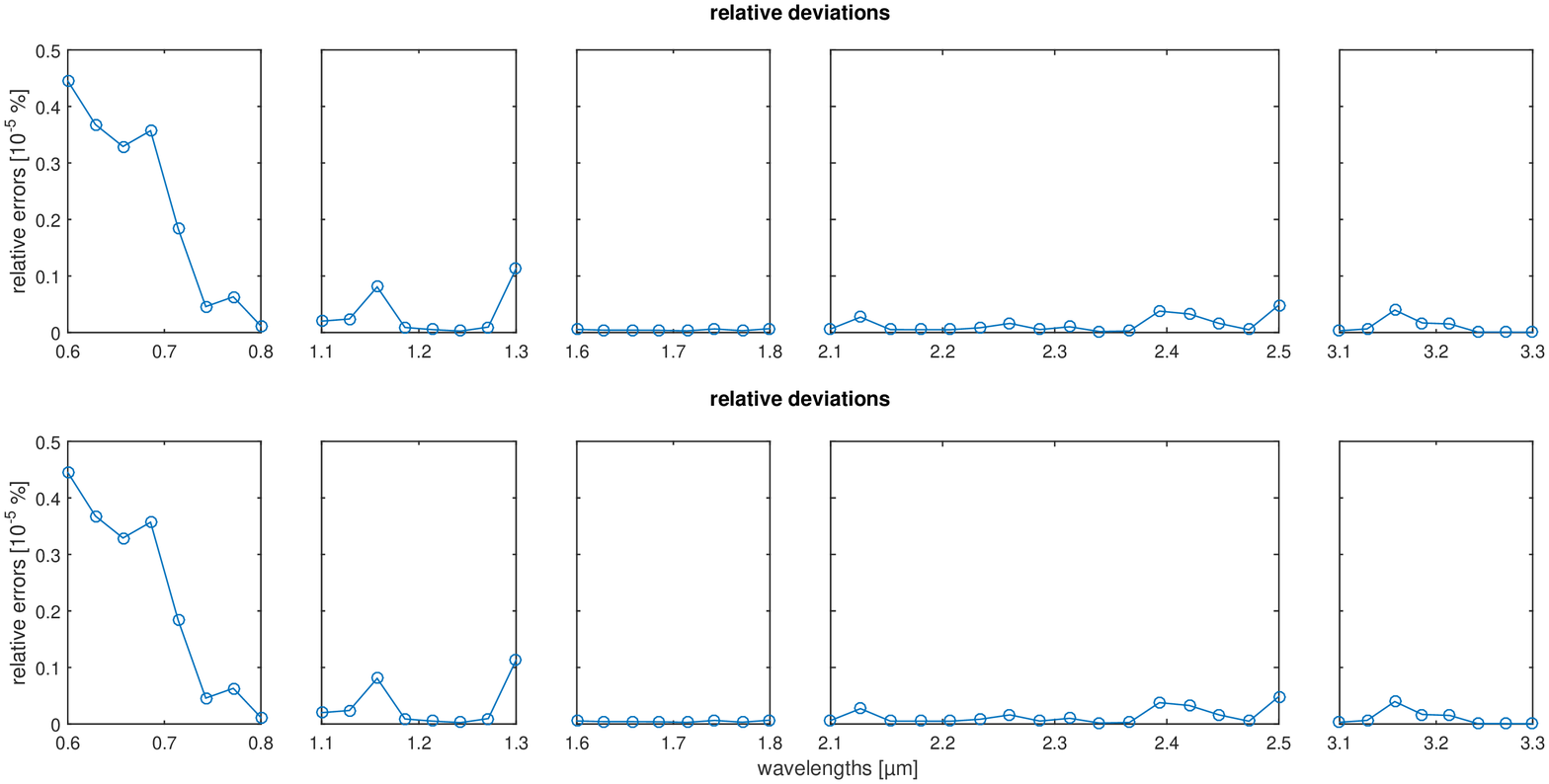}}
\includegraphics[clip,trim=0 0 0 {.5\wd0},width=1.0\textwidth]{RelDevCsI.eps}
\endgroup
\label{RelDevCsI}
\end{figure}

\subsubsection{Results for $\mathrm{H}_{2}$O}

\begin{figure}[h!]
\begingroup
\sbox0{\includegraphics[width =1.0\textwidth]{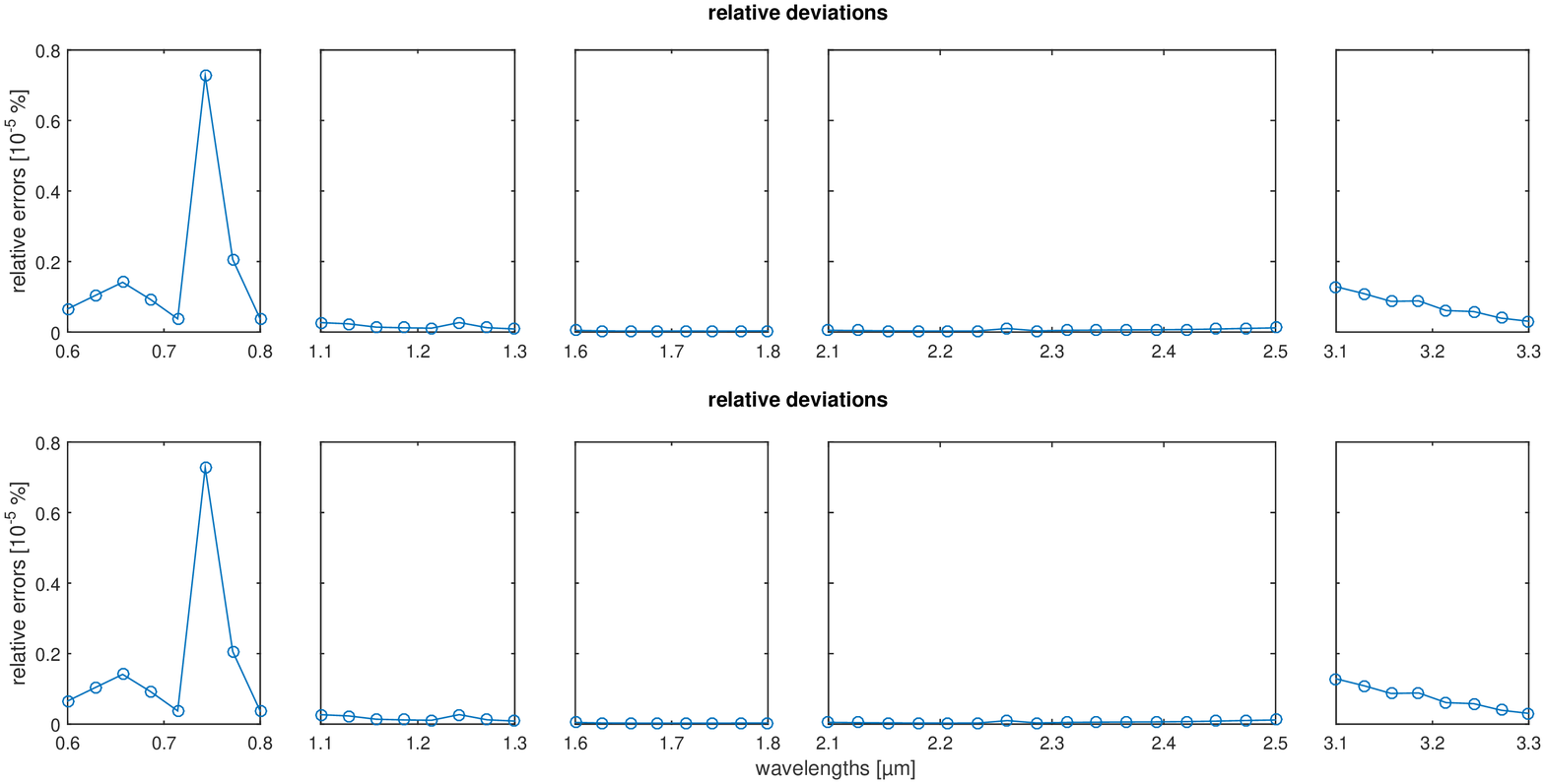}}
\includegraphics[clip,trim=0 0 0 {.5\wd0},width=1.0\textwidth]{RelDevWater.eps}
\endgroup
\label{RelDevWater}
\end{figure}

\subsection{Conclusion}

For Ag, the average total run time over all $48$ wavelengths for method 1 was $44.0167\%$ less than for method 2, for $\mathrm{H}_{2}$O $43.9808\%$ and for CsI $44.7322\%$, i.e. method 1 is almost two times faster than method 2. The results are of the same quality, since their relative deviations are just small fractions of percentages.

The continuation method approach saves run time significantly with the same quality of the results compared to using the truncation index \eqref{truncation_index} all the time. 

\section{Nonlinear Tikhonov Regularization}

So far we have solved the regression problem \eqref{nonlinreg_continuous} without any regularization, thus the obtained refractive index reconstructions might still be too error-contaminated to be of practical use. A widely used regularization strategy for nonlinear regression problems is Tikhonov regularization, which yields the regularized regression problem 
\begin{equation}
\label{nonlinreg_continuous_regularized}
\boldsymbol{x}_{\gamma} := \underset{\boldsymbol{x} \in \mathbb{R}^{D}}{\mathrm{argmin}} \|\boldsymbol{\Sigma}^{-\frac{1}{2}}(\boldsymbol{f}_{t}(\boldsymbol{x}) - \boldsymbol{f}(\boldsymbol{x}_{true}) - \boldsymbol{\delta})\|_{2}^{2} + \gamma \|\boldsymbol{x} - \boldsymbol{x}^{*}\|_{2}^{2} \quad \text{s.t.} \quad \boldsymbol{x} \in \Omega
\end{equation}
when we apply it on \eqref{nonlinreg_continuous}, cf. \cite{EKN89}. Here $\gamma$ is a regularization parameter and $\boldsymbol{x}^{*}$ is an estimate of the sought-after true solution. In many cases the unregularized problem has a whole set of minimizers, thus the vector $\boldsymbol{x}^{*}$ works also as a selection criterion. Now if a reasonable $\boldsymbol{x}^{*}$ is found, the regularization parameter $\gamma$ can be determined with the discrepancy principle, i.e. $\gamma$ is computed such that
\begin{align*}
\|\boldsymbol{\Sigma}^{-\frac{1}{2}}(\boldsymbol{f}_{t}(\boldsymbol{x}_{\gamma}) - \boldsymbol{f}(\boldsymbol{x}_{true}) - \boldsymbol{\delta})\|_{2} = R(\delta),
\end{align*} 
where $R(\delta)$ is an estimate of the residual of the ``true'' solution which depends on the noise level $\delta$. For this task monotonicity in the residual of $\boldsymbol{x}_{\gamma}$ is established in \cite{EKN89}.

The problem of finding a good estimate $\boldsymbol{x}^{*}$ still remains. In \cite{SEK93} an alternative implementable parameter choice strategy without the need of an $\boldsymbol{x}^{*}$ is derived. Applied on our problem it gives
\begin{align*}
& \gamma \big\langle \boldsymbol{\Sigma}^{-\frac{1}{2}}(\boldsymbol{f}_{t}(\boldsymbol{x}_{\gamma}) - \boldsymbol{f}(\boldsymbol{x}_{true}) - \boldsymbol{\delta}), \; \mathrm{J}_{\gamma}^{-1}(\boldsymbol{\Sigma}^{-\frac{1}{2}}(\boldsymbol{f}_{t}(\boldsymbol{x}_{\gamma}) - \boldsymbol{f}(\boldsymbol{x}_{true}) - \boldsymbol{\delta})\big\rangle = R(\delta), \\
\text{with} \quad & \mathrm{J}_{\gamma} := \gamma \mathrm{I} + \boldsymbol{\Sigma}^{-\frac{1}{2}}\mathrm{Jac}_{\boldsymbol{f}_{t}}(\boldsymbol{x}_{\gamma})\mathrm{Jac}_{\boldsymbol{f}_{t}}(\boldsymbol{x}_{\gamma})^{T}\boldsymbol{\Sigma}^{-\frac{1}{2}}.
\end{align*}
This method has the drawback that the matrix $\mathrm{J}_{\gamma}$ needs to be inverted, which may lead to instabilities.

Nevertheless the quality of the regularized solutions still depends strongly on the start values for solving \eqref{nonlinreg_continuous_regularized}. We know about our sought-after refractive indices that they form smooth curves on each of the five optical windows. The complex refractive index curves of most materials can be described using the so-called Lorentz-oscillator-model, cf. \cite{Qu10}. Here points with bigger curvatures only occur at so-called resonance frequencies corresponding to some isolated resonance wavelengths. Motivated by these facts we derive in the following a method to find reasonable start values for Phillips-Twomey-regularization out of the results of Algorithm \ref{reconstruction_algorithm}, which will be outlined in Section \ref{FurtherReg}.

\section{Finding the Smoothest Coupled Solutions}
\label{smoothest_reconstructions}

We have the problem of identifying the best approximation to the sought-after true particle material refractive index $\boldsymbol{x}_{true}$ out of a set of multiple solutions obtained with Algorithm \ref{reconstruction_algorithm}. This problem can be solved by coupling the solutions, which means that we combine solutions from neighboring wavelengths $l$ in each of the five optical windows in order to obtain a unique solution for every optical window. We know about the complex refractive index curves to be retrieved that they are smooth, hence we expect their sum of the squared second finite differences both in the real and imaginary parts to be small.

Let $l_{1}$, ..., $l_{s}$ denote the wavelengths of any of our five wavelength ranges. Let $N_{1}$, ..., $N_{s}$ be the number of solutions found for all the wavelengths. We denote with $\boldsymbol{x}^{i}_{j}$ the $j$-th solution found for wavelength $l_{i}$ for $i = 1, ..., s$ and $j = 1, ..., N_{i}$. Now we wish to find the smoothest combined solution from all possible combinations $\boldsymbol{x}^{1}_{j_{1}}$, ..., $\boldsymbol{x}^{s}_{j_{s}}$ for $j_{i} = 1, ..., N_{i}$, hence we have a total number of $\prod_{i=1}^{s}N_{i}$ combinations. Here we measure smoothness of a combination $\boldsymbol{x}^{1}_{j_{1}}$, ..., $\boldsymbol{x}^{s}_{j_{s}}$ with the sum 
\begin{equation*}
S := \sum_{i=2}^{s-1} \left(\left(\big(\boldsymbol{x}^{i-1}_{j_{i-1}}\big)_{1} - 2\big(\boldsymbol{x}^{i}_{j_{i}}\big)_{1} + \big(\boldsymbol{x}^{i+1}_{j_{i+1}}\big)_{1}\right)^{2} + \left(\big(\boldsymbol{x}^{i-1}_{j_{i-1}}\big)_{2} - 2\big(\boldsymbol{x}^{i}_{j_{i}}\big)_{2} + \big(\boldsymbol{x}^{i+1}_{j_{i+1}}\big)_{2}\right)^{2}\right)
\end{equation*}    
of its second finite differences both in the real parts $(\boldsymbol{x}^{i}_{j_{i}}\big)_{1}$ and its imaginary parts $(\boldsymbol{x}^{i}_{j_{i}}\big)_{2}$, which means that we regard a combination the smoother the smaller its sum $S$ is.

We encounter the problem that the total number of possible combinations $\prod_{i=1}^{s}N_{i}$ might get too big to iterate through all combinations in the search for the smoothest one in acceptable time. Therefore we propose a greedy algorithm, which uses each second finite difference as start point to find a smooth combination.

\begin{algorithm}[H]
\caption{Detection of the Smoothest Combination}
\label{greedy_algorithm}
\begin{algorithmic}[1]
\State $S_{min} = \infty$
\State $S_{cur} = 0$
\State $Comb = \{\}$
\State $SmoothestCombination = \{\}$
\For{$z = 2 \; \textbf{to} \; s-1$} \label{StartMainLoop} 
    \For{$c1 = 1 \; \textbf{to} \; N_{z - 1}$} \label{StartInnerLoops}
        \For{$c2 = 1 \; \textbf{to} \; N_{z}$}
            \For{$c3 = 1 \; \textbf{to} \; N_{z + 1}$} \label{EndInnerLoops}
	 \State $S_{cur} = \Big(\big(\boldsymbol{x}^{z-1}_{c1}\big)_{1} - 2\big(\boldsymbol{x}^{z}_{c2}\big)_{1} + \big(\boldsymbol{x}^{z+1}_{c3}\big)_{1}\Big)^{2} + \Big(\big(\boldsymbol{x}^{z-1}_{c1}\big)_{2} - 2\big(\boldsymbol{x}^{z}_{c2}\big)_{2} + \big(\boldsymbol{x}^{z+1}_{c3}\big)_{2}\Big)^{2}$ \label{StartInitialization}
            \State $Comb(z-1) = \boldsymbol{x}^{z-1}_{c1}$
            \State $Comb(z) = \boldsymbol{x}^{z}_{c2}$
     \algstore{myalg_b}
  \end{algorithmic}
\end{algorithm}  

	
\begin{algorithm}
  \begin{algorithmic}
      \algrestore{myalg_b}
            \State $Comb(z+1) = \boldsymbol{x}^{z+1}_{c3}$ \label{EndInitialization}
                 \For{$k = z - 2 \; \textbf{to} \; 1$} \label{StartLeftLoop}
	      \State $D_{min} = \infty$
	      \State $D_{cur} = \infty$
	      \State $\boldsymbol{x}_{min} = (0,0)^{T}$
	      \State $\boldsymbol{x}_{mid} = Comb(k+1)$
	      \State $\boldsymbol{x}_{right} = Comb(k+2)$
		\For{$j = 1 \; \textbf{to} \; N_{k}$} \label{StartNewMinSecDiff}
		\State $D_{cur} = \Big(\big(\boldsymbol{x}^{k}_{j}\big)_{1} - 2\big(\boldsymbol{x}_{mid})_{1} + \big(\boldsymbol{x}_{right}\big)_{1}\Big)^{2} \newline \hspace*{4.25cm} + \Big(\big(\boldsymbol{x}^{k}_{j}\big)_{2} - 2\big(\boldsymbol{x}_{mid}\big)_{2} + \big(\boldsymbol{x}_{right}\big)_{2}\Big)^{2}$ \vspace{0.1cm}
		     \If{$D_{cur} < D_{min}$}
		     \State $D_{min} = D_{cur}$
		     \State $\boldsymbol{x}_{min} = \boldsymbol{x}^{k}_{j}$
		     \EndIf
		\EndFor \label{EndNewMinSecDiff}
	      \State $S_{cur} = S_{cur} + D_{min}$
	      \State $Comb(k) = \boldsymbol{x}_{min}$	
	      \EndFor \label{EndLeftLoop}
                 \For{$k = z + 2 \; \textbf{to} \; s$} \label{StartRightLoop}
	      \State $D_{min} = \infty$
	      \State $D_{cur} = \infty$
	      \State $\boldsymbol{x}_{min} = (0,0)^{T}$
	      \State $\boldsymbol{x}_{mid} = Comb(k-1)$
	      \State $\boldsymbol{x}_{left} = Comb(k-2)$
		\For{$j = 1 \; \textbf{to} \; N_{k}$}
		\State $D_{cur} =  \Big(\big(\boldsymbol{x}_{left}\big)_{1} - 2\big(\boldsymbol{x}_{mid})_{1} + \big(\boldsymbol{x}^{k}_{j}\big)_{1}\Big)^{2} \newline  \hspace*{4.25cm} + \Big(\big(\boldsymbol{x}_{left}\big)_{2} - 2\big(\boldsymbol{x}_{mid}\big)_{2} + \big(\boldsymbol{x}^{k}_{j}\big)_{2}\Big)^{2}$
		     \If{$D_{cur} < D_{min}$}
		     \State $D_{min} = D_{cur}$
		     \State $\boldsymbol{x}_{min} = \boldsymbol{x}^{k}_{j}$
		     \EndIf
		\EndFor
	      \State $S_{cur} = S_{cur} + D_{min}$
	      \State $Comb(k) = \boldsymbol{x}_{min}$	
	      \EndFor \label{EndRightLoop}
            \If{$S_{cur} < S_{min}$} \label{StartSmoothest}
            \State $S_{min} = S_{cur}$
            \State $SmoothestCombination = Comb$
            \EndIf \label{EndSmoothest}
            \EndFor
        \EndFor
    \EndFor
\EndFor \label{EndMainLoop}
\end{algorithmic}
\end{algorithm}

The main loop spanning over the lines \ref{StartMainLoop} - \ref{EndMainLoop} iterates through all positions $z = 2, ..., s-1$ of a middle point for a second finite difference both in the real and imaginary part. At each position $z$ the inner loops beginning in lines \ref{StartInnerLoops} - \ref{EndInnerLoops} iterate through all possible second finite differences which can be formed out of the vectors $\boldsymbol{x}^{z-1}_{c1}$,  $\boldsymbol{x}^{z}_{c2}$ and $\boldsymbol{x}^{z+1}_{c3}$ for $c1 = 1, ..., N_{z-1}$, $c2 = 1, ..., N_{z}$ and $c3 = 1, ..., N_{z+1}$, i.e. they loop through all of their middle points and left and right neighbors at position $z$. In lines \ref{StartInitialization} - \ref{EndInitialization} the variable $S_{cur}$ is initialized with the sum of the squared second finite differences in the real and imaginary parts of the current vectors $\boldsymbol{x}^{z-1}_{c1}$,  $\boldsymbol{x}^{z}_{c2}$ and $\boldsymbol{x}^{z+1}_{c3}$ and the positions $z-1$, $z$ and $z+1$ of the array $Comb$ are filled with the current vectors. For $z \geq 3$ the loop in lines \ref{StartLeftLoop} - \ref{EndLeftLoop} successively fills the positions $k = z - 2, z - 3, ..., 1$ of the array $Comb$. At each new position $k$ the minimal sum of the two squared second finite differences in the real and imaginary part $D_{min}$ is determined in lines \ref{StartNewMinSecDiff} - \ref{EndNewMinSecDiff}, where the middle and right point are fixed and taken as the leftmost two vectors from the array $Comb$ and the right point runs through all $\boldsymbol{x}^{k}_{j}$ for $j = 1, ..., N_{k}$. After the vector $\boldsymbol{x}_{min}$ giving out of all of the $\boldsymbol{x}^{k}_{j}$'s giving the minimal sum $D_{min}$ is found, the $D_{min}$ is added to $S_{cur}$ and $\boldsymbol{x}_{min}$ is stored in $k$-th entry $Comb(k)$. In a similar way the loop in lines \ref{StartRightLoop} - \ref{EndRightLoop} succesively fills the positions $k = z + 2, z + 3, ..., s$ for $z \leq s - 2$. This time the left and middle point are fixed and taken as the rightmost two points of the array $Comb$, whereas the left point iterates through all $\boldsymbol{x}^{k}_{j}$ for $j = 1, ..., N_{k}$. Again the vector $\boldsymbol{x}_{min}$ is that one of the $\boldsymbol{x}^{k}_{j}$'s giving the minimal sum $D_{min}$ and it is stored in $Comb(k)$. As well the sum $D_{min}$ is added to $S_{cur}$.   

In the above procedure every triple of neighboring vectors from the results of Algorithm \ref{reconstruction_algorithm} is considered to possibly lie on the sought-after smoothest combination with the smallest sum of all squared second differences. The three vectors are used as start points to find a smooth combination with a greedy strategy, where only a vector is added to the current combination, if it gives the smallest sum $D_{min}$ at the left or right end of the growing set with vectors already added, until the first and last position are reached. 

Finally from all of the combinations constructed this way the smoothest one with the smallest sum $S_{min}$ out of all $S_{cur}$'s is selected in lines \ref{StartSmoothest} - \ref{EndSmoothest} to be the final output $SmoothestCombination$.

Define $N_{total} := \sum_{j=1}^{s}N_{j}$. Then the total number of operations needed for Algorithm \ref{greedy_algorithm} can be estimated by $\mathcal{O}\big(N_{total}\sum_{j=2}^{s-1}N_{j-1}N_{j}N_{j+1}\big)$ which is considerably less than the $\mathcal{O}\big(\prod_{j=1}^{s}N_{j}\big)$ operations needed by the naive method of iterating through all possible combinations.

\section{Further Regularization of Coupled Solutions}
\label{FurtherReg}

Not only for determining the smoothest refractive index curve reconstructions formed from the results of Algorithm \ref{reconstruction_algorithm} the coupled view on the solutions is beneficial - it also leads to further improvement of the results by Twomey-regularization. Let us investigate the coupled approach in a probability theoretical setting. Here we reuse the notations introduced in Section \ref{regression_with_truncated_series}, i.e. we let $\boldsymbol{x}^{1}$, ..., $\boldsymbol{x}^{s}$ denote a set of solution for any of the five optical windows. Then the joint posterior probability density of $\boldsymbol{x}^{1}$, ..., $\boldsymbol{x}^{s}$ is given by
\begin{equation}
\label{joint_density}
p(\boldsymbol{x}^{1}, ..., \boldsymbol{x}^{s} | \boldsymbol{e}^{1}, ..., \boldsymbol{e}^{s}) = \prod_{j=1}^{s}p(\boldsymbol{x}^{j} | \boldsymbol{e}^{j}) \propto \exp\left(-\textstyle{\frac{1}{2}}\sum_{j=1}^{s}\|\boldsymbol{\Sigma}_{j}^{-\frac{1}{2}}(\boldsymbol{f}_{t}^{j}(\boldsymbol{x}^{j}) - \boldsymbol{e}^{j})\|_{2}^{2}\right)\prod_{j=1}^{s}I_{\Omega}(\boldsymbol{x}^{j}),
\end{equation}
where $\boldsymbol{e}^{j}$ is the data vector for the $j$-th wavelength  $l_{j}$ having $N$ entries with $N$ being the size of the radius grid. Moreover $\boldsymbol{\Sigma}_{j}$ is the scaled covariance matrix for $l_{j}$ and $\boldsymbol{f}_{t}^{j}(\boldsymbol{x})$ is the applied model depending on $l_{j}$. Note that we initially have differing truncation indices $t_{1}$, ..., $t_{s}$. Since the coefficient functions of the truncated model function $\boldsymbol{f}_{t_{j}}(\boldsymbol{x})$ are decaying fast for each $t_{j}$, it is convenient to change to the same truncation index $t := \max\{t_{1}, ..., t_{s}\}$ for all wavelengths $l_{1}$, ..., $l_{s}$. The errors introduced by doing so are negliglible. It is easy to show that maximizing the joint density \eqref{joint_density} is equivalent to maximize all single densities $p(\boldsymbol{x}^{j} | \boldsymbol{e}^{j})$ independently, i.e. a joint MAP-estimator
\begin{equation*}
\boldsymbol{x}^{1}_{opt}, ..., \boldsymbol{x}^{s}_{opt} = \underset{\boldsymbol{x}^{1}, ..., \boldsymbol{x}^{s}}{\mathrm{argmax}}\;p(\boldsymbol{x}^{1}, ..., \boldsymbol{x}^{s} | \boldsymbol{e}^{1}, ..., \boldsymbol{e}^{s}) 
\end{equation*}
consists of the single MAP-estimators
\begin{equation*}
\boldsymbol{x}^{j}_{opt} = \underset{\boldsymbol{x}}{\mathrm{argmax}}\;p(\boldsymbol{x} | \boldsymbol{e}^{j})
\end{equation*}
for $j = 1, ..., s$. This means the the results of Algorithm \ref{reconstruction_algorithm} can be used to construct MAP-estimators for the joint posterior probabilty density. 

This behavior changes when we replace the joint prior probabilty density
\begin{equation*}
p_{prior}(\boldsymbol{x}^{1}, ..., \boldsymbol{x}^{s}) = \left(\mathrm{vol}(\Omega)\right)^{-s}\prod_{j=1}^{s}I_{\Omega}(\boldsymbol{x}^{j})
\end{equation*}
with
\begin{align*}
p_{prior}(\boldsymbol{x}^{1}, ..., \boldsymbol{x}^{s}) \propto \exp\left(-\textstyle{\frac{1}{2}}\gamma S(\boldsymbol{x}^{1}, ..., \boldsymbol{x}^{s})\right)\prod_{j=1}^{s}I_{\Omega}(\boldsymbol{x}^{j}), 
\end{align*}
where
\begin{align*}
S(\boldsymbol{x}^{1}, ..., \boldsymbol{x}^{s}) :=  & \sum_{i=2}^{s-1} \left(\Big(\big(\boldsymbol{x}^{i-1}\big)_{1} - 2\big(\boldsymbol{x}^{i}\big)_{1} + \big(\boldsymbol{x}^{i+1}\big)_{1}\Big)^{2} + \Big(\big(\boldsymbol{x}^{i-1}\big)_{2} - 2\big(\boldsymbol{x}^{i}\big)_{2} + \big(\boldsymbol{x}^{i+1}\big)_{2}\Big)^{2}\right) \\
+ \rho & \sum_{i=1}^{s}\Big(\big(\boldsymbol{x}^{i}\big)_{1}^{2} + \big(\boldsymbol{x}^{i}\big)_{2}^{2}\Big),
\end{align*} 
where $\gamma$ is a regularization parameter and $\rho$ is a parameter specifying the amount Tikhonov regularization.

In the new prior distribution we use a combination of Tikhonov and Phillips-Twomey-regularization both in the real and imaginary parts. Here we apply a small amount of Tikhonov-regularization by setting $\rho = 10^{-8}$, such that the resulting regularization operator gets regular. This means that the regularized regression problem \eqref{Twomey_regularization} can be transformed into standard Tikhonov form and that the monotonicity results from \cite{EKN89} are still valid. Each second finite difference is clearly a function of three neighboring points, therefore a decoupled computation of the joint MAP-estimator
\begin{equation}
\label{Twomey_regularization}
\boldsymbol{x}^{1}_{opt}, ..., \boldsymbol{x}^{s}_{opt} := \underset{\boldsymbol{x}^{1}, ..., \boldsymbol{x}^{s}}{\mathrm{argmax}}\;p(\boldsymbol{x}^{1}, ..., \boldsymbol{x}^{s} | \boldsymbol{e}^{1}, ..., \boldsymbol{e}^{s}) 
\end{equation}
with
\begin{equation*}
p(\boldsymbol{x}^{1}, ..., \boldsymbol{x}^{s} | \boldsymbol{e}^{1}, ..., \boldsymbol{e}^{s}) \propto \exp\left(-\textstyle{\frac{1}{2}}\sum_{j=1}^{s}\|\boldsymbol{\Sigma}_{j}^{-\frac{1}{2}}(\boldsymbol{f}_{t}^{j}(\boldsymbol{x}^{j}) - \boldsymbol{e}^{j})\|_{2}^{2}-\textstyle{\frac{1}{2}}\gamma S(\boldsymbol{x}^{1}, ..., \boldsymbol{x}^{s})\right)\prod_{j=1}^{s}I_{\Omega}(\boldsymbol{x}^{j})
\end{equation*}
for each wavelength seperately is not possible anymore after changing to the new prior density. However the result vectors $\boldsymbol{x}^{1}$, ..., $\boldsymbol{x}^{s}$ from Algorithm \ref{greedy_algorithm} form a good start vector to solve the nonlinear regression problem \eqref{Twomey_regularization}.

We selected the regularization parameter $\gamma$ using the discrepancy principle, i.e. we compute $\gamma$ such that the regularized solution 
\begin{align*}
\boldsymbol{x}^{1}_{\gamma}, ..., \boldsymbol{x}^{s}_{\gamma} := \underset{\boldsymbol{x}^{1}, ..., \boldsymbol{x}^{s}}{\mathrm{argmin}}\;\sum_{j=1}^{s}\|\boldsymbol{\Sigma}_{j}^{-\frac{1}{2}}(\boldsymbol{f}_{t}^{j}(\boldsymbol{x}^{j}) - \boldsymbol{e}^{j})\|_{2}^{2} + \gamma S(\boldsymbol{x}^{1}, ..., \boldsymbol{x}^{s}) \quad \text{s.t.} \quad \boldsymbol{x}^{j} \in \Omega, \quad j = 1, ..., s
\end{align*}
fulfills a relation of the form
\begin{align*}
\sum_{j=1}^{s}\|\boldsymbol{\Sigma}_{j}^{-\frac{1}{2}}(\boldsymbol{f}_{t}^{j}(\boldsymbol{x}^{j}_{\gamma}) - \boldsymbol{e}^{j})\|_{2}^{2} = R(\delta),
\end{align*}
where $R(\delta)$ is a proposed residual value depending on the noise level $\delta$. In \cite{AK16}  several different residual values are proposed for a fixed model discretization and a set of regularization parameters is obtained from those using the discrepancy principle. The pairings of model discretizations and regularization parameters obtained this way are compared by their Bayesian posterior probabilities. For these probabilities a set of integrals of the different model posterior densities is needed to be computed, which can be done with Monte Carlo integration methods. Due to the highly nonlinear behavior of our model $\boldsymbol{f}_{t}(\boldsymbol{x})$ such integration methods are not available here. Therefore we simplified the posterior exploration in such a way that only one residual value is proposed. 

Since each observed probability density $p(\boldsymbol{e}^{j} | \boldsymbol{x}^{j})$ for $j = 1, ..., s$ is Gaussian, the joint observed density $p(\boldsymbol{e}^{1}, ..., \boldsymbol{e}^{1} | \boldsymbol{x}^{1}, ..., \boldsymbol{x}^{s}) = \prod_{j=1}^{s}p(\boldsymbol{e}^{j} | \boldsymbol{x}^{j})$ is Gaussian as well. We have $\boldsymbol{x}^{j} \in \mathbb{R}^{2}$, thus the sum of residuals $\sum_{j=1}^{s}\|\boldsymbol{\Sigma}_{j}^{-\frac{1}{2}}(\boldsymbol{f}_{t}^{j}(\boldsymbol{x}^{j}) - \boldsymbol{e}^{j})\|_{2}^{2}$ running through all wavelengths in the optical window is $\chi^{2}(2s)$-distributed. This yields
\begin{align*}
\mathbb{E}\Big(\sum_{j=1}^{s}\|\boldsymbol{\Sigma}_{j}^{-\frac{1}{2}}(\boldsymbol{f}_{t}^{j}(\boldsymbol{x}^{j}) - \boldsymbol{e}^{j})\|_{2}^{2}\Big) = 2s.
\end{align*}
Now a widely proposed residual value for the discrepancy principle is $\tau \cdot 2s$, where $\tau = 1.1$ is the so-called Morozov safety factor. This choice is prone to under- or overregularization since the residual value corresponding to the ``true'' solution might be much smaller or bigger than $2\tau s$. Therefore we proposed a residual value which depends more dynamically on the observed behovior of the residual.

Let $\boldsymbol{x}^{1}_{0}$, ..., $\boldsymbol{x}^{s}_{0}$ denote the unregularized solutions, i.e. the results of Algorithm \ref{greedy_algorithm}. Then their squared residual is given by $R_{0} := \sum_{j=1}^{s}\|\boldsymbol{\Sigma}_{j}^{-\frac{1}{2}}(\boldsymbol{f}_{t}^{j}(\boldsymbol{x}^{j}_{0}) - \boldsymbol{e}^{j})\|_{2}^{2}$. We first proposed
\begin{equation*}
R(\delta) = \max\left\{2\tau_{1} s, \; \tau R_{0}\right\},
\end{equation*}
where we selected $\tau_{1} = 1.1$. This means that the residual of the regularized solution is beginning at $R_{0}$ at least increased by the factor $\tau_{1}$, which avoids underregularization. If it then happens that $\frac{R(\delta)}{R_{0}} > \theta$ with $\theta = 1.5$ the proposed residual is most likely too big and overregularization occurs. In this case we corrected $R(\delta)$ by setting
\begin{equation*}
R(\delta) = \max\left\{2 \tau_{2} s, \; \theta R_{0}\right\}
\end{equation*}
with $\tau_{2} = 0.9$.

\section{Numerical Results}
\label{ReconRefrac}

To see how reliable our proposed reconstruction algortihm is, we performed for each of the scatterer materials Ag, $\mathrm{H}_{2}$O a numerical study with $100$ sweeps through all $48$ wavelenghts of the five optical windows with the same settings as in Section \ref{comparison_with_trunc_ind}. We found out that the radii $r_{1} := 0.1 \; \mu$m,  $r_{2} := 0.2 \; \mu$m and $r_{3} := 0.3 \; \mu$m contain the most information about the refractive indices. This was found by keeping our $48$ wavelengths fixed and comparing the quality of inversion results under varying aerosol particle radii. Bigger radii did not improve the results in our simulations and refractive index reconstructions only using bigger radii even turned out to be too unstable. A more thorough treatment of this problem can be found in \cite{EHR11}, where a covariance eigenvalue analysis is used. Although not directly comparable with our study of uncoated particles, the coated radii $0.0975 \; \mu$m, $0.2305 \; \mu$m and $0.11 \; \mu$m carrying the most information content found in this study are roughly comparable to our radii.

We computed original spectral extinctions
\begin{equation*}
\left(\boldsymbol{e}_{true}\right)_{i,j} := \pi r_{i}^{2}\sum_{n = 1}^{N_{tr}}q_{n}(m_{med}(l_{i}), m_{part}(l_{i}), r_{j}, l_{i}), \quad i = 1, ..., 48, \quad j = 1, ..., 3 
\end{equation*}
for Ag, $\mathrm{H}_{2}$O and CsI and added zero-mean Gaussian noise to it in order to obtain the simulated noisy spectral extinctions
\begin{equation*}
\left(\boldsymbol{e}\right)_{i,j} = \left(\boldsymbol{e}_{true}\right)_{i,j} + \delta_{i,j} \quad \text{ with } \; \delta_{i,j} \sim \mathcal{N}(0,(0.05 \cdot \left(\boldsymbol{e}_{true}\right)_{i,j})^{2}), \quad i = 1, ..., 48, \quad j = 1, ..., 3.
\end{equation*}
The standard deviations were taken to be $5\%$ of the true spectral extinctions. Real experiments using $500$ wavelengths were contaminated by Gaussian noise with $30\%$ of the true spectral extinctions as standard deviations. We expect that switching to $48$ wavelengths and thus increasing the time resolution of the measurements will lower the standard deviations to a small percentage. We used a sample size of $N_{s} = 300$ to compute each mean $\left(\boldsymbol{e}_{real}\right)_{i,j}$ of noisy spectral extinctions. 

In the following the results are presented separately for each of the three materials. The uppermost plot presents the relative errors of the unregularized solutions obtained with Algorithm \ref{greedy_algorithm} from the original scatterer refractive indices. The next plot displays the run times of Algorithm \ref{reconstruction_algorithm}, which returned all local minima of \eqref{nonlinreg_continuous}. These candidate solutions served as input for Algorithm \ref{greedy_algorithm}. Then the relative errors of the regularized solutions are presented. Finally the relative errors of the average of the regularized solutions are shown.

\subsection{Results for Ag}

\subsubsection{Results of Algorithms \ref{reconstruction_algorithm} and \ref{greedy_algorithm}}

\begin{figure}[h!]
\centering
\includegraphics[width =1.0\textwidth]{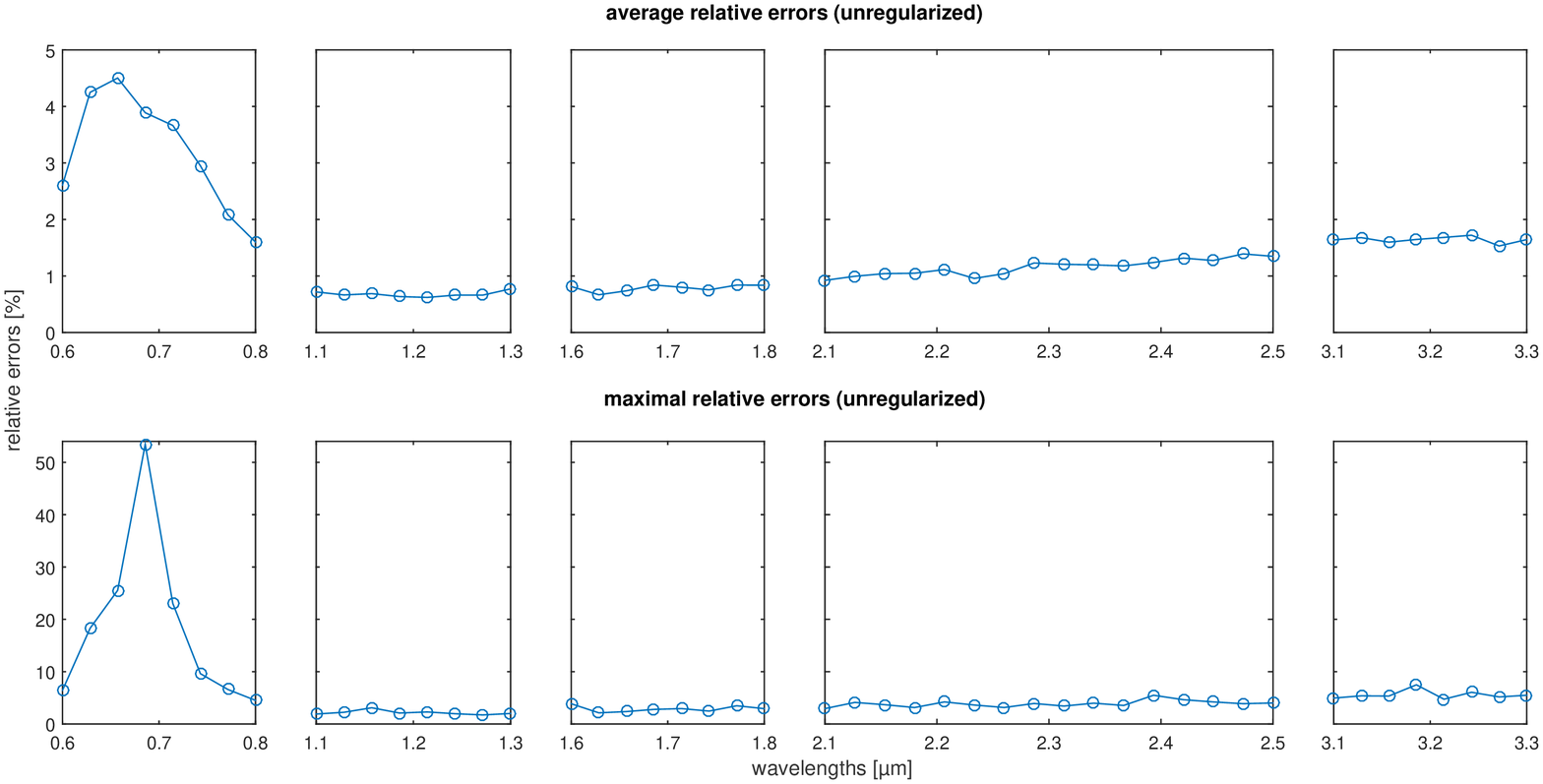}
\label{ErrUnregAvgMaxAg}
\end{figure}

\begin{figure}[h!]
\centering
\includegraphics[width =1.0\textwidth]{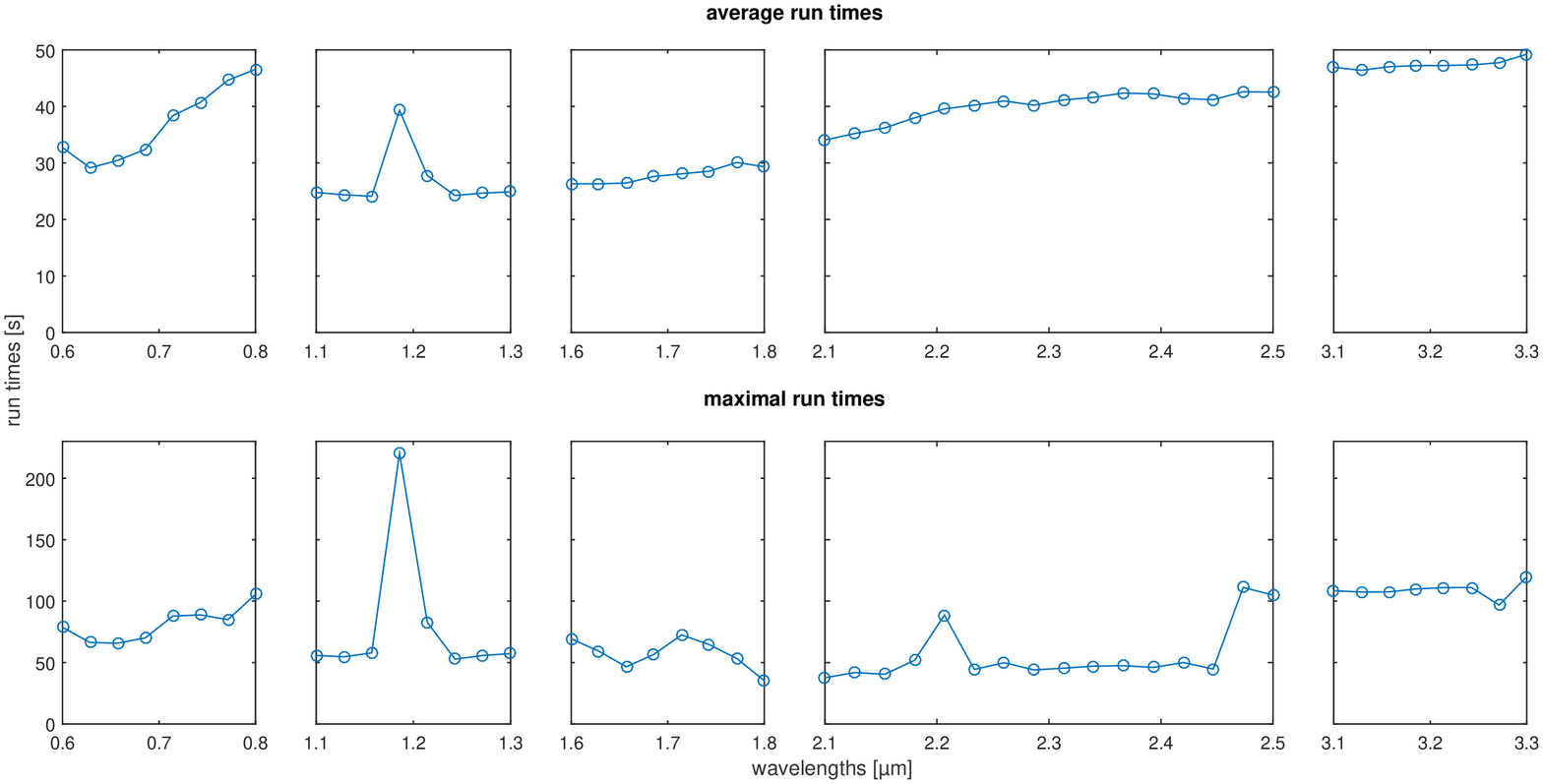}
\label{RunTimesAvgMaxAg}
\end{figure}

\newpage

\subsubsection{Relative Errors of the Regularized Solutions}

\begin{figure}[h!]
\centering
\includegraphics[width =1.0\textwidth]{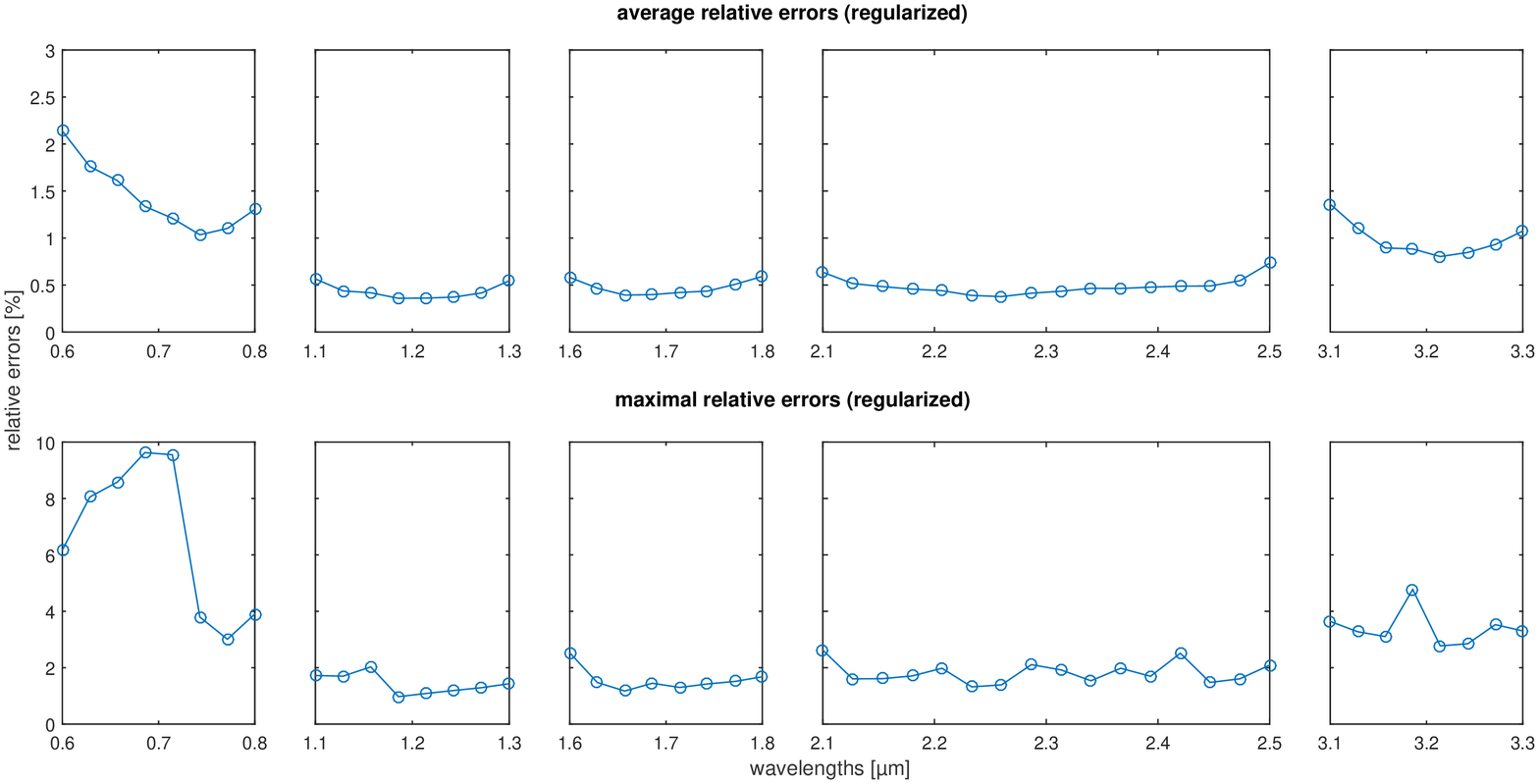}
\label{ErrRegAvgMaxAg}
\end{figure}

\subsubsection{Relative Errors of the Average of the Regularized Solutions}

\begin{figure}[h!]
\begingroup
\sbox0{\includegraphics[width =1.0\textwidth]{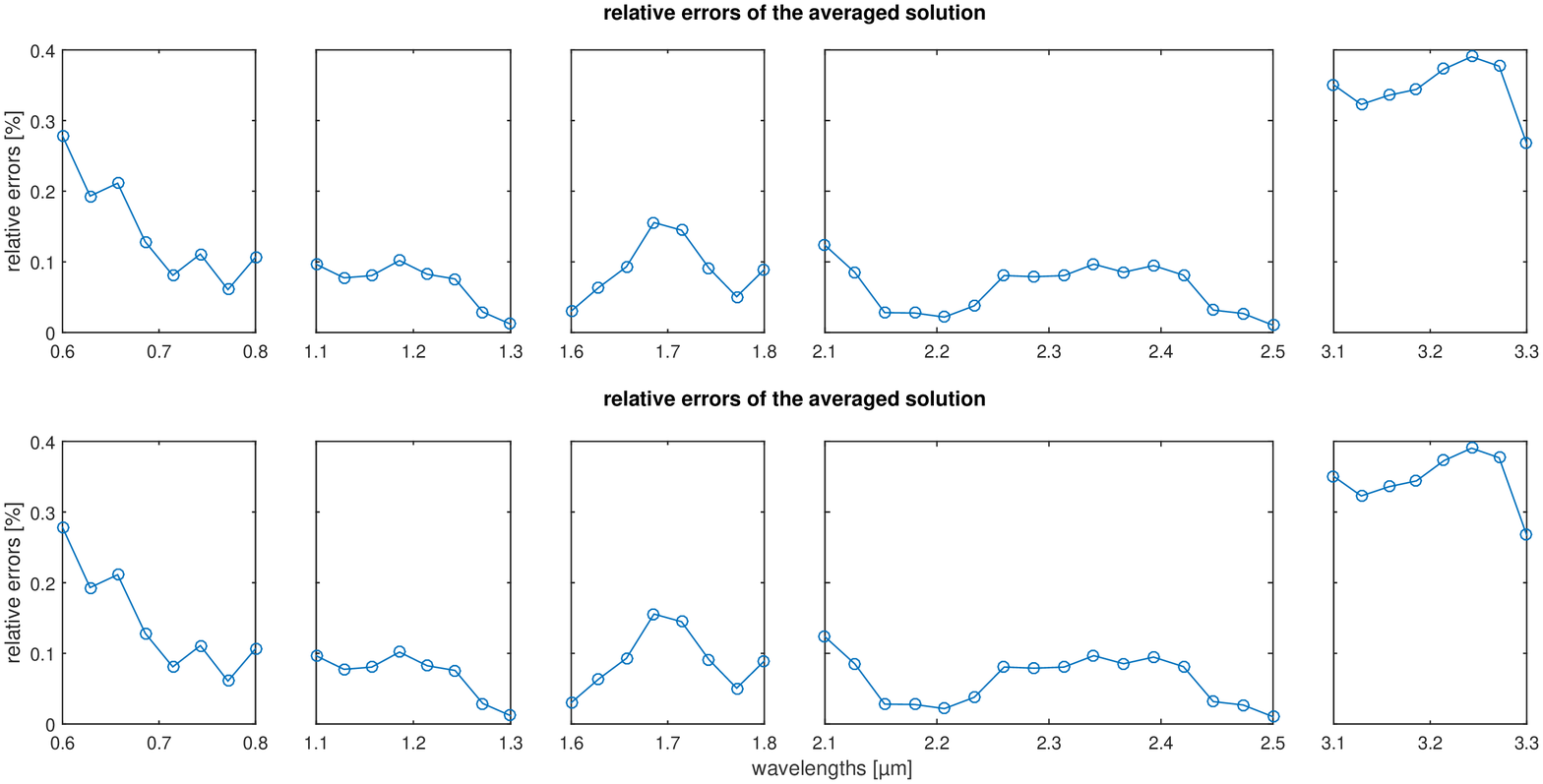}}
\includegraphics[clip,trim=0 0 0 {.5\wd0},width=1.0\textwidth]{ErrorsAveragedSolutionSilver.eps}
\endgroup
\label{ErrRegAvgAg}
\end{figure}

\newpage

\subsection{Results for CsI}

\subsubsection{Results of Algorithms \ref{reconstruction_algorithm} and \ref{greedy_algorithm}}

\begin{figure}[h!]
\centering
\includegraphics[width =1.0\textwidth]{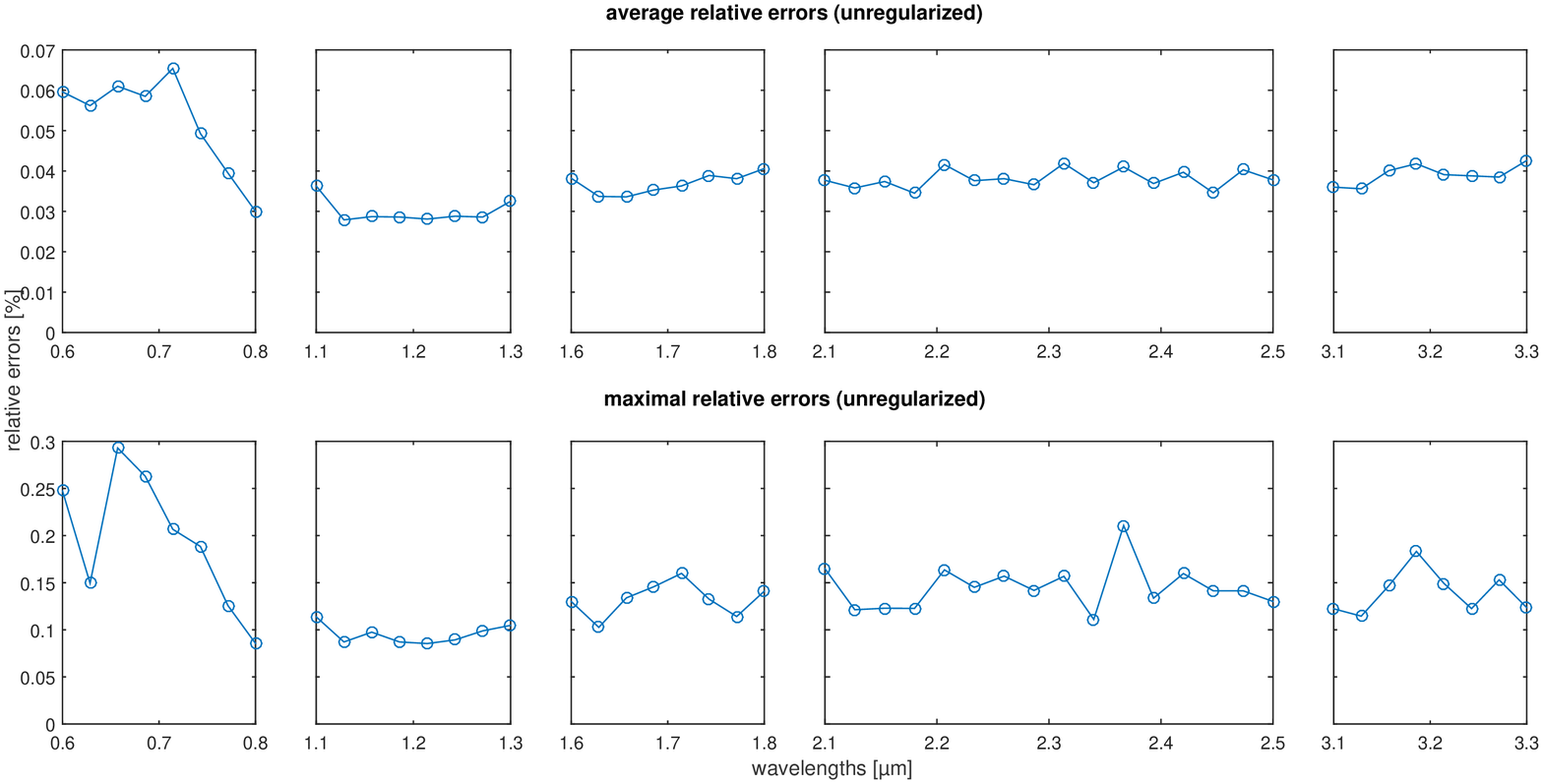}
\label{ErrUnregAvgMaxCsI}
\end{figure}

\begin{figure}[h!]
\centering
\includegraphics[width =1.0\textwidth]{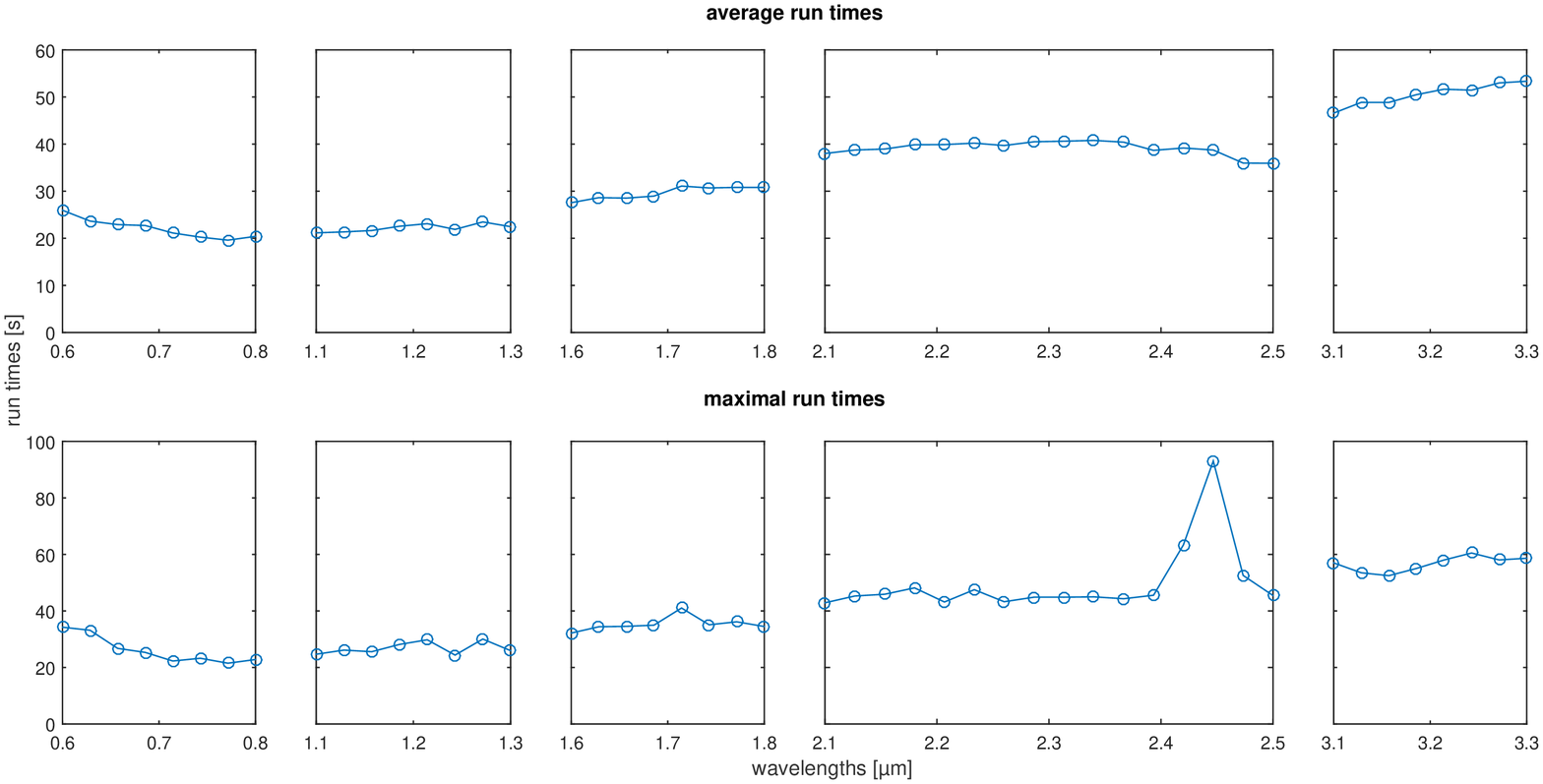}
\label{RunTimesAvgMaxCsI}
\end{figure}

\newpage

\subsubsection{Relative Errors of the Regularized Solutions}

\begin{figure}[h!]
\centering
\includegraphics[width =1.0\textwidth]{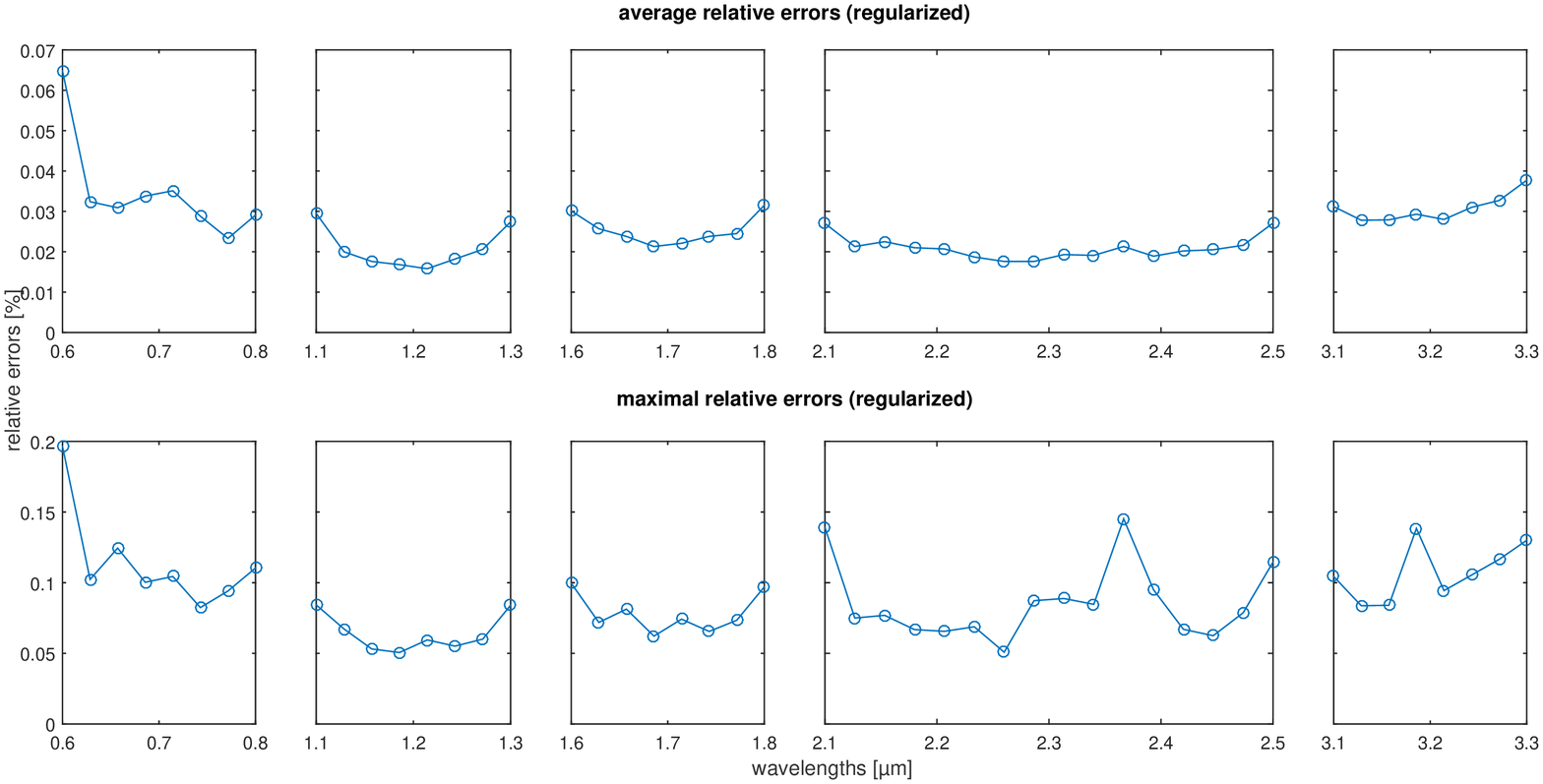}
\label{ErrRegAvgMaxCsI}
\end{figure}

\subsubsection{Relative Errors of the Average of the Regularized Solutions}

\begin{figure}[h!]
\begingroup
\sbox0{\includegraphics[width =1.0\textwidth]{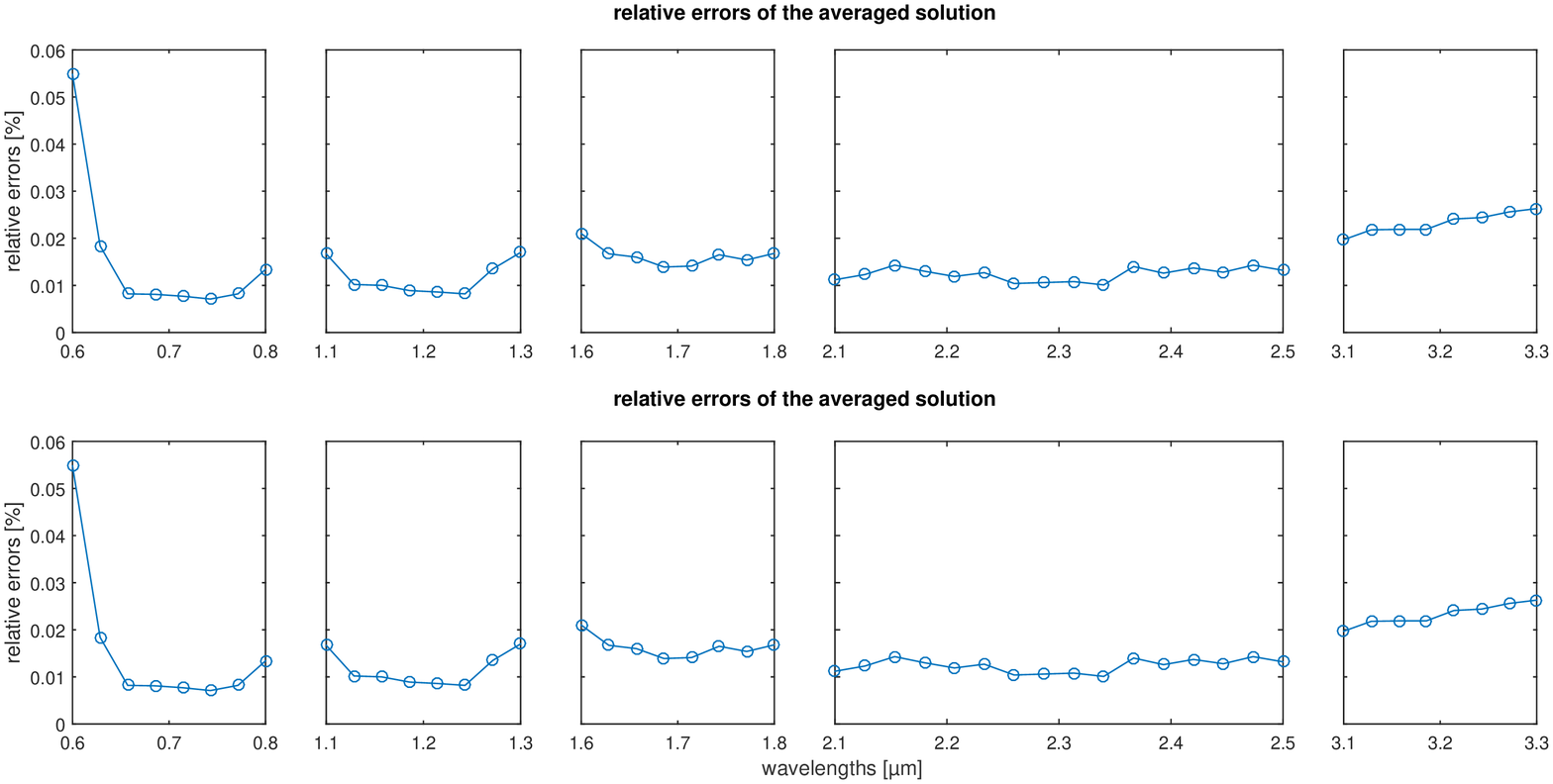}}
\includegraphics[clip,trim=0 0 0 {.5\wd0},width=1.0\textwidth]{ErrorsAveragedSolutionCsI.eps}
\endgroup
\label{ErrRegAvgCsI}
\end{figure}

\newpage

\subsection{Results for $\mathrm{H}_{2}$O}

\subsubsection{Results of Algorithms \ref{reconstruction_algorithm} and \ref{greedy_algorithm}}

\begin{figure}[h!]
\centering
\includegraphics[width =1.0\textwidth]{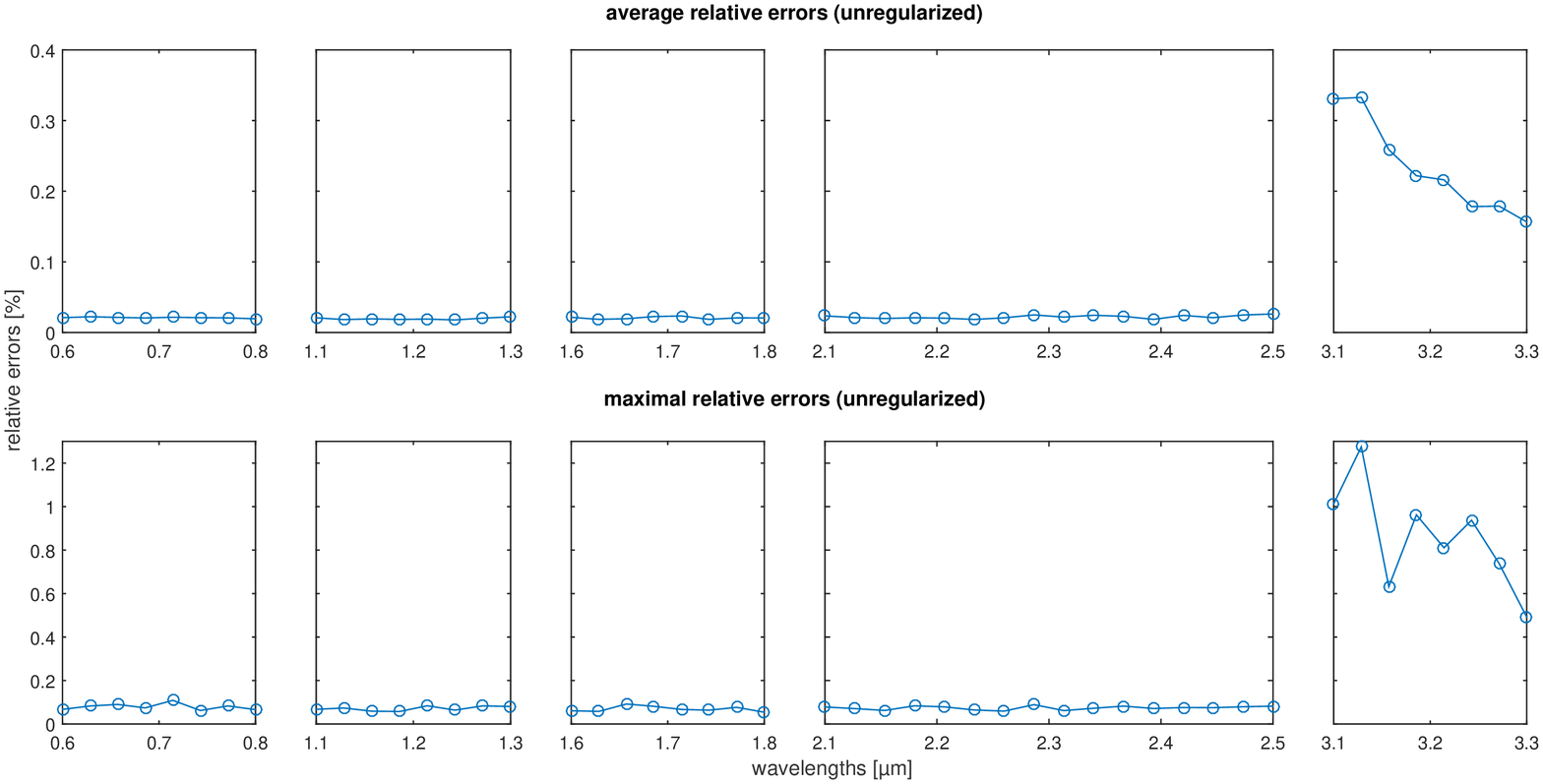}
\label{ErrUnregAvgMaxWater}
\end{figure}

\begin{figure}[h!]
\centering
\includegraphics[width =1.0\textwidth]{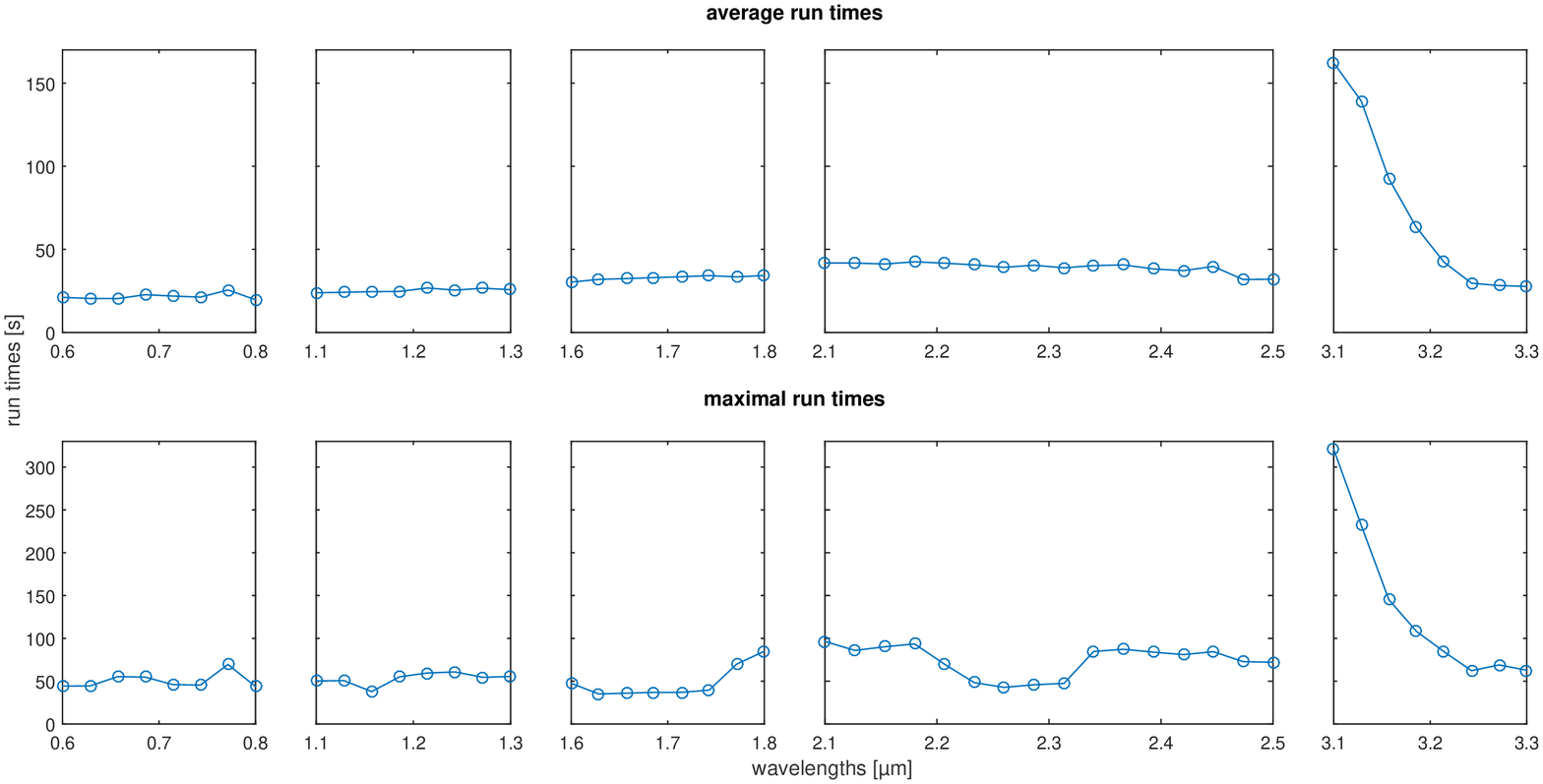}
\label{RunTimesAvgMaxWater}
\end{figure}

\newpage

\subsubsection{Relative Errors of the Regularized Solutions}

\begin{figure}[h!]
\centering
\includegraphics[width =1.0\textwidth]{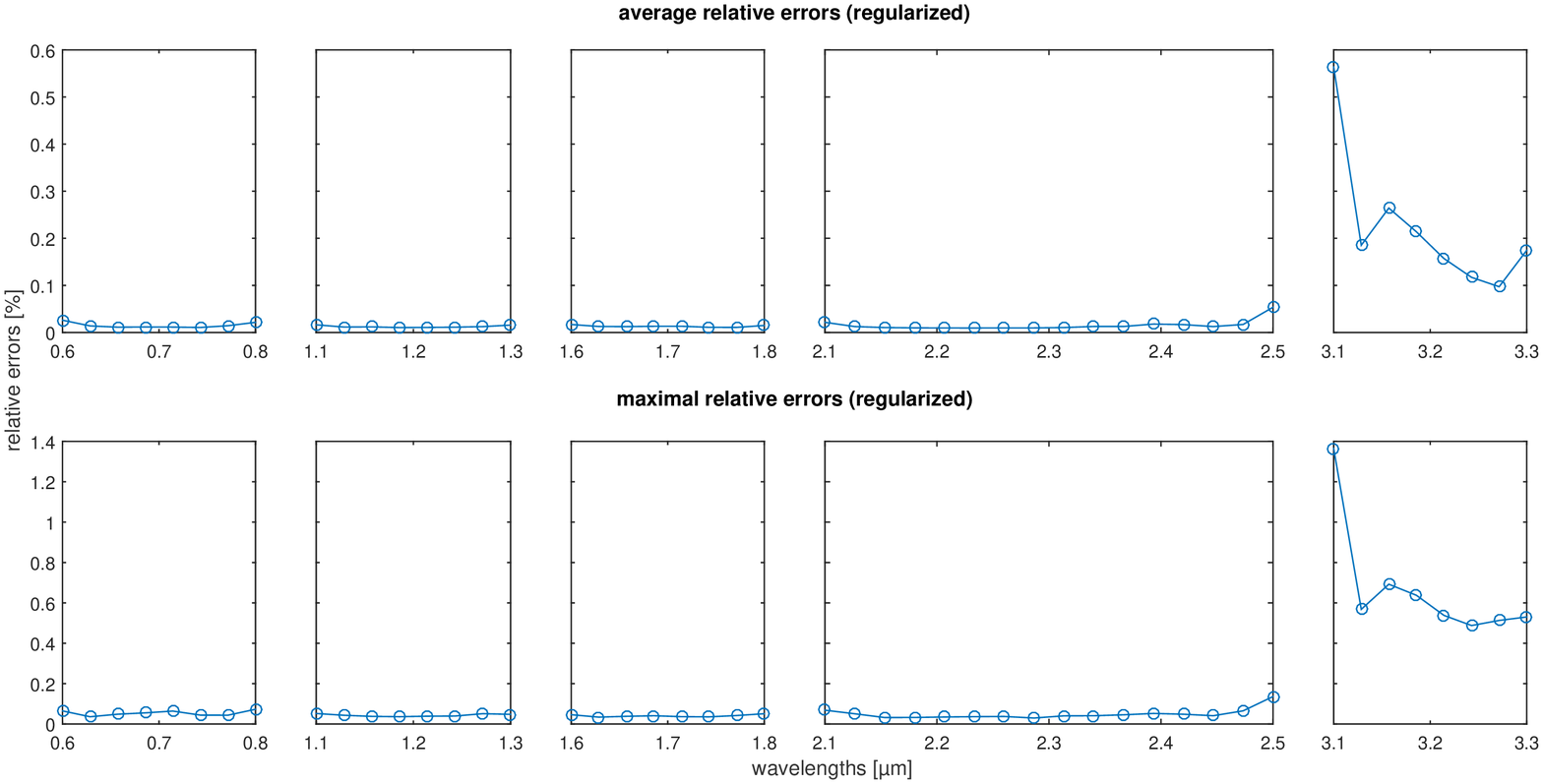}
\label{ErrRegAvgMaxWater}
\end{figure}

\subsubsection{Relative Errors of the Average of the Regularized Solutions}

\begin{figure}[h!]
\begingroup
\sbox0{\includegraphics[width =1.0\textwidth]{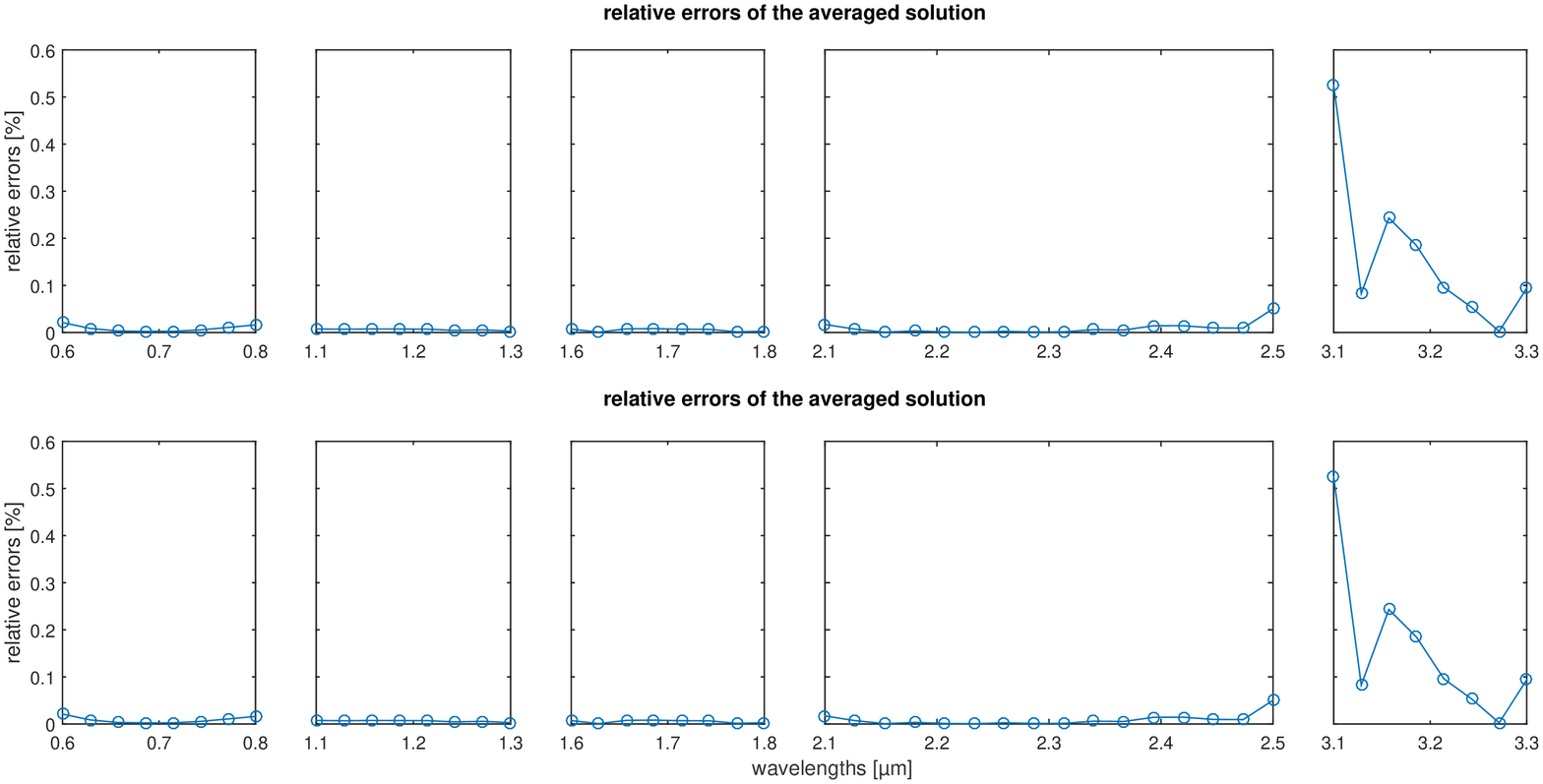}}
\includegraphics[clip,trim=0 0 0 {.5\wd0},width=1.0\textwidth]{ErrorsAveragedSolutionWater.eps}
\endgroup
\label{ErrRegAvgWater}
\end{figure}

\subsection{Conclusion}

The severest relative errors can be observed for Ag. For the initial unregularized solutions they lie between $1$ and $5\%$ on average and can go up to ca. $53\%$ in the extreme cases as one can see in the leftmost subplot for the first optical window. The run times of Algorithm \ref{reconstruction_algorithm} lie between $30$ and $50$ seconds in the average case and can rise up to $200$ seconds in the extreme cases. A typical sweep through all $48$ wavelengths needed ca. $30$ minutes in total and this value was very much the same for all three materials. For Ag the regularization procedure effectively reduced the relative errors such that they are in the range between $0.5$ and $2.2\%$ on average and are below $10\%$ in the extreme cases. Finally one can see in the last plot that the relative errors of the average of all $100$ regularized solutions are all below $0.4\%$.

For CsI the relative errors of the unregularized solutions are already quite small and lie between $0.03$ and $0.065\%$ on average and rise only up to $0.3\%$ in the extreme cases. The run times of Algorithm \ref{reconstruction_algorithm} are typically in the range from $25$ to $55$ seconds and are always below $95$ seconds. The regularization of the solutions brought only a small improvement of the results here such that the relative errors did not change much. They are still in the same range from $0.03$ and $0.065\%$ on average but only reach up to ca. $0.2\%$ now. The relative errors of the average of the $100$ regularized solutions are between $0.01$ and $0.055\%$.

Also for $\mathrm{H}_{2}$O the relative errors of the unregularized solutions are comparably small and are below $0.35\%$ on average and still below $1.3\%$ in the extreme cases. Especially the rightmost subplot for the last optical window shows the biggest relative errors, whereas for all the other optical windows the relative errors are below $0.03\%$ on average and below $0.15\%$ in the extreme cases. A similar behavior can be observed for the run times of Algorithm \ref{reconstruction_algorithm}. For the first four optical windows they are between $20$ and $45$ seconds on average and below $100$ seconds in the extreme cases, whereas for the last optical window they are between $30$ and $170$ seconds on average and can even rise up to $350$ seconds. For $\mathrm{H}_{2}$O the regularization procedure improves the relative errors only slightly for the first four optical windows and even increases them for the last optical window such that they can rise up to ca. $0.06\%$ on average and $1.4\%$ in the extreme cases. The relative errors of the average of the $100$ regularized solutions are virtually zero for the first four optical windows and below $0.55\%$ for the last optical window.

\section{Higher Noise Levels}

To see how our proposed reconstruction algortihm behaves for higher noise levels, we performed for each of the scatterer materials Ag, $\mathrm{H}_{2}$O and CsI two numerical studies with $10$ sweeps through all $48$ wavelenghts of the five optical windows with the same settings as in Section \ref{comparison_with_trunc_ind}. We computed original spectral extinctions
\begin{equation*}
\left(\boldsymbol{e}_{true}\right)_{i,j} := \pi r_{i}^{2}\sum_{n = 1}^{N_{tr}}q_{n}(m_{med}(l_{i}), m_{part}(l_{i}), r_{j}, l_{i}), \quad i = 1, ..., 48, \quad j = 1, ..., 3 
\end{equation*}
for all wavelengths  and added zero-mean Gaussian noise to it in order to obtain the simulated noisy spectral extinctions
\begin{equation*}
\left(\boldsymbol{e}\right)_{i,j} = \left(\boldsymbol{e}_{true}\right)_{i,j} + \delta_{i,j} \quad \text{ with } \; \delta_{i,j} \sim \mathcal{N}(0,(0.15 \cdot \left(\boldsymbol{e}_{true}\right)_{i,j})^{2}), \quad i = 1, ..., 48, \quad j = 1, ..., 3.
\end{equation*}
for the first study and
\begin{equation*}
\left(\boldsymbol{e}\right)_{i,j} = \left(\boldsymbol{e}_{true}\right)_{i,j} + \delta_{i,j} \quad \text{ with } \; \delta_{i,j} \sim \mathcal{N}(0,(0.3 \cdot \left(\boldsymbol{e}_{true}\right)_{i,j})^{2}), \quad i = 1, ..., 48, \quad j = 1, ..., 3.
\end{equation*}
for the second. The standard deviations were taken to be $15\%$ and $30\%$ respectively of the true spectral extinctions. We used a sample size of $N_{s} = 300$ to compute each means $\left(\boldsymbol{e}_{real}\right)_{i,j}$ of noisy spectral extinctions. 

For brevity we only present the relative errors of the average of the $10$ regularized solutions.

\subsection{Results for Ag}

\subsubsection{Standard Deviation of $15\%$}

\begin{figure}[h!]
\begingroup
\sbox0{\includegraphics[width =1.0\textwidth]{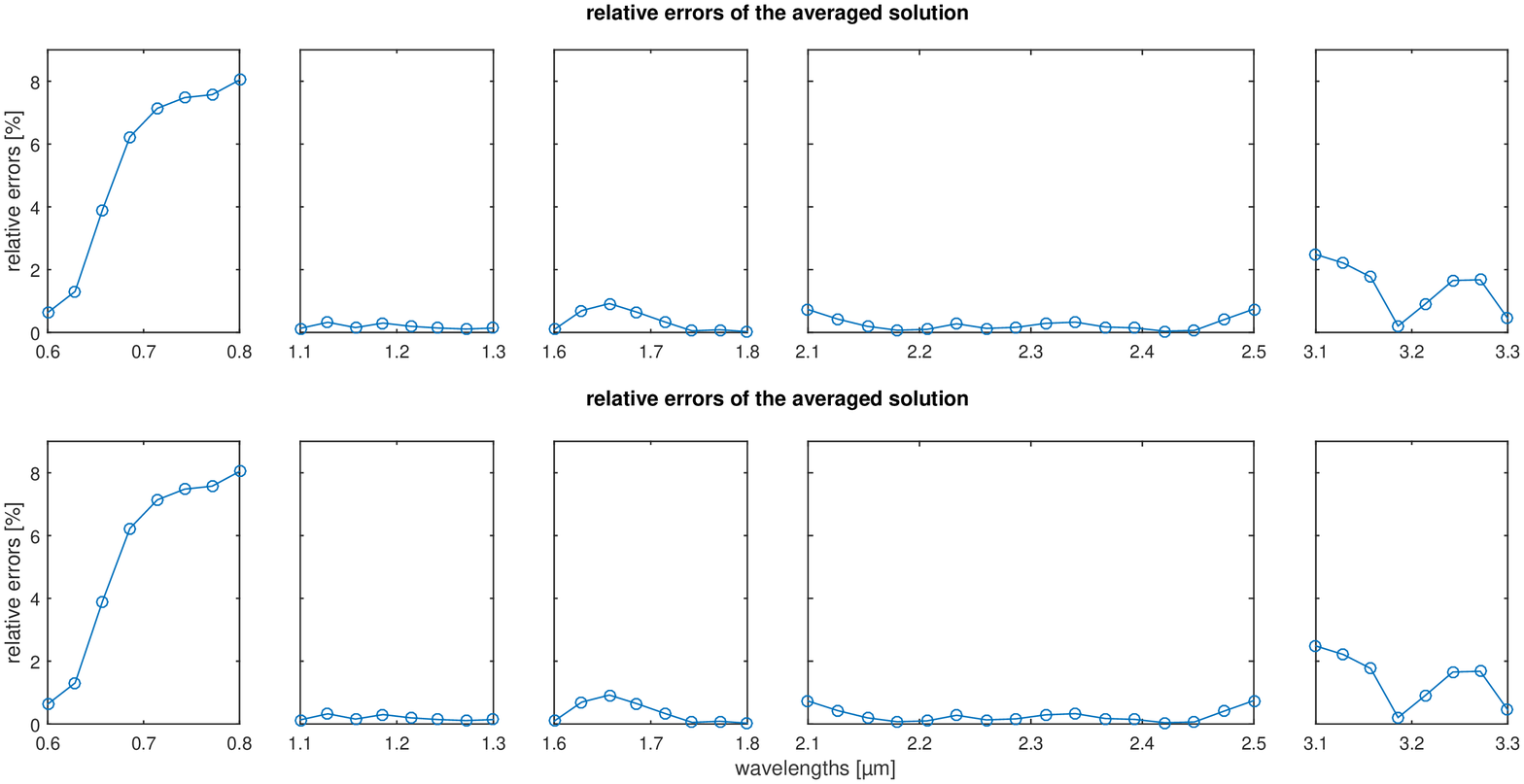}}
\includegraphics[clip,trim=0 0 0 {.5\wd0},width=1.0\textwidth]{ErrorsAveragedSolutionSilver_0dot15.eps}
\endgroup
\label{ErrRegAvgSilver_0.15}
\end{figure}

\newpage

\subsubsection{Standard Deviation of $30\%$}

\begin{figure}[h!]
\begingroup
\sbox0{\includegraphics[width =1.0\textwidth]{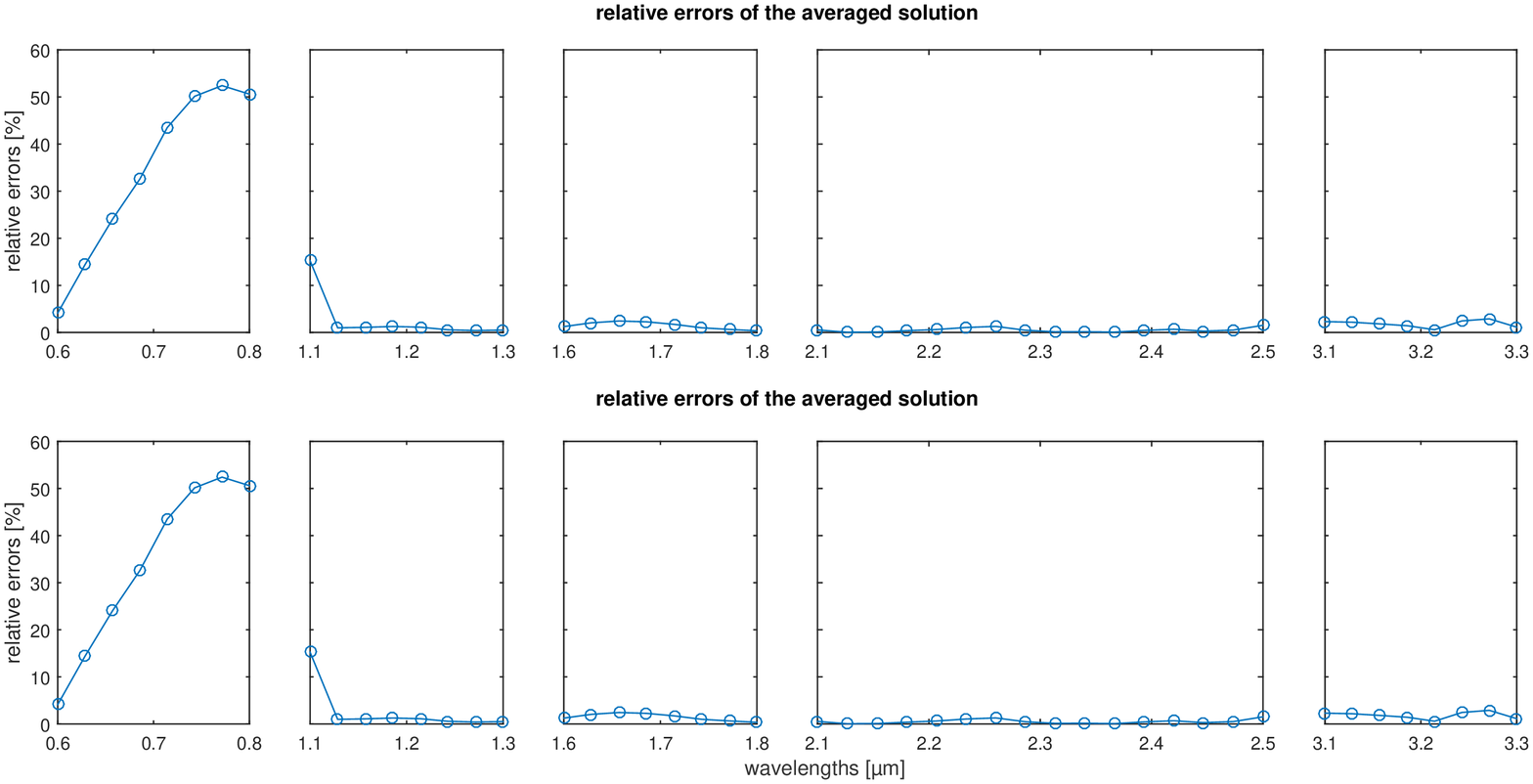}}
\includegraphics[clip,trim=0 0 0 {.5\wd0},width=1.0\textwidth]{ErrorsAveragedSolutionSilver_0dot3.eps}
\endgroup
\label{ErrRegAvgSilver_0.3}
\end{figure}

\subsection{Results for CsI}

\subsubsection{Standard Deviation of $15\%$}

\begin{figure}[h!]
\begingroup
\sbox0{\includegraphics[width =1.0\textwidth]{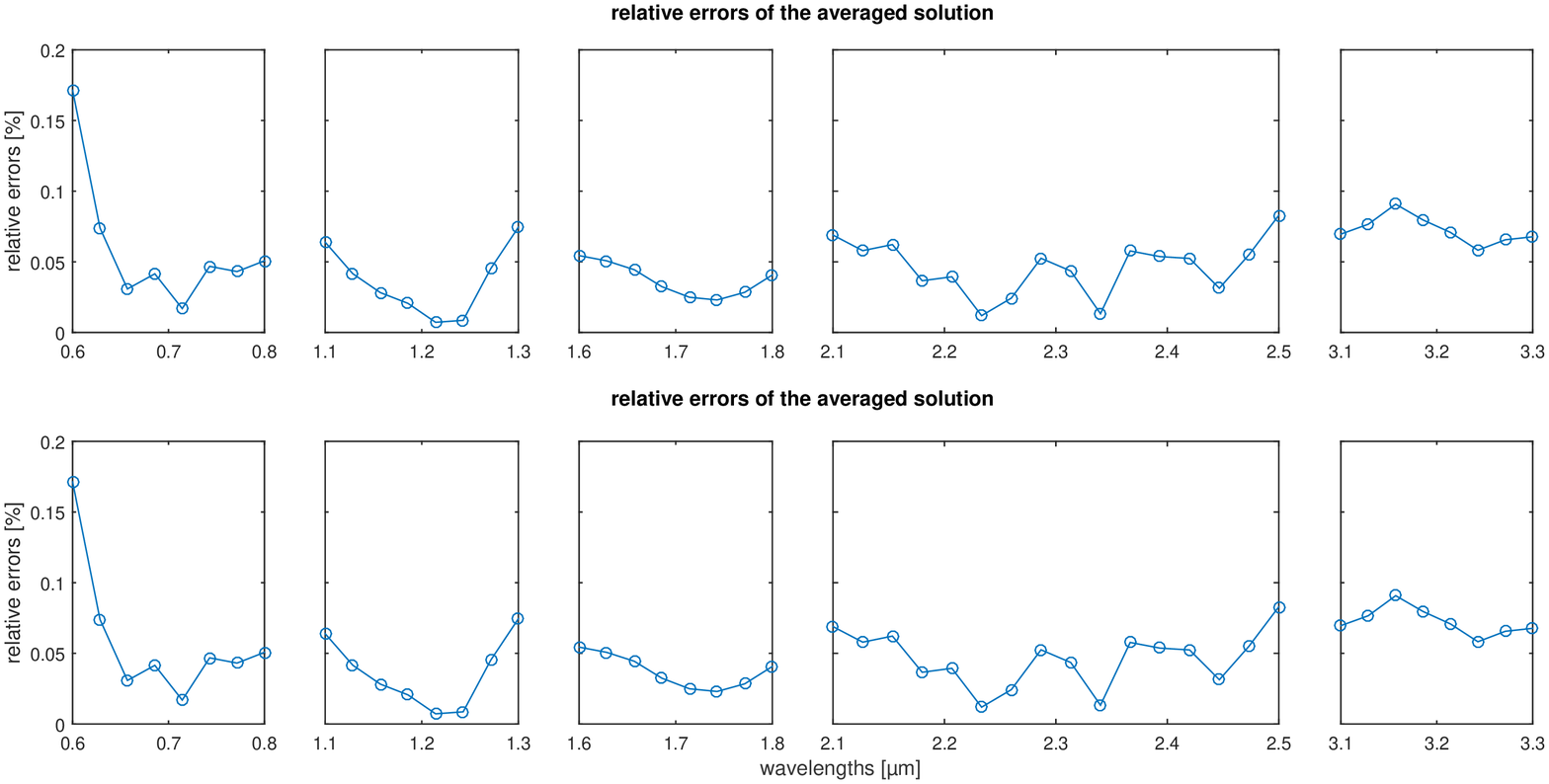}}
\includegraphics[clip,trim=0 0 0 {.5\wd0},width=1.0\textwidth]{ErrorsAveragedSolutionCsI_0dot15.eps}
\endgroup
\label{ErrRegAvgCsI_0.15}
\end{figure}

\subsubsection{Standard Deviation of $30\%$}

\begin{figure}[h!]
\begingroup
\sbox0{\includegraphics[width =1.0\textwidth]{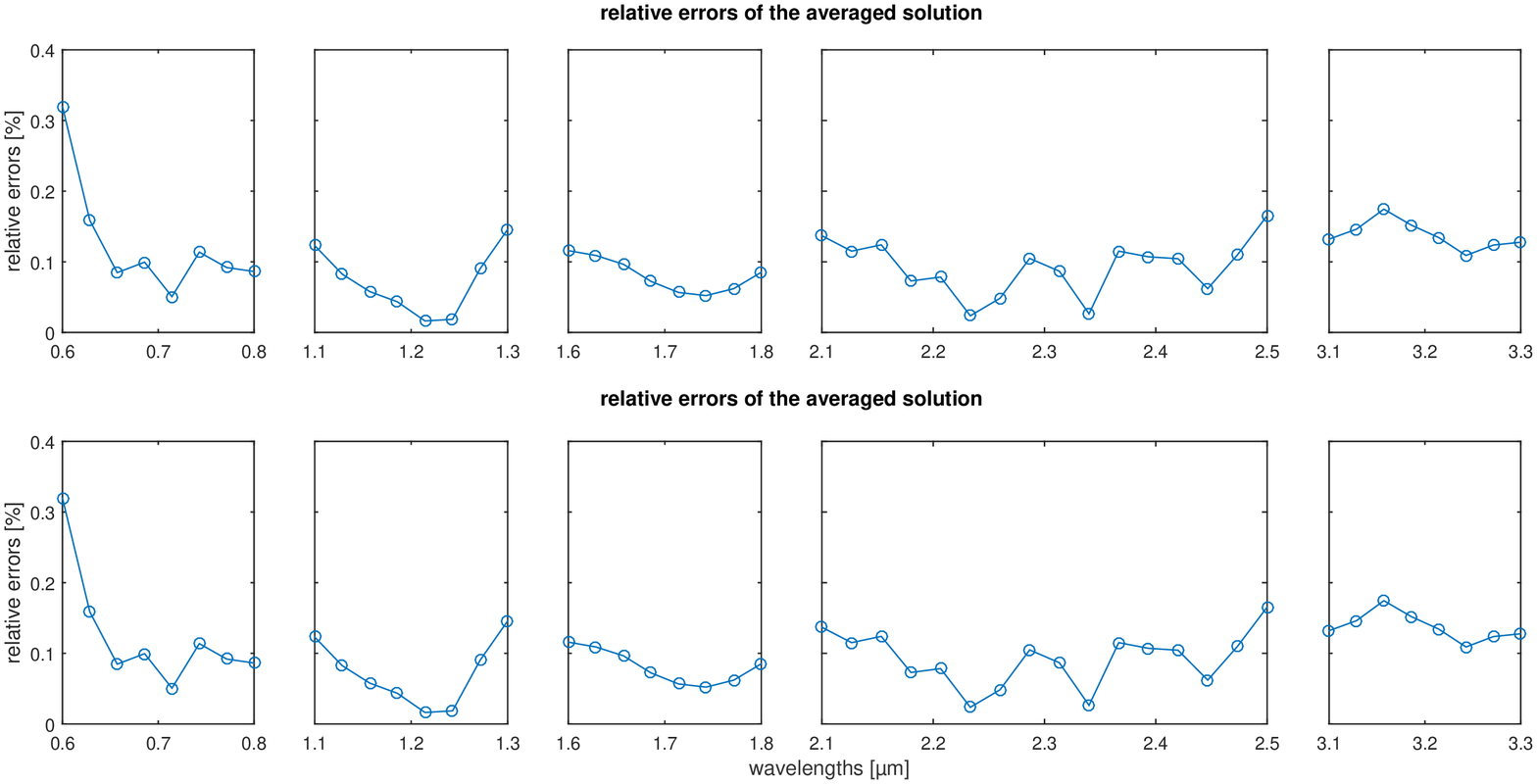}}
\includegraphics[clip,trim=0 0 0 {.5\wd0},width=1.0\textwidth]{ErrorsAveragedSolutionCsI_0dot3.eps}
\endgroup
\label{ErrRegAvgCsI_0.3}
\end{figure}

\subsection{Results for $\mathrm{H}_{2}$O}

\subsubsection{Standard Deviation of $15\%$}

\begin{figure}[h!]
\begingroup
\sbox0{\includegraphics[width =1.0\textwidth]{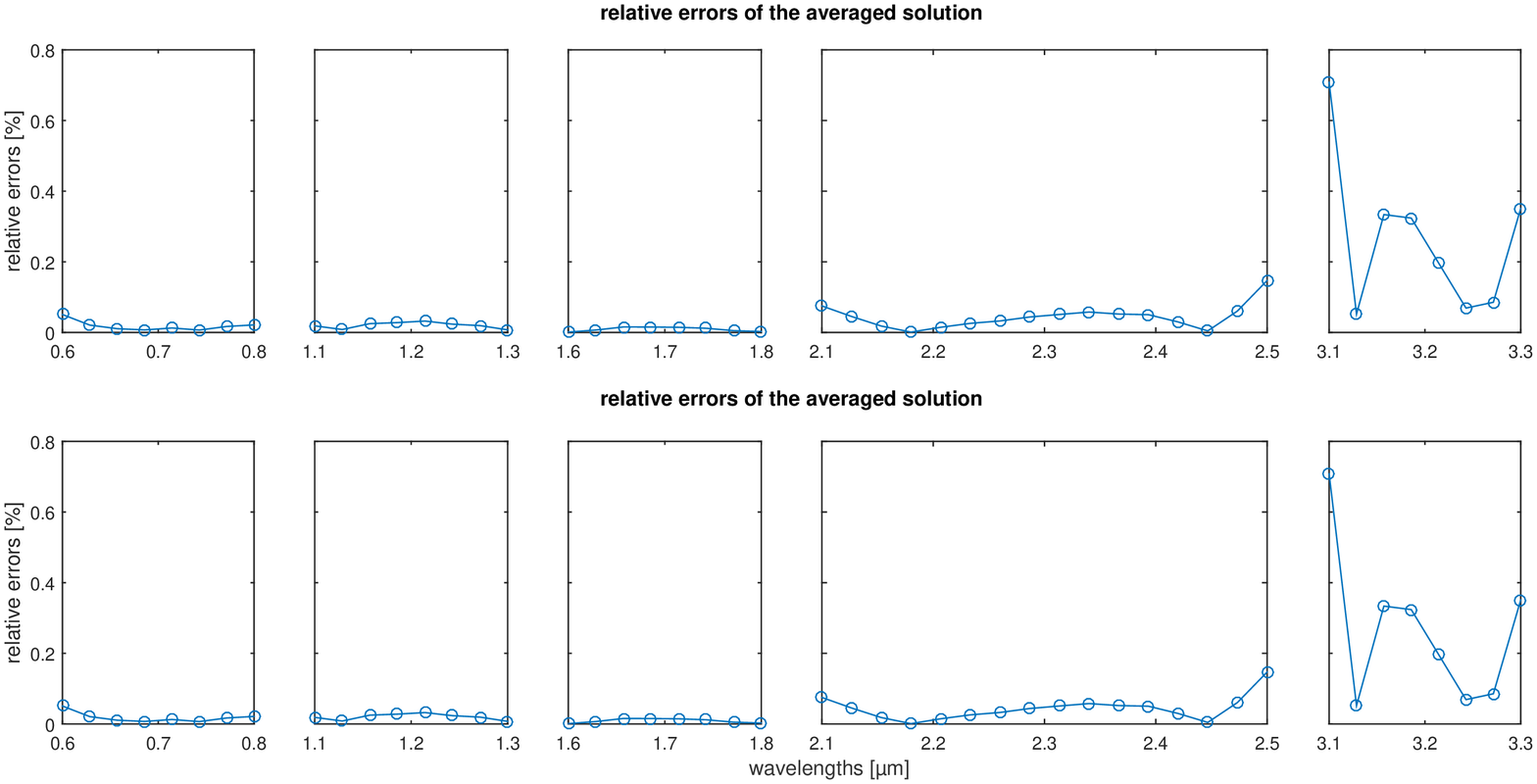}}
\includegraphics[clip,trim=0 0 0 {.5\wd0},width=1.0\textwidth]{ErrorsAveragedSolutionWater_0dot15.eps}
\endgroup
\label{ErrRegAvgWater_0.15}
\end{figure}

\newpage

\subsubsection{Standard Deviation of $30\%$}

\begin{figure}[h!]
\begingroup
\sbox0{\includegraphics[width =1.0\textwidth]{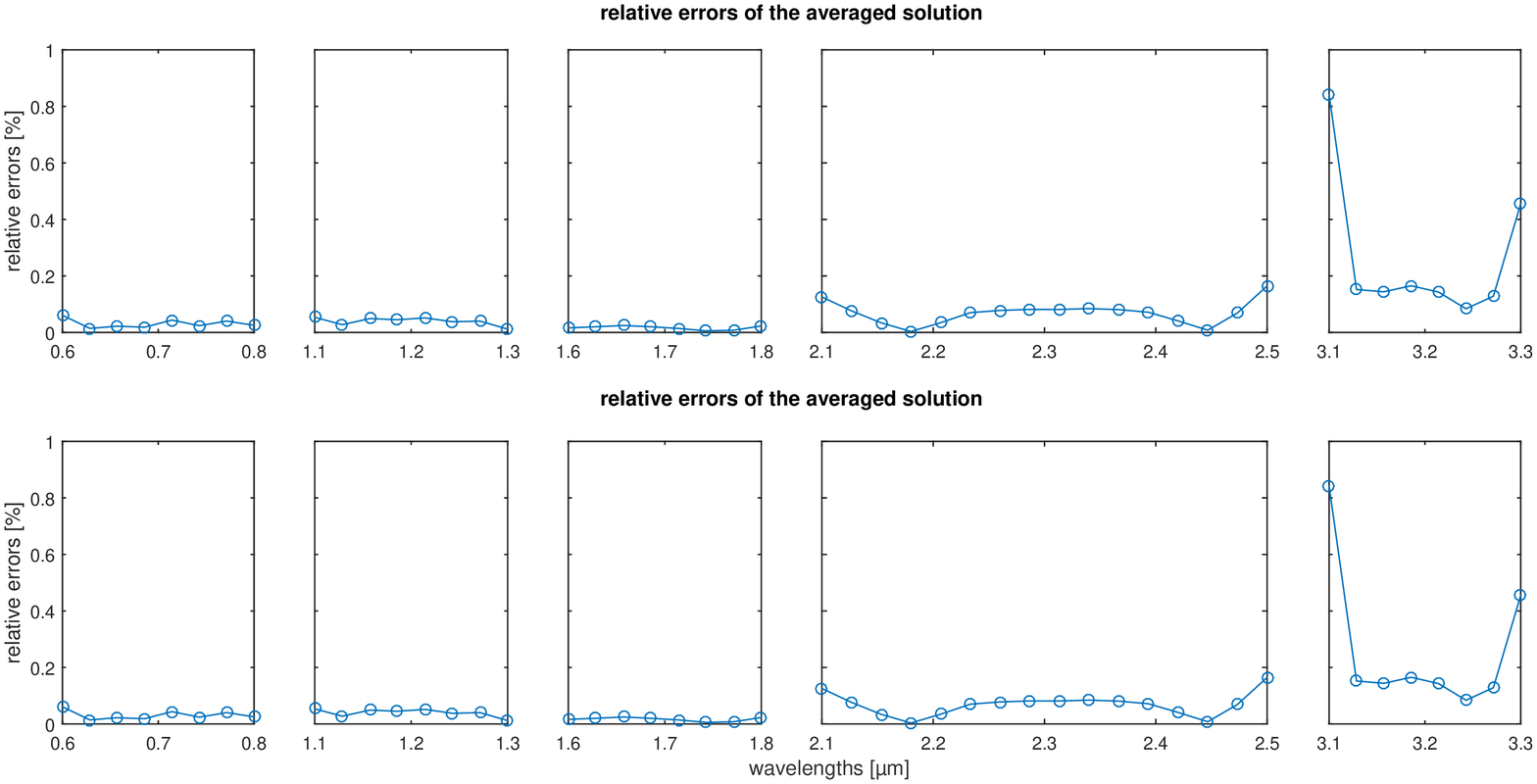}}
\includegraphics[clip,trim=0 0 0 {.5\wd0},width=1.0\textwidth]{ErrorsAveragedSolutionWater_0dot3.eps}
\endgroup
\label{ErrRegAvgWater_0.3}
\end{figure}

\subsection{Conclusion}

Whereas the relative errors for CsI and $\mathrm{H}_{2}$O are still below $1\%$, they can rise up to ca. $53\%$ for Ag. Therefore the reconstructed refractive indices for Ag under this noise level are most likely not of practical use. This shows that the FASP measurements of monodisperse aerosols must be sufficiently accurate in order to retrieve the scatterer refractive indices from them.

\section{Numerical Study}

We performed four numerical studies for two-component aerosols with log-normal, RRSB and Hedrih model size distributions as outlined in \cite{AK16}. The aerosol particles were assumed to be homogeneously internally mixed, such that only one effective refractive index was retrieved. One component of the simulated aerosols was $\mathrm{H}_{2}$O with volume fractions of $0$, $11$, $22$, $33$, $44$, $56$, $67$, $78$, $89$ and $100\%$. In the first two studies we simulated mixtures of $\mathrm{H}_{2}$O and CsI, where we used the original aerosol component refractive indices for the first study. For the second study we used the average of the $100$ regularized solutions from Section \ref{ReconRefrac}. We did the same for the third and fourth study, but here we simulated mixtures of $\mathrm{H}_{2}$O and Ag. In the third study we utilized the original aerosol component refractive indices and for the fourth the average of the $100$ regularized solutions from Section \ref{ReconRefrac}.

We applied the same reconstruction methods described in \cite{AK16} under the same settings, i.e. for each reconstruction we generated $300$ artificial noisy measurements for all $48$ wavelengths, where the measurement error was simulated as additive zero-mean Gaussian noise. For each wavelength, the standard deviations were taken as $5\%$ of the solutions of the forward problem. In \cite{AK16} three different regularization methods, namely Tikhonov, minimal first differences and Twomey regularization, were compared and their results turned out to be very similar. Therefore we only used Tikhonov regularization in the following. The results for the first study were directly adopted from \cite{AK16}.

For every inversion we computed the $L^{2}$-error of the obtained reconstruction relative to the original size distribution and measured the total run time needed for the inversion. The computations were performed on a notebook with a $2.27$ GHz CPU and $3.87$ GB accessible primary memory.

\newpage

\section{Results for Mixtures of $\mathrm{H}_{2}$O and CsI}

\subsubsection{Noise-free Refractive Indices}

\begin{table}[h!]
\centering
\begin{tabular}{||m{2.55cm}||m{2.5cm}|m{2.5cm}|m{2.5cm}||}
\hline
\multirow{4}{*}{\vspace{0.9cm} original $\mathrm{H}_{2}$O} & \multicolumn{3}{c||}{} \tabularnewline[-2.5ex]
& \multicolumn{3}{c||}{average $L^{2}$-errors (\%)} \tabularnewline
& \multicolumn{3}{c||}{} \tabularnewline[-2.5ex]
\cline{2-4}
& \centering Log-Normal & \centering RRSB & \centering Hedrih \tabularnewline
\vspace{-0.6cm} volume percent & \centering Distribution & \centering Distribution & \centering Distribution \tabularnewline
\hline
\hline
0 \% & \centering 33.5018 & \centering 38.4939 & \centering 17.2492 \tabularnewline
\hline
11 \% & \centering 30.0691 & \centering 31.7809 & \centering 17.5805 \tabularnewline
\hline
22 \% & \centering 28.8329 & \centering 30.4269 & \centering 17.3665 \tabularnewline
\hline
33 \% & \centering 24.9249 & \centering 27.7925 & \centering 15.0919 \tabularnewline
\hline
44 \% & \centering 23.4809 & \centering 24.2137 & \centering 16.9073 \tabularnewline
\hline
56 \% & \centering 22.6134 & \centering 21.5396 & \centering 15.5153 \tabularnewline
\hline
67 \% & \centering 20.4546 & \centering 19.5181 & \centering 14.5155 \tabularnewline
\hline
78 \% & \centering 18.9590 & \centering 16.9927 & \centering 16.7670 \tabularnewline
\hline
89 \% & \centering 18.5494 & \centering 14.4005 & \centering 13.3015 \tabularnewline
\hline
100 \% & \centering 17.9452 & \centering 12.0441 & \centering 11.3083 \tabularnewline
\hline
\end{tabular}
\end{table}

\begin{table}[h!]
\centering
\begin{tabular}{||m{2.55cm}||m{2.5cm}|m{2.5cm}|m{2.5cm}||}
\hline
\multirow{4}{*}{\vspace{0.9cm} original $\mathrm{H}_{2}$O} & \multicolumn{3}{c||}{} \tabularnewline[-2.5ex]
& \multicolumn{3}{c||}{average fraction deviation (\%)} \tabularnewline
& \multicolumn{3}{c||}{} \tabularnewline[-2.5ex]
\cline{2-4}
& \centering Log-Normal & \centering RRSB & \centering Hedrih \tabularnewline
\vspace{-0.6cm} volume percent & \centering Distribution & \centering Distribution & \centering Distribution \tabularnewline
\hline
\hline
0 \% & \centering 10.7150 & \centering 6.5800 & \centering 5.0800 \tabularnewline
\hline
11 \% & \centering 7.3700 & \centering 5.3000 & \centering 6.3050 \tabularnewline
\hline
22 \% & \centering 6.1750 & \centering 4.7100 & \centering 4.2600 \tabularnewline
\hline
33 \% & \centering 4.3700 & \centering 3.7300 & \centering 4.4350 \tabularnewline
\hline
44 \% & \centering 3.9550 & \centering 3.7300 & \centering 3.2100 \tabularnewline
\hline
56 \% & \centering 3.2000 & \centering 3.1650 & \centering 3.3750 \tabularnewline
\hline
67 \% & \centering 2.6050 & \centering 2.3650 & \centering 1.8700 \tabularnewline
\hline
78 \% & \centering 2.2350 & \centering 1.8050 & \centering 2.7750 \tabularnewline
\hline
89 \% & \centering 2.2000 & \centering 1.3750 & \centering 1.9150 \tabularnewline
\hline
100 \% & \centering 1.4150 & \centering 0.4650 & \centering 1.0050 \tabularnewline
\hline
\end{tabular}
\end{table}

\subsubsection{Noisy Refractive Indices}

\begin{table}[h!]
\centering
\begin{tabular}{||m{2.55cm}||m{2.5cm}|m{2.5cm}|m{2.5cm}||}
\hline
\multirow{4}{*}{\vspace{0.9cm} original $\mathrm{H}_{2}$O} & \multicolumn{3}{c||}{} \tabularnewline[-2.5ex]
& \multicolumn{3}{c||}{average $L^{2}$-errors (\%)} \tabularnewline
& \multicolumn{3}{c||}{} \tabularnewline[-2.5ex]
\cline{2-4}
& \centering Log-Normal & \centering RRSB & \centering Hedrih \tabularnewline
\vspace{-0.6cm} volume percent & \centering Distribution & \centering Distribution & \centering Distribution \tabularnewline
\hline
\hline
0 \% & \centering 32.8763 & \centering 42.3379 & \centering 18.1280 \tabularnewline
\hline
11 \% & \centering 31.3243 & \centering 31.8214 & \centering 17.1997 \tabularnewline
\hline
22 \% & \centering 27.6456 & \centering 32.6273 & \centering 17.8791 \tabularnewline
\hline
33 \% & \centering 24.5193 & \centering 27.1816 & \centering 15.4878 \tabularnewline
\hline
44 \% & \centering 23.4746 & \centering 24.2689 & \centering 16.0616 \tabularnewline
\hline
56 \% & \centering 22.8027 & \centering 21.3374 & \centering 15.1618 \tabularnewline
\hline
67 \% & \centering 19.6712 & \centering 17.9423 & \centering 14.4443 \tabularnewline
\hline
78 \% & \centering 18.4366 & \centering 16.5406 & \centering 17.5635 \tabularnewline
\hline
89 \% & \centering 18.4948 & \centering 14.1970 & \centering 13.8089 \tabularnewline
\hline
100 \% & \centering 18.4893 & \centering 11.2708 & \centering 12.4506 \tabularnewline
\hline
\end{tabular}
\end{table}

\newpage

\begin{table}[h!]
\centering
\begin{tabular}{||m{2.55cm}||m{2.5cm}|m{2.5cm}|m{2.5cm}||}
\hline
\multirow{4}{*}{\vspace{0.9cm} original $\mathrm{H}_{2}$O} & \multicolumn{3}{c||}{} \tabularnewline[-2.5ex]
& \multicolumn{3}{c||}{average fraction deviation (\%)} \tabularnewline
& \multicolumn{3}{c||}{} \tabularnewline[-2.5ex]
\cline{2-4}
& \centering Log-Normal & \centering RRSB & \centering Hedrih \tabularnewline
\vspace{-0.6cm} volume percent & \centering Distribution & \centering Distribution & \centering Distribution \tabularnewline
\hline
\hline
0 \% & \centering 10.9350 & \centering 7.6900 & \centering 4.8150 \tabularnewline
\hline
11 \% & \centering 6.6400 & \centering 4.8600 & \centering 6.5900 \tabularnewline
\hline
22 \% & \centering 5.6600 & \centering 4.8650 & \centering 5.0150 \tabularnewline
\hline
33 \% & \centering 4.1400 & \centering 3.9450 & \centering 4.1950 \tabularnewline
\hline
44 \% & \centering 3.8600 & \centering 3.4400 & \centering 3.4650 \tabularnewline
\hline
56 \% & \centering 3.7550 & \centering 2.8200 & \centering 3.3550 \tabularnewline
\hline
67 \% & \centering 2.3500 & \centering 2.3900 & \centering 1.9550 \tabularnewline
\hline
78 \% & \centering 2.0250 & \centering 1.9600 & \centering 2.8300 \tabularnewline
\hline
89 \% & \centering 2.1350 & \centering 1.2600 & \centering 1.9050 \tabularnewline
\hline
100 \% & \centering 1.6100 & \centering 0.3350 & \centering 1.0200 \tabularnewline
\hline
\end{tabular}
\end{table}

\subsection{Results for Mixtures of $\mathrm{H}_{2}$O and Ag}

\subsubsection{Noise-free Refractive Indices}

\begin{table}[h!]
\centering
\begin{tabular}{||m{2.55cm}||m{2.5cm}|m{2.5cm}|m{2.5cm}||}
\hline
\multirow{4}{*}{\vspace{0.9cm} original $\mathrm{H}_{2}$O} & \multicolumn{3}{c||}{} \tabularnewline[-2.5ex]
& \multicolumn{3}{c||}{average $L^{2}$-errors (\%)} \tabularnewline
& \multicolumn{3}{c||}{} \tabularnewline[-2.5ex]
\cline{2-4}
& \centering Log-Normal & \centering RRSB & \centering Hedrih \tabularnewline
\vspace{-0.6cm} volume percent & \centering Distribution & \centering Distribution & \centering Distribution \tabularnewline
\hline
\hline
0 \% & \centering 61.2258 & \centering 67.4898 & \centering 54.9872 \tabularnewline
\hline
11 \% & \centering 47.6091 & \centering 59.5032 & \centering 49.4157 \tabularnewline
\hline
22 \% & \centering 37.1819 & \centering 67.8529 & \centering 36.8106 \tabularnewline
\hline
33 \% & \centering 57.4058 & \centering 69.0461 & \centering 46.9500 \tabularnewline
\hline
44 \% & \centering 54.5266 & \centering 71.4585 & \centering 51.6323 \tabularnewline
\hline
56 \% & \centering 45.6020 & \centering 67.5288 & \centering 36.6557 \tabularnewline
\hline
67 \% & \centering 36.0602 & \centering 53.0986 & \centering 22.8959 \tabularnewline
\hline
78 \% & \centering 33.9722 & \centering 43.2794 & \centering 16.1026 \tabularnewline
\hline
89 \% & \centering 24.1114 & \centering 27.5773 & \centering 16.0464 \tabularnewline
\hline
100 \% & \centering 21.0693 & \centering 11.1298 & \centering 11.5123 \tabularnewline
\hline
\end{tabular}
\end{table}

\begin{table}[h!]
\centering
\begin{tabular}{||m{2.55cm}||m{2.5cm}|m{2.5cm}|m{2.5cm}||}
\hline
\multirow{4}{*}{\vspace{0.9cm} original $\mathrm{H}_{2}$O} & \multicolumn{3}{c||}{} \tabularnewline[-2.5ex]
& \multicolumn{3}{c||}{average fraction deviation (\%)} \tabularnewline
& \multicolumn{3}{c||}{} \tabularnewline[-2.5ex]
\cline{2-4}
& \centering Log-Normal & \centering RRSB & \centering Hedrih \tabularnewline
\vspace{-0.6cm} volume percent & \centering Distribution & \centering Distribution & \centering Distribution \tabularnewline
\hline
\hline
0 \% & \centering 0 & \centering 0 & \centering 0 \tabularnewline
\hline
11 \% & \centering 0.2650 & \centering 0.3500 & \centering 0.2950 \tabularnewline
\hline
22 \% & \centering 0.5350 & \centering 3.6650 & \centering 1.8050 \tabularnewline
\hline
33 \% & \centering 18.5000 & \centering 16.4750 & \centering 10.4350 \tabularnewline
\hline
44 \% & \centering 13.1200 & \centering 13.5400 & \centering 14.5050 \tabularnewline
\hline
56 \% & \centering 8.6300 & \centering 9.6450 & \centering 8.0700 \tabularnewline
\hline
67 \% & \centering 6.5950 & \centering 6.8300 & \centering 5.7750 \tabularnewline
\hline
78 \% & \centering 3.4650 & \centering 3.2900 & \centering 1.9750 \tabularnewline
\hline
89 \% & \centering 1.3150 & \centering 1.5450 & \centering 1.0600 \tabularnewline
\hline
100 \% & \centering 0.4750 & \centering 0.0600 & \centering 0.1600 \tabularnewline
\hline
\end{tabular}
\end{table}

\newpage

\subsubsection{Noisy Refractive Indices}

\begin{table}[h!]
\centering
\begin{tabular}{||m{2.55cm}||m{2.5cm}|m{2.5cm}|m{2.5cm}||}
\hline
\multirow{4}{*}{\vspace{0.9cm} original $\mathrm{H}_{2}$O} & \multicolumn{3}{c||}{} \tabularnewline[-2.5ex]
& \multicolumn{3}{c||}{average $L^{2}$-errors (\%)} \tabularnewline
& \multicolumn{3}{c||}{} \tabularnewline[-2.5ex]
\cline{2-4}
& \centering Log-Normal & \centering RRSB & \centering Hedrih \tabularnewline
\vspace{-0.6cm} volume percent & \centering Distribution & \centering Distribution & \centering Distribution \tabularnewline
\hline
\hline
0 \% & \centering 62.0576 & \centering 67.5179 & \centering 57.9628 \tabularnewline
\hline
11 \% & \centering 46.9321 & \centering 59.9884 & \centering 51.1218 \tabularnewline
\hline
22 \% & \centering 38.6799 & \centering 65.8397 & \centering 36.0770 \tabularnewline
\hline
33 \% & \centering 56.2737 & \centering 69.8686 & \centering 47.1554 \tabularnewline
\hline
44 \% & \centering 55.2457 & \centering 70.4010 & \centering 53.9243 \tabularnewline
\hline
56 \% & \centering 44.6657 & \centering 67.1041 & \centering 37.0255 \tabularnewline
\hline
67 \% & \centering 37.8969 & \centering 55.6161 & \centering 23.9009 \tabularnewline
\hline
78 \% & \centering 39.5606 & \centering 45.8995 & \centering 16.6339 \tabularnewline
\hline
89 \% & \centering 24.1072 & \centering 28.0702 & \centering 16.2226 \tabularnewline
\hline
100 \% & \centering 18.3892 & \centering 9.8479 & \centering 11.4723 \tabularnewline
\hline
\end{tabular}
\end{table}

\begin{table}[h!]
\centering
\begin{tabular}{||m{2.55cm}||m{2.5cm}|m{2.5cm}|m{2.5cm}||}
\hline
\multirow{4}{*}{\vspace{0.9cm} original $\mathrm{H}_{2}$O} & \multicolumn{3}{c||}{} \tabularnewline[-2.5ex]
& \multicolumn{3}{c||}{average fraction deviation (\%)} \tabularnewline
& \multicolumn{3}{c||}{} \tabularnewline[-2.5ex]
\cline{2-4}
& \centering Log-Normal & \centering RRSB & \centering Hedrih \tabularnewline
\vspace{-0.6cm} volume percent & \centering Distribution & \centering Distribution & \centering Distribution \tabularnewline
\hline
\hline
0 \% & \centering 0 & \centering 0 & \centering 0 \tabularnewline
\hline
11 \% & \centering 0.2700 & \centering 0.3950 & \centering 0.3200 \tabularnewline
\hline
22 \% & \centering 0.4450 & \centering 2.9900 & \centering 1.0750 \tabularnewline
\hline
33 \% & \centering 16.0650 & \centering 16.4650 & \centering 10.6400 \tabularnewline
\hline
44 \% & \centering 13.1950 & \centering 14.2150 & \centering 15.4550 \tabularnewline
\hline
56 \% & \centering 8.2350 & \centering 9.0350 & \centering 8.8350 \tabularnewline
\hline
67 \% & \centering 6.4650 & \centering 7.1100 & \centering 5.5450 \tabularnewline
\hline
78 \% & \centering 3.8000 & \centering 3.1700 & \centering 1.8200 \tabularnewline
\hline
89 \% & \centering 1.3250 & \centering 1.5600 & \centering 1.0900 \tabularnewline
\hline 
100 \% & \centering 0.3000 & \centering 0.0200 & \centering 0.0800 \tabularnewline
\hline
\end{tabular}
\end{table}

\subsection{Conclusion}

The resuts of the first and second study only differ by ca. $4\%$ at most and behave very similarly. The same is for the third and fourth study. These numerical results indicate that $100$ FASP measurement sweeps consisting of $300$ single measurements with an accuracy as in Section \ref{ReconRefrac} are sufficient to determine aerosol refractive indices in such a quality, that they are suitable for particle size distribution reconstructions for two-component homogeneously internally mixed aerosols using the FASP. The particle radii of the three monodisperse aerosols generated for the refractive indices retrieval need to be $0.1 \; \mu$m, $0.2 \; \mu$m and $0.3 \; \mu$m respectively. 

\section{Outlook}

It is of interest to investigate if the methods derived in this study can be extended to the case of core-plus-shell aerosols.

\section{Acknowledgement}

I thank Graham Alldredge, Ph.D. for proofreading the manuscript, useful recommendations and fruitful discussions on this topic. 

\begin{figure}[htbp]
\begin{minipage}{0.75\textwidth}
This work is sponsored by the German Federal Ministry of Education and Research (BMBF) under the contract number 02NUK022A.

Responsibility for the content of this report lies with the authors.
\end{minipage}
\hfill
\begin{minipage}{0.225\textwidth}
\includegraphics[width=\textwidth]{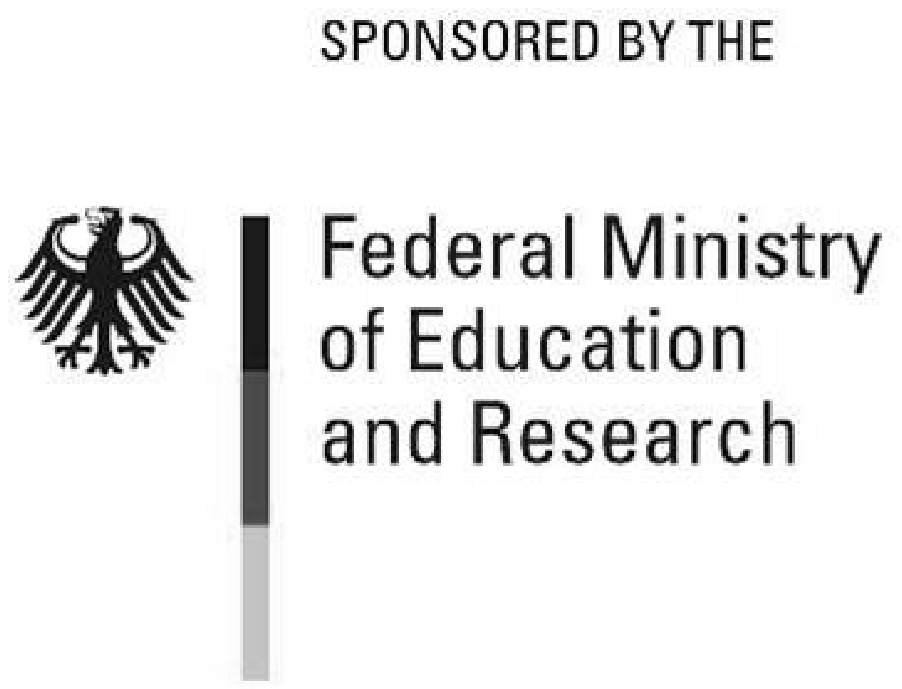}
\end{minipage}
\end{figure}

\bibliographystyle{ieeetr}
\bibliography{literature}

\end{document}